\documentclass[review,3p]{IEEEtran}

\usepackage{mathrsfs,latexsym,float,enumerate,amssymb,amsmath,bm,color,lineno}
\usepackage{amscd,verbatim,cuted,dsfont}
\usepackage{graphicx}
\usepackage[colorlinks=true,linkcolor=blue,citecolor=blue]{hyperref}
\usepackage{booktabs}
\usepackage{subfigure}
\usepackage{multirow}
\usepackage{bbding}
\usepackage{stfloats}
\usepackage{threeparttable}

\newtheorem{theorem}{Theorem}
\newtheorem{lemma}{Lemma}

\newtheorem{definition}{Definition}
\newtheorem{example}{Example}

\newtheorem{property}{Property}

\newtheorem{remark}{Remark}

\usepackage{pifont}

\begin{document}

\title{Strict Intuitionistic Fuzzy Distance/Similarity Measures Based on
Jensen-Shannon Divergence}

\author{Xinxing Wu, Zhiyi Zhu, Guanrong Chen, Tao Wang, Peide Liu 
\thanks{X. Wu is with 
(1) the School of Sciences, Southwest Petroleum University, Chengdu, Sichuan 610500, China;
(2) the Institute for Artificial Intelligence, Southwest Petroleum University, Chengdu, Sichuan 610500, China.
(3) the Zhuhai College of Jilin University, Zhuhai, Guangdong 519041, China.
e-mail: (wuxinxing5201314@163.com).}

\thanks{Z. Zhu is with the School of Sciences, Southwest Petroleum University, Chengdu, Sichuan 610500, China.
e-mail: (zhuzhiyi2019@163.com).}

\thanks{G. Chen is with the Department of Electrical Engineering, City University of
Hong Kong, Hong Kong SAR, China.
e-mail: (eegchen@cityu.edu.hk).}

\thanks{T. Wang is with the School of Sciences, Southwest Petroleum University, Chengdu, Sichuan 610500, China.
e-mail: (wt1994math@163.com).}

\thanks{P. Liu is with the School of Management Science and Engineering, Shandong University
of Finance and Economics, Jinan, Shandong 250014, China.
e-mail: (peide.liu@gmail.com).}

\thanks{All correspondences should be addressed to X. Wu.}
\thanks{This work was supported by the National Natural Science Foundation of China
(Nos.~11601449 and 71771140) and the Key Natural Science Foundation of Universities in
Guangdong Province (No. 2019KZDXM027).}
}

\markboth{}%
\markboth{\large This work has been submitted to the IEEE for possible publication. Copyright may be transferred without notice,
after which this version may no longer be accessible.}

\maketitle

\begin{abstract}
Being a pair of dual concepts, the normalized distance and similarity measures are very important
tools for decision-making and pattern recognition under intuitionistic fuzzy sets framework. To
be more effective for decision-making and pattern recognition applications, a good normalized distance
measure should ensure that its dual similarity measure satisfies the axiomatic definition.
In this paper, we first construct some examples to illustrate that the dual similarity measures of two nonlinear
distance measures introduced in [A distance measure for intuitionistic fuzzy sets and its
application to pattern classification problems, \emph{IEEE Trans. Syst., Man, Cybern., Syst.},
vol.~51, no.~6, pp. 3980--3992, 2021] and [Intuitionistic fuzzy sets: spherical representation
and distances, \emph{Int. J. Intell. Syst.}, vol.~24, no.~4, pp. 399--420, 2009]
do not meet the axiomatic definition of intuitionistic fuzzy similarity measure. We show that
(1) they cannot effectively distinguish some intuitionistic fuzzy values (IFVs) with obvious size
relationship; (2) except for the endpoints, there exist infinitely many pairs of IFVs,
where the maximum distance 1 can be achieved under these two distances; leading to
counter-intuitive results. To overcome these drawbacks, we introduce the concepts of
strict intuitionistic fuzzy distance measure (SIFDisM) and strict intuitionistic fuzzy similarity
measure (SIFSimM), and propose an improved intuitionistic fuzzy distance measure based on
Jensen-Shannon divergence. We prove that (1) it is a SIFDisM; (2) its dual similarity measure
is a SIFSimM; (3) its induced entropy is an intuitionistic fuzzy entropy. Comparative analysis
and numerical examples demonstrate that our proposed distance measure is completely superior to
the existing ones.
\end{abstract}


\IEEEpeerreviewmaketitle

\section{Introduction}
\IEEEPARstart{T}{o} deal with the ubiquitous uncertainty and fuzziness in real life, Zadeh
(1965)~\cite{Za1965} presented the fuzzy set (FS) theory by applying membership degrees to
measure the importance of a fuzzy element. Zadeh's FS generalized the concept of crisp sets,
which is described by a characteristic function that can take any value in the interval
$[0, 1]$. However, duo to the limitation of a membership function that only indicates two
(supporting and opposing) states of fuzziness, the fuzzy sets cannot express the neutral
state of ``this and also that". As a remedy to this, Atanassov \cite{Ata1986} generalized
Zadeh's fuzzy set by introducing the concept of intuitionistic fuzzy sets (IFSs) (see
also~\cite{Ata1999}), which was characterized by a membership function and a non-membership
function simultaneously, satisfying the condition that the sum of the membership
degree and the non-membership degree at every point is less than or equal to $1$.
Moreover, in dealing with many practical problems, it is not appropriate for experts
to precisely express their decisions with crisp numbers due to
the complexity and uncertainty of available information. Based on this observation, Atanassov and
Gargov~\cite{AG1989} further extended IFSs to interval-valued intuitionistic fuzzy sets (IVIFSs)
(see also \cite{Ata2020}), replacing the membership degree and the non-membership degree by some
closed intervals in $[0, 1]$.

To assess the differences of IFSs, the normalized IF distance measure (IFDisM) and
the IF similarity measure (IFSimM), being a pair of dual concepts (see \cite{DGM2017}), are important tools
for decision-making (\cite{Xu2007a,CCL2016a,ZZ2021,JLZJHZY2022}),
pattern recognition~(\cite{LC2002,LS2003,HYHL2012,PHK2013,CCL2016,Ngu2016,JJLY2019,Xiao2021}),
medical diagnosis~(\cite{LZ2018,BZLW2019,JLZJHZY2022}), clustering analysis (\cite{LSM2018,BZLW2019}),
image processing (\cite{LLC2014,JLZJHZY2022}), and data mining under intuitionistic fuzzy sets framework.
Motivated by the similarity measure for Zadeh's fuzzy sets, Li and Cheng~\cite{LC2002}
introduced the concept of similarity measure (SimM) for IFSs, which was improved by Mitchell~\cite{Mit2003},
and applied to pattern recognition problems. Recently, Szmidt~\cite{Sz2014}
presented an overview on IFDisMs and IFSimMs. The research on IFDisMs and IFSimMs
focuses on two aspects, one is the two-dimensional (2-D) representation of IFSs, which only considers
the membership and nonmembership degrees, and the other is the three-dimensional (3-D) representation
of IFSs, which simultaneously considers the membership, nonmembership, and indeterminacy degrees.
However, noticeably, because the indeterminacy degree is uniquely determined by the membership
and nonmembership degrees, the space of all IFVs is essentially a 2-D topological structure. This
indicates a serious problem that many existing distance and similarity measures via 3-D representation
of IFSs, including Euclidean similarity measure in~\cite{Sz2014} and Minkowski similarity measure
in (\cite{Xu2007a,Li2014}), do not meet the axiomatic definition of IFSimMs (see \cite[Section~3]{WCZ2022}).

For 2-D IFDisMs and IFSimMs, Szmidt and Kacprzyk~\cite{SK2000} presented the normalized Hamming distance
and normalized Euclidean distance for IFSs. Grzegorzewski~\cite{Gr2004} proposed an IFDisM
by using the Hausdorff metric for closed intervals (also see~\cite{HY2004}). Wang and Xin~\cite{WX2005}
introduced the axiomatic definition of distance measure as a dual concept of the similarity measure for
IFSs and constructed some new 2-D IFDisMs by combining the 2-D Hamming IFDisM~\cite{SK2000} and the
2-D Hausdorff IFDisM~\cite{Gr2004}. Further, based on IFDisMs introduced by
Wang and Xin~\cite{WX2005}, Xu and Chen~\cite{XC2008} presented some new IFSimMs.

For 3-D IFDisMs and IFSimMs, Xu and Chen~\cite{XC2008} suggested some IFSimMs based on the idea of
TOPSIS of Hwang and Yoon~\cite{HY1981}. Wu et al.~\cite{WCZ2022} pointed out that the IFSimM in
\cite[Eq.~14]{XC2008} based on the normalized Euclidean distance does not satisfy the axiomatic
definition of IFSimMs. Motivated by the divergence measure in information theory, Chen et al.~\cite{CCL2016}
introduced a novel IFSimM using centroid points of transformed right-angled TrFNs. Yang and Chiclana
\cite{YC2009} constructed a nonlinear spherical distance measure for IFSs by using the `$\arccos$' function.
Hung and Yang~\cite{HY2008} constructed a $J_{\gamma}$-divergence ($\gamma>0$) for IFSs and proved
that it satisfies the axiomatic definition of IFDisM of Wang and Xin~\cite{WX2005}
for $\gamma\in [1, 2]$. Joshi and Kumar~\cite{JK2018} obtained a dissimilarity Jensen-Shannon
divergence measure for IFSs. In fact, this dissimilarity measure is equivalent to the $J_{1}$-divergence
in \cite{HY2008}. Recently, Xiao~\cite{Xiao2021} proposed a new IFDisM $\bm{d}_{\chi}$ based
on the Jensen-Shannon divergence and showed that this distance measure
is superior to those in \cite{SK2000,Gr2004,WX2005,HPK2012,YC2012,SMLZC2018,SWQH2019}.
By direct calculation, it is easy to see that $\ln 2 \cdot \bm{d}_{\chi}=\sqrt{J_{1}}$,
i.e., Xiao's distance measure $\bm{d}_{\chi}$ is a special case of Hung and Yang's
$J_{\gamma}$-divergence (see Section~\ref{Sec-III}). For more results on the IFDisMs
and IFSimMs, see \cite{Sz2014,Li2014,XC2012}. It should be noted that because the
axiomatic definition (S4) of IFSimMs requires the admissibility with
Atanassov's partial order `$\subset$', and that Atanassov's partial order only indicates the
size relationship for membership degrees and non-membership degrees, implying that many existing
3-D IFSimMs may violate the axiomatic definition (S4) of IFSimMs.

We will provide Examples~\ref{Exm-Wu-a}--\ref{Exm-Wu-2-Sec2} below to show that Xiao's distance measure
\cite{Xiao2021}, Hung and Yang's $J_{1}$-divergence \cite{HY2008}, and Joshi and Kumar's
 dissimilarity Jensen-Shannon divergence measure in ~\cite{JK2018} have the follows
three drawbacks: (1) its dual similarity measure does not satisfy the
axiomatic definition (S4) of Definition~\ref{Def-Li-Cheng-IFVs}; (2) it cannot effectively
distinguish some IFVs/IFSs with obvious size relationship; (3) Except for the endpoints
$\langle 1, 0\rangle$ and $\langle 0, 1 \rangle$, there exist infinitely many pairs of
IFVs, where the maximum distance $1$ can be achieved under these distances.
Meanwhile, observing from Examples~\ref{Exm-8-YC} and \ref{Exm-Wu-2-1-Sec5} below,
Yang and Chiclana's spherical distance measure in \cite{YC2009} has the same drawbacks.
To distinguish IFVs more effectively and overcome the above three drawbacks for
IFDisMs/IFSimMs, we introduce the concepts of strict intuitionistic
fuzzy similarity measure (SIFSimM) and strict intuitionistic fuzzy distance measure
(SIFDisM) and construct a novel IFDisM based on Jensen-Shannon divergence measure.
We demonstrate that its dual similarity measure is a SIFSimM and its induced entropy measure
meets the axiomatic definition of IF entropy. Moreover, we present some comparative analysis
to illustrate that our proposed distance measure is completely superior to the existing
IFDisMs; in particular, it is much better than Xiao's distance measure in
\cite{Xiao2021}, 
Hung and Yang's $J_{\alpha}$-divergence in~\cite{HY2008}, Joshi and Kumar's
dissimilarity divergence in ~\cite{JK2018}, and Yang and Chiclana's spherical
distance in \cite{YC2009}. Finally, to demonstrate the effectiveness  of
our proposed IFSimM, we apply it to a practical pattern recognition problem.


\section{Preliminaries}

\subsection{Intuitionistic fuzzy set (IFS)}

\begin{definition}[{\textrm{\protect\cite[Definition~1.1]{Ata1999}}}]
Let $X$ be the universe of discourse (UOD). An \textit{intuitionistic fuzzy set}~(IFS)
$I$ in $X$ is defined as
\begin{equation}
{\small I=\left\{\frac{\langle \mu_{_{I}}(x), \nu_{_{I}}(x)\rangle}{x}
~\Big|~ x\in X\right\},}
\end{equation}
where the functions
$\mu_{_{I}}: X \longrightarrow [0,1]$
and
$\nu_{_{I}}: X \longrightarrow [0,1]$
define the \textit{degree of membership} and the
\textit{degree of non-membership} of the element $x \in X$ to the set $I$,
respectively, and for every $x\in X$,
$\mu_{_{I}}(x)+\nu_{_{I}}(x)\leq 1.$
\end{definition}

Let $\mathrm{IFS}(X)$ denote the set of all IFSs in $X$.
For $I\in \mathrm{IFS}(X)$, the \textit{indeterminacy degree} $\pi_{_{I}}
(x)$ of an element $x$ belonging to $I$ is defined by $\pi_{_I}(x)=1-\mu_{_I}(x)-\nu_{_I}(x)$.
Clearly, each subset $A$ of $X$ can be expressed as an IFS, {\small $A=\left\{
\frac{\langle \bm{1}_{A}(x), \bm{1}_{X\backslash A}(x)\rangle}{x} \mid x\in X\right\}$},
which is also called a \textit{crisp set}.
In \cite{Xu2007,XC2012}, the pair $\langle\mu_{_I}(x), \nu_{_I}(x)\rangle$ is
called an \textit{ intuitionistic fuzzy value} (IFV) or an \textit{intuitionistic fuzzy number} (IFN).
For convenience, use $\alpha=\langle \mu_{\alpha}, \nu_{\alpha}\rangle$ to represent an IFV $\alpha$,
which satisfies $\mu_{\alpha}\in [0, 1]$, $\nu_{\alpha}\in [0, 1]$, and $0\leq \mu_{\alpha}
+\nu_{\alpha}\leq 1$. Let $\Theta$ denote the set of all IFVs, i.e.,
$\Theta=\{\langle \mu, \nu \rangle\in [0, 1]^{2} \mid \mu+\nu \leq 1\}$.
For $\alpha=\langle \mu_{\alpha}, \nu_{\alpha}\rangle\in \Theta$, the
\textit{complement} $\alpha^{\complement}$
of $\alpha$ is $\alpha^{\complement}=\langle \nu_{\alpha}, \mu_{\alpha}\rangle$.

Atanassov's order `$\subset$' \cite{Ata1999}, defined by the property that $\alpha\subset \beta$ if and only if
$\alpha\cap \beta=\alpha$, is a partial order on $\Theta$. Clearly, $\alpha\subset \beta$
if and only if $\mu_{\alpha}\leq \mu_{\beta}$ and $\nu_{\alpha}\geq \nu_{\beta}$. The order `$\subsetneqq$' on
$\Theta$ is defined by the property $\alpha\subsetneqq \beta$ if and only if $\alpha\subset \beta$ and $\alpha\neq \beta$.

\subsection{Similarity/Distance measures for IFSs}
\begin{definition}[{\textrm{\protect\cite{XC2012,Li2014}}}]
\label{Def-Li-Cheng-IFVs}
A mapping $\mathbf{S}: \Theta\times \Theta\longrightarrow [0, 1]$ is
called an \textit{intuitionistic fuzzy similarity measure} (IFSimM) on $\Theta$
if it satisfies the following conditions:
for any $\alpha_1$, $\alpha_2$, $\alpha_3\in \Theta$,
\begin{enumerate}
  \item[(S1)] $0\leq \mathbf{S}(\alpha_1, \alpha_2)\leq 1$.
  \item[(S2)] $\mathbf{S}(\alpha_1, \alpha_2)=1$ if and only if $\alpha_1=\alpha_2$.
  \item[(S3)] $\mathbf{S}(\alpha_1, \alpha_2)=\mathbf{S}(\alpha_2, \alpha_1)$.
  \item[(S4)] If $\alpha_1\subset \alpha_2\subset \alpha_3$, then
  $\mathbf{S}(\alpha_1, \alpha_3)\leq \mathbf{S}(\alpha_1, \alpha_2)$
  and $\mathbf{S}(\alpha_1, \alpha_3) \leq \mathbf{S}(\alpha_2, \alpha_3)$.
\end{enumerate}
\end{definition}

\begin{definition}[{\textrm{\protect\cite{XC2012,Li2014}}}]
\label{Def-Li-Cheng-IFSs}
Let $X$ be a UOD. A mapping $\mathbf{S}: \mathrm{IFS}(X)\times \mathrm{IFS}(X)
\longrightarrow [0, 1]$ is called an
\textit{intuitionistic fuzzy similarity measure} (IFSimM)
on $\mathrm{IFS}(X)$ if it satisfies the following conditions:
for any $I_1$, $I_2$, $I_3\in \mathrm{IFS}(X)$,
\begin{enumerate}
  \item[(S1)] $0\leq \mathbf{S}(I_1, I_2)\leq 1$.
  \item[(S2)] $\mathbf{S}(I_1, I_2)=1$ if and only if $I_1=I_2$.
  \item[(S3)] $\mathbf{S}(I_1, I_2)=\mathbf{S}(I_2, I_1)$.
  \item[(S4)] If $I_1\subset I_2\subset I_3$, then $\mathbf{S}(I_1, I_3)\leq \mathbf{S}(I_1, I_2)$
  and $\mathbf{S}(I_1, I_3)$ $\leq \mathbf{S}(I_2, I_3)$.
\end{enumerate}
\end{definition}

To effectively compare and distinguish IFVs and IFSs, we introduce the
concept of strict intuitionistic fuzzy similarity/distance measure as follows.

\begin{definition}
A mapping $\mathbf{S}: \Theta\times \Theta\longrightarrow [0, 1]$ is
called a \textit{strict intuitionistic fuzzy similarity measure} (SIFSimM) on $\Theta$
if, for any $\alpha_1$, $\alpha_2$, $\alpha_3\in \Theta$,
it satisfies (Sl)--(S3) in Definition~\ref{Def-Li-Cheng-IFVs} and
(S4$^{\prime}$) and (S5) described by
\begin{enumerate}
\item[(S4$^{\prime}$)] (Strict distinctiveness) If $\alpha_1\subsetneqq \alpha_2\subsetneqq \alpha_3$,
then $\mathbf{S}(\alpha_1, \alpha_3)< \mathbf{S}(\alpha_1, \alpha_2)$
  and $\mathbf{S}(\alpha_1, \alpha_3) < \mathbf{S}(\alpha_2, \alpha_3)$.

\item[(S5)] (Extreme dissimilarity on endpoints) $\mathbf{S}(\alpha_1, \alpha_2)=0$ if and only if
($\alpha_1=\langle 0, 1\rangle$ and $\alpha_2=\langle 1, 0 \rangle$) or
($\alpha_1=\langle 1, 0\rangle$ and $\alpha_2=\langle 0, 1 \rangle$).
\end{enumerate}
\end{definition}

Property (S4$^{\prime}$) indicates that the similarity measure $\mathbf{S}$
can strictly distinguish every pair of different IFVs under Atanassov-order `$\subset$'.
Property (S5) indicates that it is extremely unsimilar (similarity measure is zero)
for a pair of IFVs depending only on two endpoints.

\begin{definition}
Let $X$ be a UOD. A mapping $\mathbf{S}: \mathrm{IFS}(X)\times \mathrm{IFS}(X)
\longrightarrow [0, 1]$ is called a
\textit{strict intuitionistic fuzzy similarity measure} (SIFSimM) on $\mathrm{IFS}(X)$
if, for any $I_1$, $I_2$, $I_3\in \mathrm{IFS}(X)$,
it satisfies (Sl)--(S3) in Definition~\ref{Def-Li-Cheng-IFSs} and
(S4$^{\prime}$) and (S5) described by
\begin{enumerate}
\item[(S4$^{\prime}$)] If $I_1\subsetneqq I_2\subsetneqq I_3$,
then $\mathbf{S}(I_1, I_3)< \mathbf{S}(I_1, I_2)$
  and $\mathbf{S}(I_1, I_3)$ $< \mathbf{S}(I_2, I_3)$.

\item[(S5)] $\mathbf{S}(I_1, I_2)=0$ if and only if, for any $x\in X$,
($I_1(x)=\langle 0, 1\rangle$ and $I_2(x)=\langle 1, 0 \rangle$) or
($I_1(x)=\langle 1, 0\rangle$ and $I_2(x)=\langle 0, 1 \rangle$).
\end{enumerate}
\end{definition}

\begin{remark}
Property (S5) can be equivalently expressed as that $\mathbf{S}(I_1, I_2)=0$ if and only if
$I_1$ is a crisp set and $I_{1}=I_{2}^{\complement}$.
\end{remark}

Dually, a mapping $d: \mathrm{IFS}(X)\times \mathrm{IFS}(X)
\longrightarrow [0, 1]$ is called a \textit{strict intuitionistic fuzzy
distance measure} (SIFDisM) on $\mathrm{IFS}(X)$ if it satisfies the following conditions:
\begin{enumerate}
  \item[(D1)] The mapping $d$ is a distance measure on $\mathrm{IFS}(X)$
  (see \cite[Definition 3.1]{Sz2014});
  \item[(D2)] The mapping $\mathbf{S}(\alpha, \beta)=1-d(\alpha, \beta)$
  is a SIFSimM on $\mathrm{IFS}(X)$.
\end{enumerate}

\subsection{Entropy measure for IFSs}

Entropy is an important information measure. Szmidt and Kacprzyk
\cite{SK2001} gave the axiomatic definition of entropy measure
for IFSs as follows:

\begin{definition}[{\textrm{\protect\cite{SK2001}}}]
A mapping $E: \Theta\longrightarrow [0, 1]$ is called an \textit{intuitionistic
fuzzy entropy measure} (IFEM) on $\Theta$ if it satisfies the following conditions:
for any $\alpha$, $\beta \in \Theta$,
\begin{enumerate}
  \item[(E1)] $E(\alpha)=0$ if and only if $\alpha=\langle 1, 0\rangle$ or $\alpha=\langle 0, 1\rangle$.
  \item[(E2)] $E(\alpha)=1$ if and only if $\mu_{\alpha}=\nu_{\alpha}$.
  \item[(E3)] $E(\alpha)=E(\alpha^{\complement})$.
  \item[(E4)] $E(\alpha)\leq E(\beta)$ whenever it holds either
  $\mu_{\alpha}\leq \mu_{\beta} \leq \nu_{\beta} \leq \nu_{\alpha}$ or
  $\mu_{\alpha}\geq \mu_{\beta} \geq \nu_{\beta} \geq \nu_{\alpha}$.
\end{enumerate}
\end{definition}

\begin{definition}[{\textrm{\protect\cite{SK2001}}}]
Let $X$ be a UOD. A mapping $E: \mathrm{IFS}(X)\longrightarrow [0, 1]$ is
called an \textit{intuitionistic fuzzy entropy measure} (IFEM) on $\mathrm{IFS}(X)$ if
it satisfies the following conditions: for any $I_1$, $I_2 \in \mathrm{IFS}(X)$,
\begin{enumerate}
  \item[(E1)] $E(I_1)=0$ if and only if $I_1$ is a crisp set.
  \item[(E2)] $E(I_1)=1$ if and only if, for any $x\in X$, $\mu_{I_1}(x)=\nu_{I_1}(x)$.
  \item[(E3)] $E(I_1)=E(I_1^{\complement})$.
  \item[(E4)] $E(I_1)\leq E(I_2)$ if, for any $x\in X$, it holds either
  $\mu_{I_1}(x)\leq \mu_{I_2}(x)\leq \nu_{I_2}(x) \leq \nu_{I_1}(x)$ or
  $\mu_{I_1}(x)\geq \mu_{I_2}(x)\geq \nu_{I_2}(x) \geq \nu_{I_1}(x)$.
\end{enumerate}
\end{definition}

\section{The drawbacks of Xiao's distance measure $d_{\widetilde{\chi}}$}
\label{Sec-III}

Let $X=\{x_1, x_2, \ldots, x_n\}$ be a finite UOD and
{\small $I_1=\Big\{\frac{\alpha_{j}^{(1)}}{x_j} \mid
1\leq j\leq n, \alpha_{j}^{(1)}\in \Theta \Big\}$} and {\small $I_2=
\Big\{\frac{\alpha_{j}^{(2)}}{x_j} \mid 1\leq j\leq n,
\alpha_{j}^{(2)}\in \Theta \Big\}$} be two IFSs on $X$.
Based on Jensen-Shannon divergence, Xiao~\cite{Xiao2021}
 introduced a new distance measure $\bm{d}_{\widetilde{\chi}}$
 for IFSs as follows:
\begin{equation}
\label{Eq-Xiao-1}
\begin{split}
& \bm{d}_{\widetilde{\chi}}(I_1, I_2)\\
=&
\frac{1}{n}\sum_{j=1}^{n}
\left[\frac{1}{2}\left( \mu_{I_1}(x_j)  \cdot  \log_{2}
\frac{2 \mu_{I_1}(x_j)}{\mu_{I_1}(x_j)+\mu_{I_2}(x_j)}\right.\right.\\
& \quad \quad \quad +\mu_{I_2}(x_j) \cdot \log_{2}\frac{2 \mu_{I_2}(x_j)}{\mu_{I_1}(x_j)+\mu_{I_2}(x_j)}\\
& \quad \quad \quad +\nu_{I_1}(x_j) \cdot \log_{2}\frac{2 \nu_{I_1}(x_j)}{\nu_{I_1}(x_j)+\nu_{I_2}(x_j)}\\
& \quad \quad \quad +\nu_{I_2}(x_j) \cdot \log_{2}\frac{2 \nu_{I_2}(x_j)}{\nu_{I_1}(x_j)+\nu_{I_2}(x_j)}\\
& \quad \quad \quad +\pi_{I_1}(x_j) \cdot \log_{2}\frac{2 \pi_{I_1}(x_j)}{\pi_{I_1}(x_j)+\pi_{I_2}(x_j)}\\
& \quad \quad \quad \left.\left. +\pi_{I_2}(x_j) \cdot \log_{2}\frac{2 \pi_{I_2}(x_j)}{\pi_{I_1}(x_j)
+\pi_{I_2}(x_j)}\right)\right]^{0.5},
\end{split}
\end{equation}
where $\alpha_{j}^{(1)}=\langle \mu_{I_1}(x_i), \nu_{I_1}(x_j)\rangle$
and $\alpha_{j}^{(2)}=\langle \mu_{I_2}(x_i),$ $\nu_{I_2}(x_j)\rangle$.
Its dual similarity measure $\mathbf{S}_{\widetilde{\chi}}$ is defined by
$\mathbf{S}_{\widetilde{\chi}}(I_1, I_2)=1-\bm{d}_{\widetilde{\chi}}(I_1, I_2)$.

The following examples show that Xiao's distance measure defined by Eq.~\eqref{Eq-Xiao-1}
has some drawbacks: (1) its dual similarity measure $\mathbf{S}_{\widetilde{\chi}}$
does not satisfy the axiomatic definition (S4) of Definition~\ref{Def-Li-Cheng-IFVs};
(2) it cannot effectively distinguish some IFVs/IFSs with obvious size relationship;
(3) Except for the endpoints $\langle 1, 0\rangle$ and $\langle 0, 1 \rangle$,
there exist infinitely many pairs such that the
maximum distance $1$ can be achieved under this distance.

\begin{example}
\label{Exm-Wu-a}
Let the UOD $X=\{x\}$ and choose $I_{1}=\left\{\frac{\langle 0.33, 0.36\rangle}{x}\right\}$,
$I_{2}=\left\{\frac{\langle \frac{1}{3}, \frac{1}{3}\rangle}{x}\right\}$, and
$I_{3}=\left\{\frac{\langle 0.334, 0.333333 \rangle}{x}\right\}\in \mathrm{IFS}(X)$.
Clearly, $I_1\subset I_2 \subset I_3$. Meanwhile, by Eq.~\eqref{Eq-Xiao-1}
and direct calculation, we have $\mathbf{S}_{\widetilde{\chi}}(I_1, I_2)
=1-\bm{d}_{\widetilde{\chi}}(I_1, I_2)=0.9738972$ and $\mathbf{S}_{\widetilde{\chi}}(I_1, I_3)=1-
\bm{d}_{\widetilde{\chi}}(I_1, I_3)=0.9741713$, and thus $\mathbf{S}_{\widetilde{\chi}}(I_1, I_2)
<\mathbf{S}_{\widetilde{\chi}}(I_1, I_3)$. This contradicts the axiomatic definition
(S4) of Definition~\ref{Def-Li-Cheng-IFVs} because $I_1\subset I_2 \subset I_3$.

Furthermore, take
$I_{2}^{(\lambda)}=\left\{\frac{\langle \frac{1}{3}, \lambda\rangle}{x}\right\}$ and
$I_{3}^{(\lambda)}=\left\{\frac{\langle 0.334, \lambda \rangle}{x}\right\}\in \mathrm{IFS}(X)$
in Example~\ref{Exm-Wu-a}. By varying the parameter $\lambda$ from $0$ to $0.36$,
Fig.~\ref{Fig-Exm-Wu-a} visualizes the changing trend of distances between $I_{2}^{(\lambda)}$
and $I_1$, and between $I_{3}^{(\lambda)}$ and $I_1$, based on the distance
measure $\bm{d}_{\widetilde{\chi}}$. Observing from Fig.~\ref{Fig-Exm-Wu-a}, it is clear that
there exists $\lambda^{*}\in (0, 0.36)$ (the intersection of the two curves in Fig.~\ref{Fig-Exm-Wu-a})
such that, for any $0<\lambda<\lambda^{*}$, it holds that
$\bm{d}_{\widetilde{\chi}}(I_{1}, I_{2}^{(\lambda)})> \bm{d}_{\widetilde{\chi}}(I_{1},
I_{3}^{(\lambda)})$. This further illustrates the unreasonableness of the distance
measure $\bm{d}_{\widetilde{\chi}}$ since $I_{1}  \subset I_{2}^{(\lambda)}
\subset I_{3}^{(\lambda)}$.

\begin{figure}[H]
\centering
{\scalebox{0.33}{\includegraphics[]{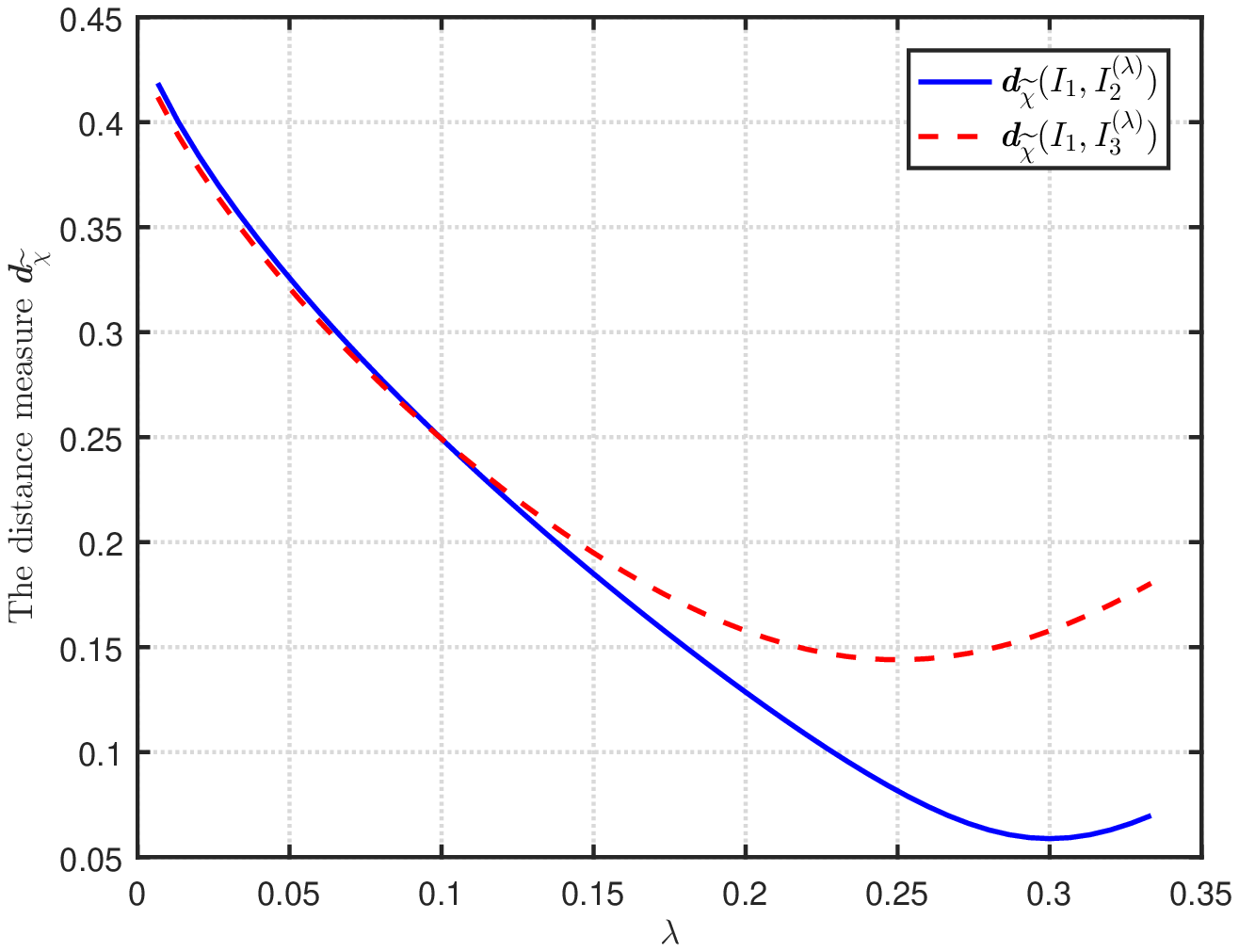}}}
\caption{The distance measures $\bm{d}_{\widetilde{\chi}}
(I_{1}, I_{2}^{(\lambda)})$ and $\bm{d}_{\widetilde{\chi}}
(I_{1}, I_{3}^{(\lambda)})$ in Example~\ref{Exm-Wu-a}}
\label{Fig-Exm-Wu-a}
\end{figure}
\end{example}

Hung and Yang~\cite{HY2004} introduced the $J_{\gamma}$-divergence
to measure the difference between two IFVs as follows: for $\alpha$,
$\beta \in \Theta$,
\begin{itemize}
  \item If $\gamma\in (0, 1) \cup (1, +\infty)$, $J_{\gamma}(\alpha, \beta)
  =\frac{-1}{\gamma-1}\big[(\frac{\mu_{\alpha}+\mu_{\beta}}{2})^{\gamma}
  -\frac{1}{2}(\mu_{\alpha}^{\gamma}+\mu_{\beta}^{\gamma})+
  (\frac{\nu_{\alpha}+\nu_{\beta}}{2})^{\gamma}
  -\frac{1}{2}(\nu_{\alpha}^{\gamma}+\nu_{\beta}^{\gamma})+
  (\frac{\pi_{\alpha}+\pi_{\beta}}{2})^{\gamma}
  -\frac{1}{2}(\pi_{\alpha}^{\gamma}+\pi_{\beta}^{\gamma})\big]$;
  \item If $\gamma=1$, $J_{\gamma}(\alpha, \beta)
  =\frac{-1}{2}\big[(\mu_{\alpha}+\mu_{\beta})
  \ln (\frac{\mu_{\alpha}+\mu_{\beta}}{2})-\mu_{\alpha}\cdot \ln \mu_{\alpha}
  -\mu_{\beta}\cdot \ln \mu_{\beta}+ (\nu_{\alpha}+\nu_{\beta})
  \ln (\frac{\nu_{\alpha}+\nu_{\beta}}{2})-\nu_{\alpha}\cdot \ln \nu_{\alpha}
  -\nu_{\beta}\cdot \ln \nu_{\beta}+ (\pi_{\alpha}+\pi_{\beta})
  \ln (\frac{\pi_{\alpha}+\pi_{\beta}}{2})-\pi_{\alpha}\cdot \ln \pi_{\alpha}
  -\pi_{\beta}\cdot \ln \pi_{\beta}\big]$.
\end{itemize}

\begin{remark}
Joshi and Kumar \cite{JK2018} introduced an IF dissimilarity divergence
$\widehat{DJS}$, which is equivalent to the $J_{1}$-divergence of Hung and Yang~\cite{HY2004}.
By direct calculation, it can be verified that $\sqrt{\widehat{DJS}(\alpha, \beta)}
=\sqrt{J_{1}(\alpha, \beta)}=\ln 2 \cdot \bm{d}_{\chi}(\alpha, \beta)$, and thus
Example~\ref{Exm-Wu-a} also indicates that \cite[Theorem 1~(D3)]{HY2004}
and \cite[Subsection~3.3]{JK2018} do not hold.
\end{remark}

The following example illustrates the unreasonableness
of Xiao's distance measure $\bm{d}_{\widetilde{\chi}}$
from another perspective.

\begin{example}
\label{Exm-Wu-a1}
Let the UOD $X=\{x\}$ and choose $I_{1}=
\left\{\frac{\langle \frac{1}{3}, \frac{1}{3} \rangle}{x}\right\}$
and $I_{2}^{(\lambda)}=\left\{\frac{\langle \lambda, 0.00001\rangle}{x}\right\}\in \mathrm{IFS}(X)$.
By varying the parameter $\lambda$ from $\frac{1}{3}$ to $1$,
Fig.~\ref{Fig-Exm-Wu-a1} visualizes the changing trend of distances between $I_{2}^{(\lambda)}$
and $I_1$ by using Xiao's distance measure $\bm{d}_{\widetilde{\chi}}$.
Observing from Fig.~\ref{Fig-Exm-Wu-a1}, it is clear that, for any $\lambda_1$,
$\lambda_2\in (\frac{1}{3}, 0.5)$ with $\lambda_1<\lambda_2$, it holds that
$\bm{d}_{\widetilde{\chi}}(I_{1}, I_{2}^{(\lambda_1)})> \bm{d}_{\widetilde{\chi}}(I_{1},
I_{2}^{(\lambda_2)})$, which contradicts the fact that $I_{1}\subset I_{2}^{\lambda_1}
\subset I_{2}^{\lambda_2}$. This also illustrates the unreasonableness of the distance
measure $\bm{d}_{\widetilde{\chi}}$.
\begin{figure}[H]
\centering
{\scalebox{0.33}{\includegraphics[]{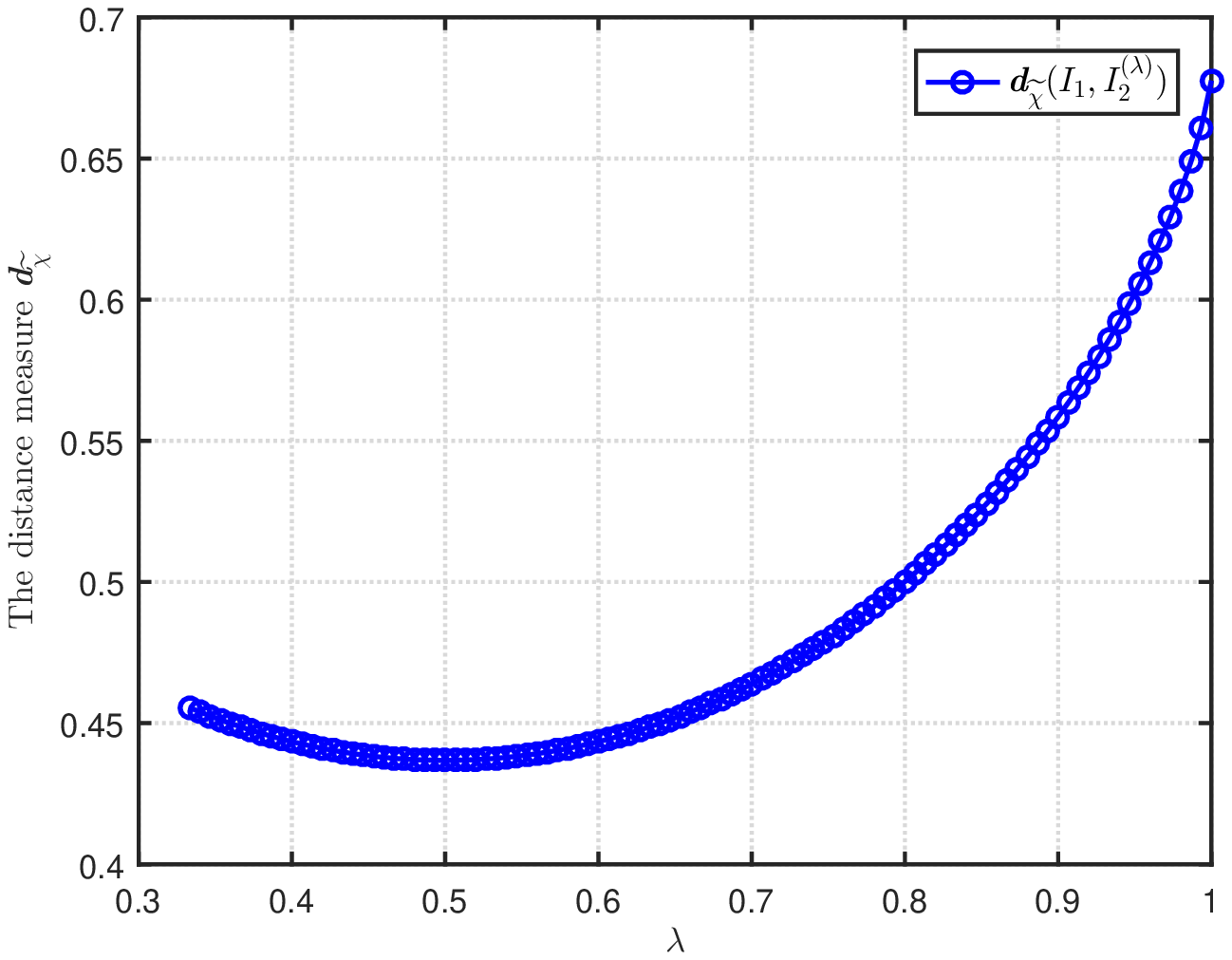}}}
\caption{The distance measure $\bm{d}_{\widetilde{\chi}}
(I_{1}, I_{2}^{(\lambda)})$ in Example~\ref{Exm-Wu-a1}}
\label{Fig-Exm-Wu-a1}
\end{figure}
\end{example}

\begin{example}
\label{Exm-Wu-2-Sec2}
Assume that the IFSs $I_1$, $I_1^{\prime}$, $I_2$, and $I_3$ on UOD $X=\{x\}$ are given by
$I_1=\left\{\frac{\langle 1, 0\rangle}{x}\right\}$,
$I_1^{\prime}=\left\{\frac{\langle 0, 1\rangle}{x}\right\}$,
$I_2^{(\lambda)}=\left\{\frac{\langle \lambda, 0\rangle}{x}\right\}$,
and
$I_3^{(\lambda)}=\left\{\frac{\langle \lambda, 1-\lambda\rangle}{x}\right\}$, where
$0\leq \lambda\leq 1$.
By direct calculation, it follows from Eq.~\eqref{Eq-Xiao-1}
that
{\small\begin{equation}
\label{eq-Exm-2-1-Sec2}
\begin{split}
& \bm{d}_{\widetilde{\chi}}(I_1, I_2^{(\lambda)})\\
=& \sqrt{\frac{1}{2}\left[\log_{2}\frac{2}{1+\lambda}
+\lambda\log_{2}\frac{2\lambda}{1+\lambda}+(1-\lambda)\log_{2}\frac{2(1-\lambda)}{1-\lambda}\right]}\\
=& \sqrt{\frac{1}{2}\left[\log_{2}\frac{2}{1+\lambda}
+\lambda\log_{2}\frac{2\lambda}{1+\lambda}+(1-\lambda)\right]}\\
=& \bm{d}_{\widetilde{\chi}}(I_1, I_3^{(\lambda)}),
\end{split}
\end{equation}}
and
{\small \begin{equation}
\label{eq-Exm-2-2-Sec2}
\begin{split}
&\bm{d}_{\widetilde{\chi}}(I_1^{\prime}, I_2^{(\lambda)})\\
=& \sqrt{\frac{1}{2}\left[\lambda \log_{2}\frac{2\lambda}{\lambda}
+1\cdot \log_{2} \frac{2}{1}+(1-\lambda)\log_{2}\frac{2(1-\lambda)}{1-\lambda}\right]}\\
=& 1.
\end{split}
\end{equation}}

(i) For $\lambda\neq 1$, it is clear that $I_{2}^{(\lambda)}\supsetneqq I_{3}^{(\lambda)}$.
It follows from Eq.~\eqref{eq-Exm-2-1-Sec2} that Xiao's distance measure
$\bm{d}_{\widetilde{\chi}}$ cannot effectively distinguish $I_{2}^{(\lambda)}$
and $I_{3}^{(\lambda)}$ form $I_{1}$, outputting an unreasonable result with
$I_{1}\supset I_{2}^{(\lambda)}\supsetneqq I_{3}^{(\lambda)}$.

(ii) For $0\leq \lambda_1< \lambda_2\leq 1$, it is clear that
$I_{2}^{(\lambda_1)}\subsetneqq I_{2}^{(\lambda_2)}$. It follows from
Eq.~\eqref{eq-Exm-2-2-Sec2} that Xiao's distance measure $\bm{d}_{\widetilde{\chi}}$
cannot effectively distinguish $I_{2}^{(\lambda_1)}$ and $I_{2}^{(\lambda_2)}$ with
$\lambda_1\neq \lambda_2$ from $I_{1}^{\prime}$, also outputting an unreasonable result
with $I_{2}^{(\lambda_1)}\subsetneqq I_{2}^{(\lambda_2)}\subset I_{1}$.
Moreover, it follows from Eq.~\eqref{eq-Exm-2-2-Sec2} that there exist infinite
IFSs $I_{2}^{(\lambda)}$ ($\lambda\in [0, 1]$) such that the distance
from $I_1^{\prime}$ is equal to the maximum value of $1$.
\end{example}

\section{The drawbacks of Yang and Chiclana's spherical distance $d_{_{\mathrm{YC}}}$}
\label{Sec-IV}

Let $X=\{x_1, x_2, \ldots, x_n\}$ be a finite UOD and
$I_1=\Big\{\frac{\alpha_{j}^{(1)}}{x_j} \mid
1\leq j\leq n, \alpha_{j}^{(1)}\in \Theta \Big\}$ and $I_2=
\Big\{\frac{\alpha_{j}^{(2)}}{x_j} \mid 1\leq j\leq n,
\alpha_{j}^{(2)}\in \Theta \Big\}$ be two IFSs on $X$.
Based on the $\arccos$ function, Yang and Chiclana~\cite{YC2009}
introduced a spherical distance $d_{_{\mathrm{YC}}}$
for IFSs as follows:
\begin{small}
\begin{equation}
\label{Eq-YC-1}
\begin{split}
& d_{_{\mathrm{YC}}}(I_1, I_2)\\
=&
\frac{2}{n\pi}\sum_{j=1}^{n}
\arccos\left(\sqrt{\mu_{I_1}(x_{j})\mu_{I_{2}}(x_{j})}+
\sqrt{\nu_{I_1}(x_{j})\nu_{I_{2}}(x_{j})}\right.\\
& \quad \quad \quad \quad \quad \quad \left.+
\sqrt{\pi_{I_1}(x_{j})\pi_{I_{2}}(x_{j})}\right),
\end{split}
\end{equation}
\end{small}
where $\alpha_{j}^{(1)}=\langle \mu_{I_1}(x_j), \nu_{I_1}(x_j)\rangle$
and $\alpha_{j}^{(2)}=\langle \mu_{I_2}(x_j),$ $\nu_{I_2}(x_j)\rangle$.
Its dual similarity measure $\mathbf{S}_{_{\mathrm{YC}}}$ is defined by
$\mathbf{S}_{_{\mathrm{YC}}}(I_1, I_2)=1-d_{_{\mathrm{YC}}}(I_1, I_2)$.

The following examples show that Yang and Chiclana's spherical distance
has the same drawbacks as Xiao's distance measure.

\begin{example}
\label{Exm-8-YC}
Let $X=\{x\}$ and choose $I_{1}=\left\{\frac{\langle 0.5, 0.5\rangle}{x}\right\}$,
$I_{2}=\left\{\frac{\langle 0.6, 0.3\rangle}{x}\right\}$, and
$I_{3}=\left\{\frac{\langle 0.7, 0.3 \rangle}{x}\right\}\in \mathrm{IFS}(X)$.
Clearly, $I_1\subset I_2 \subset I_3$. Meanwhile, by Eq.~\eqref{Eq-YC-1}
and direct calculation, we have $\mathbf{S}_{_{\mathrm{YC}}}(I_{1}, I_{2})=
1-d_{_{\mathrm{YC}}}(I_{1}, I_{2})=1-\frac{2}{\pi}\arccos (\sqrt{0.3}+\sqrt{0.15})
<1-\frac{2}{\pi}\arccos (\sqrt{0.35}+\sqrt{0.15})=1-d_{_{\mathrm{YC}}}(I_{1}, I_{3})
=\mathbf{S}_{_{\mathrm{YC}}}(I_{1}, I_{3})$. This contradicts the axiomatic definition
(S4) of Definition~\ref{Def-Li-Cheng-IFVs} because $I_1\subset I_2 \subset I_3$.
\end{example}

\begin{example}
\label{Exm-YC-2-Sec2}
Assume that IFSs $I_1$, $I_2$, and $I_3$ on UOD $X=\{x\}$ are given by
$I_1=\left\{\frac{\langle 1, 0\rangle}{x}\right\}$,
$I_2^{(\lambda)}=\left\{\frac{\langle \lambda, 0\rangle}{x}\right\}$,
$\tilde{I}_{2}^{(\lambda)}=\left\{\frac{\langle 0, \lambda\rangle}{x}\right\}$, and
$I_3^{(\lambda)}=\left\{\frac{\langle \lambda, 1-\lambda\rangle}{x}\right\}$,
 where
$0\leq \lambda\leq 1$.
By direct calculation, we have
\begin{equation}
\label{eq-Exm-YC-1}
d_{_{\mathrm{YC}}}(I_1, I_2^{(\lambda)})=d_{_{\mathrm{YC}}}(I_1, I_3^{(\lambda)})
=\frac{2}{\pi}\arccos(\sqrt{\lambda}),
\end{equation}
and
\begin{equation}
\label{eq-Exm-YC-2}
d_{_{\mathrm{YC}}}(I_1, \tilde{I}_2^{(\lambda)})=1.
\end{equation}

(i) For $\lambda\neq 1$, it is clear that $I_{2}^{(\lambda)}\supsetneqq I_{3}^{(\lambda)}$.
It follows from Eq.~\eqref{eq-Exm-YC-1} that Yang and Chiclana's spherical distance
$d_{_{\mathrm{YC}}}$ cannot effectively distinguish $I_{2}^{(\lambda)}$
and $I_{3}^{(\lambda)}$ form $I_{1}$, outputting a unreasonable result with
$I_{1}\supset I_{2}^{(\lambda)}\supsetneqq I_{3}^{(\lambda)}$.

(ii) For $\lambda \neq 0$, it is clear that
$\tilde{I}_{2}^{(\lambda)}\subsetneqq I_{2}^{(\lambda)}$. it follows from
Eq.~\eqref{eq-Exm-YC-2} that Yang and Chiclana's spherical distance
$d_{_{\mathrm{YC}}}$ cannot effectively distinguish $\tilde{I}_{2}^{(\lambda)}$
and $I_{2}^{(\lambda)}$ with $\lambda \neq 0$ from $I_{1}$, also outputting
a unreasonable result with $\tilde{I}_{2}^{(\lambda)}\subsetneqq I_{2}^{(\lambda)}
\subset I_{1}$. Meanwhile, it follows from Eq.~\eqref{eq-Exm-YC-2} that there
exist infinite IFSs $\tilde{I}_{2}^{(\lambda)}$ ($\lambda\in [0, 1]$), where the
distance from $I_1$ is equal to the maximum value $1$.
\end{example}

\section{A novel SIFDisM/SIFSimM based on Jensen-Shannon divergence}
\label{Sec-V}

In this section, we propose a new strict distance measure and a new
strict similarity measure for IFVs and IFSs based on Jensen-Shannon divergence,
which can overcome the drawbacks of Xiao's distance measure and Hung and Yang's
$J_{\gamma}$-divergence discussed in the previous section.

\subsection{A new distance/similarity measure on IFVs}

Clearly, an IFV $\alpha=\langle \mu, \nu \rangle$ can be equivalently
expressed as an interval $[\nu, 1-\mu]$, and thus we can use
$-(\nu \ln \nu+(1-\mu) \ln (1-\mu))$ to express the Shannon entropy $H(\alpha)$ of
$\alpha$. For $\alpha=\langle \mu_{\alpha}, \nu_{\alpha} \rangle$,
$\beta=\langle \mu_{\beta}, \nu_{\beta} \rangle \in \Theta$, by
applying Jensen-Shannon divergence, we define the \textit{Jensen-Shannon IF divergence measure} $\bm{JS}_{_{\mathrm{IF}}}(\alpha, \beta)$
between $\alpha$ and $\beta$ as follows:
\begin{equation}
\label{Eq-Wu-1-*}
\begin{split}
\bm{JS}_{_{\mathrm{IF}}}(\alpha, \beta)=&
H(\frac{\alpha+\beta}{2})-\frac{1}{2}H(\alpha)-\frac{1}{2}H(\beta)\\
=& \frac{1}{2}\left[ (1-\mu_{\alpha}) \cdot  \ln
\frac{2 (1-\mu_{\alpha})}{(1-\mu_{\alpha})+(1-\mu_{\beta})}\right.\\
& \quad +(1-\mu_{\beta}) \cdot \ln\frac{2 (1-\mu_{\beta})}{(1-\mu_{\alpha})+(1-\mu_{\beta})}\\
& \left. \quad +\nu_{\alpha} \cdot \ln\frac{2 \nu_{\alpha}}{\nu_{\alpha}+\nu_{\beta}}
 +\nu_{\beta} \cdot \ln\frac{2 \nu_{\beta}}{\nu_{\alpha}+\nu_{\beta}}\right],
\end{split}
\end{equation}
where $0 \cdot \ln 0=0\cdot \ln \frac{0}{0+0}=0$ and $\frac{\alpha+\beta}{2}
=\langle \frac{\mu_{\alpha}+\mu_{\beta}}{2}, \frac{\nu_{\alpha}+\nu_{\beta}}{2}\rangle$.

Let $\mathds{R}^{+}=\{x\in \mathds{R}\mid x\geq 0\}$. To obtain a new metric
for probability distributions, Endres and Schindelin~\cite{ES2003} introduced
the following function $L$ and presented some basic properties.

\begin{definition}
[{\textrm{\protect\cite[Definition~1]{ES2003}}}]
Define the function $L(p, q): \mathds{R}^{+}\times
\mathds{R}^{+} \longrightarrow \mathds{R}^{+}$
by
$$
L(p, q)=p\cdot \log_{2}\frac{2p}{p+q}+q\cdot \log_{2}\frac{2q}{p+q}.
$$
\end{definition}

\begin{remark}
Endres and Schindelin~\cite{ES2003} proved that the function $L(\_)$ is
well defined, i.e., $L(p, q)\geq 0$ holds for all $p$, $q\in \mathds{R}^{+}$.
\end{remark}

\begin{lemma}
[{\textrm{\protect\cite[Lemma~2]{ES2003}}}]
\label{Lemma-ES}
Let $p$, $q$, $r\in \mathds{R}^{+}$. Then,
$$
\sqrt{L(p, q)}\leq \sqrt{L(p, r)}+\sqrt{L(r, q)}.
$$
\end{lemma}

For simplicity of presentation, denote $\mathcal{Z}(\alpha, \beta)=L(1-\mu_{\alpha},
1-\mu_{\beta})+L(\nu_{\alpha}, \nu_{\beta})$.

Define a function $\zeta(x)=x \cdot \log_{2}(2x)+(1-x)\cdot \log_{2}(2(1-x))$
($x\in (0, 1)$). By direct derivation, we have
$$
\zeta^{\prime}(x)=\log_{2}\frac{x}{1-x}=
\begin{cases}
< 0, & 0<x< 0.5, \\
> 0, & 0.5< x<1,
\end{cases}
$$
implying that $\zeta(\_)$ is strictly decreasing on $(0, 0.5)$ and
strictly increasing on $(0.5, 1)$. Thus, for any $x\in (0, 1)$,
\begin{equation}
\label{Zata-Value-*}
0=\zeta(0.5)\leq \zeta(x)\leq \max\{\lim_{x\to 0^{+}}\zeta(x),
\lim_{x\to 1^{-}}\zeta(x)\}=1\footnote{Because we define $0\cdot \log_{2}0 =0$, one has
$\lim_{x\to 0^{+}}\zeta(x)=1=\zeta(0)$ and $\lim_{x\to 1^{-}}\zeta(x)=1=\zeta(1)$.}.
\end{equation}
Meanwhile, by direct calculation, we get
\begin{equation}
\label{Z-Function-*}
\begin{split}
& \mathcal{Z}(\alpha, \beta)\\
=& (1-\mu_{\alpha}) \cdot  \log_{2}
\frac{2 (1-\mu_{\alpha})}{(1-\mu_{\alpha})+(1-\mu_{\beta})}\\
&  +(1-\mu_{\beta}) \cdot \log_{2}\frac{2 (1-\mu_{\beta})}{(1-\mu_{\alpha})+(1-\mu_{\beta})}\\
&  +\nu_{\alpha} \cdot \log_{2}\frac{2 \nu_{\alpha}}{\nu_{\alpha}+\nu_{\beta}}
 +\nu_{\beta} \cdot \log_{2}\frac{2 \nu_{\beta}}{\nu_{\alpha}+\nu_{\beta}} \\
=& ((1-\mu_{\alpha})+(1-\mu_{\beta}))\cdot \zeta\left(\frac{(1-\mu_{\alpha})}{(1-\mu_{\alpha})+(1-\mu_{\beta})}\right)\\
& + (\nu_{\alpha}+\nu_{\beta})\cdot \zeta\left(\frac{\nu_{\alpha}}{\nu_{\alpha}+\nu_{\beta}}\right).
\end{split}
\end{equation}
This, together with Eq.~\eqref{Zata-Value-*}, implies that
$\mathcal{Z}(\alpha, \beta) \geq 0$, and thus $\bm{JS}_{_{\mathrm{IF}}}(\alpha, \beta)
=\frac{\ln 2}{2}\cdot \mathcal{Z}(\alpha, \beta) \geq 0$ by applying Eq.~\eqref{Eq-Wu-1-*}.

Applying the square root of $\bm{JS}_{_{\mathrm{IF}}}$, we define the \textit{normalized
Jensen-Shannon IF divergence measure} $\overline{\bm{JS}}_{_{\mathrm{IF}}}(\alpha, \beta)$
between $\alpha$ and $\beta$ as follows:
{\small \begin{equation}
\label{Eq-Wu-1a-*}
\begin{split}
& \overline{\bm{JS}}_{_{\mathrm{IF}}}(\alpha, \beta)=
\sqrt{\frac{\bm{JS}_{_{\mathrm{IF}}}(\alpha, \beta)}{\ln 2}}\\
=& \left[\frac{1}{2}\left( (1-\mu_{\alpha}) \cdot  \log_2
\frac{2 (1-\mu_{\alpha})}{(1-\mu_{\alpha})+(1-\mu_{\beta})}\right.\right.\\
& \quad +(1-\mu_{\beta}) \cdot \log_2\frac{2 (1-\mu_{\beta})}{(1-\mu_{\alpha})+(1-\mu_{\beta})}\\
& \left. \left. \quad +\nu_{\alpha} \cdot \log_2 \frac{2 \nu_{\alpha}}{\nu_{\alpha}+\nu_{\beta}}
 +\nu_{\beta} \cdot \log_2 \frac{2 \nu_{\beta}}{\nu_{\alpha}+\nu_{\beta}}\right)\right]^{0.5}.
\end{split}
\end{equation}}
Clearly, $\overline{\bm{JS}}_{_{\mathrm{IF}}}(\alpha, \beta)
=\sqrt{\frac{1}{2}(L(1-\mu_{\alpha}, 1-\mu_{\beta})+L(\nu_{\alpha}, \nu_{\beta}))}$
$=\sqrt{\frac{\mathcal{Z}(\alpha, \beta)}{2}}$.

The function $\overline{\bm{JS}}_{_{\mathrm{IF}}}$ has the following desirable properties.

\begin{property}
\label{Pro-1-IFV}
$\overline{\bm{JS}}_{_{\mathrm{IF}}}(\alpha, \beta)
=\overline{\bm{JS}}_{_{\mathrm{IF}}}(\beta, \alpha)$.
\end{property}

\begin{IEEEproof}
It follows directly from Eq.~\eqref{Eq-Wu-1a-*}.
\end{IEEEproof}

\begin{lemma}
\label{Lemma-1-*}
Fix an IFV $\alpha=\langle \mu_{\alpha}, \nu_{\alpha}\rangle\in \Theta$. Then, for any $\beta_1$,
$\beta_2\in \Theta$ with $\alpha\subset \beta_1\subset \beta_2$, we have
$\mathcal{Z}(\alpha, \beta_1)=\mathcal{Z}(\beta_1, \alpha)\leq
\mathcal{Z}(\alpha, \beta_2)=\mathcal{Z}(\beta_2, \alpha)$
and $\overline{\bm{JS}}_{_{\mathrm{IF}}}(\alpha, \beta_1)\leq
\overline{\bm{JS}}_{_{\mathrm{IF}}}(\alpha, \beta_2)$.
\end{lemma}

\begin{IEEEproof}
For any $\beta=\langle \mu, \nu\rangle \in \Theta$, it is clear that
\begin{equation}
\label{eq-Wu-2-*}
\begin{split}
\mathcal{Z}(\alpha, \beta)
=& \left[ (1-\mu_{\alpha}) \cdot  \log_2
\frac{2 (1-\mu_{\alpha})}{(1-\mu_{\alpha})+(1-\mu)}\right.\\
& \quad +(1-\mu) \cdot \log_2\frac{2 (1-\mu)}{(1-\mu_{\alpha})+(1-\mu)}\\
& \left. \quad +\nu_{\alpha} \cdot \log_2 \frac{2 \nu_{\alpha}}{\nu_{\alpha}+\nu}
 +\nu \cdot \log_2 \frac{2 \nu}{\nu_{\alpha}+\nu}\right].
\end{split}
\end{equation}
By Eq.~\eqref{eq-Wu-2-*} and direct calculation, we have that, for
$\mu\geq \mu_{\alpha}$ and $\nu\leq \nu_{\alpha}$,
\begin{equation}
\label{eq-Wu-3-*}
\begin{split}
\frac{\partial \mathcal{Z}}{\partial \mu}=&
\frac{1-\mu_{\alpha}}{\ln 2}\cdot \frac{1}{(1-\mu_{\alpha})+(1-\mu)}\\
& -\log_{2}\frac{2(1-\mu)}{(1-\mu)+(1-\mu_{\alpha})}\\
& -\frac{1-\mu_{\alpha}}{\ln 2}\cdot \frac{1}{(1-\mu_{\alpha})+(1-\mu)}\\
=& -\log_{2}\frac{2(1-\mu)}{(1-\mu)+(1-\mu_{\alpha})}\geq 0,
\end{split}
\end{equation}
and
\begin{equation}
\label{eq-Wu-4-*}
\begin{split}
\frac{\partial \mathcal{Z}}{\partial \nu}=&
-\frac{\nu_{\alpha}}{(\nu_{\alpha}+\nu)\cdot \ln 2}\\
& +\log_{2}\frac{2\nu}{\nu_{\alpha}+\nu}
+\frac{\nu_{\alpha}}{(\nu_{\alpha}+\nu)\cdot \ln 2}\\
=&\log_{2}\frac{2\nu}{\nu_{\alpha}+\nu}\leq 0.
\end{split}
\end{equation}
Let $\beta_1=\langle \mu_{\beta_1}, \nu_{\beta_1}\rangle$ and $\beta_2=\langle \mu_{\beta_2},
\nu_{\beta_2}\rangle$. From $\alpha\subset \beta_1\subset \beta_2$,
it follows that $\mu_{\alpha}\leq \mu_{\beta_1} \leq \mu_{\beta_2}$ and
$\nu_{\alpha}\geq \nu_{\beta_1} \geq \nu_{\beta_2}$. Together with Eqs.~\eqref{eq-Wu-3-*}
and \eqref{eq-Wu-4-*}, we get
$\mathcal{Z}(\alpha, \beta_1)
\leq \mathcal{Z}(\alpha, \langle \mu_{\beta_1}, \nu_{\beta_2}\rangle)
\leq \mathcal{Z}(\alpha, \langle \mu_{\beta_2}, \nu_{\beta_2}\rangle)
=\mathcal{Z}(\alpha, \beta_2)$.
\end{IEEEproof}

\begin{lemma}
\label{Lemma-1-*-strict}
Fix an IFV $\alpha=\langle \mu_{\alpha}, \nu_{\alpha}\rangle\in \Theta$. Then, for any $\beta_1$,
$\beta_2\in \Theta$ with $\alpha\subsetneqq \beta_1\subsetneqq \beta_2$, we have
$\mathcal{Z}(\alpha, \beta_1)=\mathcal{Z}(\beta_1, \alpha)<
\mathcal{Z}(\alpha, \beta_2)=\mathcal{Z}(\beta_2, \alpha)$
and $\overline{\bm{JS}}_{_{\mathrm{IF}}}(\alpha, \beta_1)<
\overline{\bm{JS}}_{_{\mathrm{IF}}}(\alpha, \beta_2)$.
\end{lemma}

\begin{IEEEproof}
For any $\beta=\langle \mu, \nu\rangle \in \Theta$ with
$\alpha\subsetneqq \beta$, by Eqs.~\eqref{eq-Wu-3-*}
and \eqref{eq-Wu-4-*}, we have

(1) If $\mu> \mu_{\alpha}$ and $\nu\leq \nu_{\alpha}$, then
\begin{equation}
\label{Eq-strict-1}
\frac{\partial \mathcal{Z}}{\partial \mu}
=-\log_{2}\frac{2(1-\mu)}{(1-\mu)+(1-\mu_{\alpha})}> 0.
\end{equation}

(2) If $\mu\geq \mu_{\alpha}$ and $\nu< \nu_{\alpha}$, then
\begin{equation}
\label{Eq-strict-4}
\frac{\partial \mathcal{Z}}{\partial \nu}
=\log_{2}\frac{2\nu}{\nu_{\alpha}+\nu}< 0.
\end{equation}
For $\beta_1=\langle \mu_{\beta_1}, \nu_{\beta_1}\rangle$ and $\beta_2=\langle \mu_{\beta_2},
\nu_{\beta_2}\rangle$ with $\alpha\subsetneqq \beta_1\subsetneqq \beta_2$,
consider the following cases:

2.1) If $\mu_{\alpha}< \mu_{\beta_1}< \mu_{\beta_2}$ and
$\nu_{\alpha}\geq \nu_{\beta_1} \geq \nu_{\beta_2}$, then, by Eqs.~\eqref{eq-Wu-4-*}
and \eqref{Eq-strict-1}, we have
\begin{align*}
& \mathcal{Z}(\alpha, \beta_1)=
\mathcal{Z}(\alpha, \langle \mu_{\beta_1}, \nu_{\beta_1}\rangle)\\
<& \mathcal{Z}(\alpha, \langle \mu_{\beta_2}, \nu_{\beta_1}\rangle) \quad
\text{(by Eq.~\eqref{Eq-strict-1})}\\
\leq & \mathcal{Z}(\alpha, \langle \mu_{\beta_2}, \nu_{\beta_2}\rangle)=
\mathcal{Z}(\alpha, \beta_2) \quad \text{(by Eq.~\eqref{eq-Wu-4-*})}.
\end{align*}

2.2) If $\mu_{\alpha}\leq \mu_{\beta_1} \leq \mu_{\beta_2}$ and
$\nu_{\alpha}> \nu_{\beta_1}> \nu_{\beta_2}$, then, by Eqs.~\eqref{eq-Wu-3-*}
and \eqref{Eq-strict-4}, we have
\begin{align*}
& \mathcal{Z}(\alpha, \beta_1)=
\mathcal{Z}(\alpha, \langle \mu_{\beta_1}, \nu_{\beta_1}\rangle)\\
<& \mathcal{Z}(\alpha, \langle \mu_{\beta_1}, \nu_{\beta_2}\rangle) \quad
\text{(by Eq.~\eqref{Eq-strict-4})}\\
\leq & \mathcal{Z}(\alpha, \langle \mu_{\beta_2}, \nu_{\beta_2}\rangle)=
\mathcal{Z}(\alpha, \beta_2) \quad \text{(by Eq.~\eqref{eq-Wu-3-*})}.
\end{align*}

2.3) If $\mu_{\alpha}< \mu_{\beta_1} \leq \mu_{\beta_2}$ and
$\nu_{\alpha}\geq \nu_{\beta_1}> \nu_{\beta_2}$, then, by
Eqs.~\eqref{eq-Wu-3-*} and \eqref{Eq-strict-4}, we have
\begin{align*}
& \mathcal{Z}(\alpha, \beta_1)=
\mathcal{Z}(\alpha, \langle \mu_{\beta_1}, \nu_{\beta_1}\rangle)\\
\leq & \mathcal{Z}(\alpha, \langle \mu_{\beta_2}, \nu_{\beta_1}\rangle) \quad
\text{(by Eqs.~\eqref{eq-Wu-3-*})}\\
< & \mathcal{Z}(\alpha, \langle \mu_{\beta_2}, \nu_{\beta_2}\rangle)=
\mathcal{Z}(\alpha, \beta_2) \quad \text{(by Eq.~\eqref{Eq-strict-4})}.
\end{align*}

2.4) If $\mu_{\alpha}\leq \mu_{\beta_1} < \mu_{\beta_2}$ and
$\nu_{\alpha}> \nu_{\beta_1} \geq \nu_{\beta_2}$, then, by
Eqs.~\eqref{eq-Wu-4-*} and \eqref{Eq-strict-1}, we have
\begin{align*}
& \mathcal{Z}(\alpha, \beta_1)=
\mathcal{Z}(\alpha, \langle \mu_{\beta_1}, \nu_{\beta_1}\rangle)\\
< & \mathcal{Z}(\alpha, \langle \mu_{\beta_2}, \nu_{\beta_1}\rangle) \quad
\text{(by Eqs.~\eqref{Eq-strict-1})}\\
\leq & \mathcal{Z}(\alpha, \langle \mu_{\beta_2}, \nu_{\beta_2}\rangle)=
\mathcal{Z}(\alpha, \beta_2) \quad \text{(by Eq.~\eqref{eq-Wu-4-*})}.
\end{align*}
\end{IEEEproof}

Similarly, we have the following results.

\begin{lemma}
\label{Lemma-2-*}
Fix an IFV $\alpha=\langle \mu_1, \nu_1\rangle\in \Theta$. Then, for any $\beta_1$,
$\beta_2\in \Theta$ with $\beta_2 \subset \beta_1 \subset \alpha$, we have
$\mathcal{Z}(\alpha, \beta_1)=\mathcal{Z}(\beta_1, \alpha)
\leq \mathcal{Z}(\alpha, \beta_2)=\mathcal{Z}(\beta_2, \alpha)$
and $\overline{\bm{JS}}_{_{\mathrm{IF}}}(\alpha, \beta_1)\leq
\overline{\bm{JS}}_{_{\mathrm{IF}}}(\alpha, \beta_2)$.
\end{lemma}

\begin{lemma}
\label{Lemma-2-*-strict}
Fix an IFV $\alpha=\langle \mu_1, \nu_1\rangle\in \Theta$. Then, for any $\beta_1$,
$\beta_2\in \Theta$ with $\beta_2\subsetneqq \beta_1 \subsetneqq \alpha$, we have
$\mathcal{Z}(\alpha, \beta_1)=\mathcal{Z}(\beta_1, \alpha)
< \mathcal{Z}(\alpha, \beta_2)=\mathcal{Z}(\beta_2, \alpha)$
and $\overline{\bm{JS}}_{_{\mathrm{IF}}}(\alpha, \beta_1)<
\overline{\bm{JS}}_{_{\mathrm{IF}}}(\alpha, \beta_2)$.
\end{lemma}

\begin{property}
\label{Pro-2-IFV}
$0\leq \overline{\bm{JS}}_{_{\mathrm{IF}}}(\alpha, \beta) \leq 1$.
\end{property}

\begin{IEEEproof}
Clearly, $\overline{\bm{JS}}_{_{\mathrm{IF}}}(\alpha, \beta)
\geq 0$. It suffices to check that $\mathcal{Z}(\alpha, \beta)\leq 2$ by
$\overline{\bm{JS}}_{_{\mathrm{IF}}}(\alpha, \beta) =
\sqrt{\frac{1}{2}\mathcal{Z}(\alpha, \beta)}$.

Fix an IFV $\beta=\langle \mu_{\beta}, \nu_{\beta}\rangle \in \Theta$.
To prove that $\mathcal{Z} (\alpha, \beta)\leq 2$ holds for all
$\langle \mu_{\alpha}, \mu_{\alpha}\rangle\in \Theta$, we consider the following
four cases:

(1) If $\mu_{\alpha}\geq \mu_{\beta}$ and $\nu_{\alpha} \leq \nu_{\beta}$, then
$\beta \subset \alpha \subset \langle 1, 0 \rangle$. This, together with
 Lemma~\ref{Lemma-1-*}, implies that $\mathcal{Z}(\alpha, \beta)
\leq \mathcal{Z}(\langle 1, 0\rangle, \beta)=(1-\mu_{\beta})\cdot \log_{2}\frac{2(1-\mu_{\beta})}
{1-\mu_{\beta}}+\nu_{\beta}\cdot \log_{2}\frac{2\nu_{\beta}}{\nu_{\beta}}
\leq (1-\mu_{\beta})+\nu_{\beta} \leq 2$;

(2) If $\mu_{\alpha}\leq \mu_{\beta}$ and $\nu_{\alpha} \geq \nu_{\beta}$, then
$\langle 0, 1 \rangle \subset \alpha\subset \beta$. This, together with
Lemma~\ref{Lemma-2-*}, implies that $\mathcal{Z}(\alpha, \beta) \leq \mathcal{Z}(\langle 0, 1\rangle, \beta)
=\log_{2}\frac{2}{2-\mu_{\beta}}+\underbrace{(1-\mu_{\beta})\cdot \log_{2}\frac{2(1-\mu_{\beta})}
{2-\mu_{\beta}}}\limits_{\leq 0}+\log_{2}\frac{2}{1+\nu_{\beta}}+\underbrace{\nu_{\beta}\cdot
\log_{2}\frac{2\nu_{\beta}}{1+\nu_{\beta}}}\limits_{\leq 0}\leq \log_{2}\frac{2}{2-\mu_{\beta}}
+\log_{2}\frac{2}{1+\nu_{\beta}} \leq 2$;

(3) If $\mu_{\alpha}\leq \mu_{\beta}$ and $\nu_{\alpha}\leq \nu_{\beta}$, then
$(1-\mu_{\beta}) \cdot \log_{2}\frac{2 (1-\mu_{\beta})}{(1-\mu_{\alpha})+(1-\mu_{\beta})}\leq 0$
and $\nu_{\alpha} \cdot \log_{2}\frac{2 \nu_{\alpha}}{\nu_{\alpha}+\nu_{\beta}}
\leq 0$, and thus $\mathcal{Z}(\alpha, \beta) \leq (1-\mu_{\alpha})  \cdot
\underbrace{\log_{2}\frac{2 (1-\mu_{\alpha})}{(1-\mu_{\alpha})+(1-\mu_{\beta})}}\limits_{\leq 1}
+\nu_{\beta} \cdot \underbrace{\log_{2}\frac{2 \nu_{\beta}}{\nu_{\alpha}+\nu_{\beta}}}\limits_{\leq 1}
\leq (1-\mu_{\alpha})+\nu_{\beta}\leq 2$;

(4) If $\mu_{\alpha}\geq \mu_{\beta}$ and $\nu_{\alpha}\geq \nu_{\beta}$, then
$(1-\mu_{\alpha})  \cdot  \log_{2} \frac{2 (1-\mu_{\alpha})}{(1-\mu_{\alpha})+(1-\mu_{\beta})}\leq 0$
and $\nu_{\beta} \cdot \log_{2}\frac{2 \nu_{\beta}}{\nu_{\alpha}+\nu_{\beta}}
\leq 0$, and thus $\mathcal{Z}(\alpha, \beta) \leq (1-\mu_{\beta}) \cdot
\underbrace{\log_{2}\frac{2 (1-\mu_{\beta})}{(1-\mu_{\alpha})+(1-\mu_{\beta})}}\limits_{\leq 1}
+\nu_{\alpha} \cdot \underbrace{\log_{2}\frac{2 \nu_{\alpha}}{\nu_{\alpha}+\nu_{\beta}}}\limits_{\leq 1}
\leq (1-\mu_{\beta})+\nu_{\alpha} \leq 2$.

Summing up the above shows that
$0\leq \mathcal{Z}(\alpha, \beta)
\leq 2$.
\end{IEEEproof}

\begin{property}
\label{Pro-3-IFV}
$\overline{\bm{JS}}_{_{\mathrm{IF}}}(\alpha, \beta)=1$ if and only if
($\alpha=\langle 0, 1\rangle$ and $\beta=\langle 1, 0 \rangle$) or
($\alpha=\langle 1, 0\rangle$ and $\beta=\langle 0, 1 \rangle$).
\end{property}

\begin{IEEEproof}
{\it Sufficiency.}
By direct calculation and Eq.~\eqref{Eq-Wu-1a-*}, it follows that
$\overline{\bm{JS}}_{_{\mathrm{IF}}}(\langle 0, 1\rangle,
\langle 1, 0\rangle)=\overline{\bm{JS}}_{_{\mathrm{IF}}}(\langle 1, 0\rangle,
\langle 0, 1\rangle)=1$.

{\it Necessity.} Fix $\beta=\langle \mu_{\beta}, \nu_{\beta} \rangle \in
\Theta \backslash \{\langle 0, 1\rangle, \langle 1, 0\rangle\}$.
For any $\alpha \in \Theta \backslash \{\langle 0, 1\rangle, \langle 1, 0\rangle\}$,
according to the proof of Property~\ref{Pro-2-IFV}, consider the following
four cases:

(1) If $\mu_{\alpha}\geq \mu_{\beta}$ and $\nu_{\alpha} \leq \nu_{\beta}$,
by $\nu_{\beta}<1$, we have
$\mathcal{Z}(\alpha, \beta)\leq (1-\mu_{\beta})+\nu_{\beta}
<(1-\mu_{\beta})+1\leq 2$, i.e., $\mathcal{Z}(\alpha, \beta)<2$, and thus
$\overline{\bm{JS}}_{_{\mathrm{IF}}}(\alpha, \beta)=
\sqrt{\frac{\mathcal{Z}(\alpha, \beta)}{2}}<1$;

(2) If $\mu_{\alpha}\leq \mu_{\beta}$ and $\nu_{\alpha} \geq \nu_{\beta}$,
by $\mu_{\beta}<1$, we have
$\mathcal{Z}(\alpha, \beta)\leq \log_{2}\frac{2}{2-\mu_{\beta}}
+\log_{2}\frac{2}{1+\nu_{\beta}}<1+\log_{2}\frac{2}{1+\nu_{\beta}}
\leq 2$, i.e., $\mathcal{Z}(\alpha, \beta)<2$, and thus
$\overline{\bm{JS}}_{_{\mathrm{IF}}}(\alpha, \beta)=
\sqrt{\frac{\mathcal{Z}(\alpha, \beta)}{2}}<1$;

(3) If $\mu_{\alpha}\leq \mu_{\beta}$ and $\nu_{\alpha}\leq \nu_{\beta}$,
by $\nu_{\beta}<1$, we have $\mathcal{Z}(\alpha, \beta)
\leq (1-\mu_{\alpha})+\nu_{\beta}< (1-\mu_{\alpha})+1
\leq 2$, i.e., $\mathcal{Z}(\alpha, \beta)<2$, and thus
$\overline{\bm{JS}}_{_{\mathrm{IF}}}(\alpha, \beta)=
\sqrt{\frac{\mathcal{Z}(\alpha, \beta)}{2}}<1$;

(4) If $\mu_{\alpha}\geq \mu_{\beta}$ and $\nu_{\alpha}\geq \nu_{\beta}$,
by $\nu_{\alpha}<1$, we have $\mathcal{Z}(\alpha, \beta)
\leq (1-\mu_{\beta})+\nu_{\alpha}< (1-\mu_{\beta})+1 \leq 2$,
i.e., $\mathcal{Z}(\alpha, \beta)<2$, and thus
$\overline{\bm{JS}}_{_{\mathrm{IF}}}(\alpha, \beta)=
\sqrt{\frac{\mathcal{Z}(\alpha, \beta)}{2}}<1$.

Together with $\overline{\bm{JS}}_{_{\mathrm{IF}}}
(\langle 0, 1\rangle, \langle 0, 1\rangle)=
\overline{\bm{JS}}_{_{\mathrm{IF}}}
(\langle 1, 0\rangle, \langle 1, 0\rangle)=0$,
from $\overline{\bm{JS}}_{_{\mathrm{IF}}}(\alpha, \beta)=1$,
it follows that $\alpha$, $\beta\in \{\langle 0, 1\rangle,
\langle 1, 0\rangle\}$ and $\alpha \neq \beta$.
\end{IEEEproof}

\begin{property}
\label{Pro-4-IFV}
$\overline{\bm{JS}}_{_{\mathrm{IF}}}(\alpha, \beta)=0$ if and only if $\alpha=\beta$.
\end{property}

\begin{IEEEproof}
{\it Sufficiency.} By Eq.~\eqref{Eq-Wu-1a-*}, it is clear that
$\overline{\bm{JS}}_{_{\mathrm{IF}}}(\alpha, \beta)=0$ if
$\alpha=\beta$.

{\it Necessity.} Assume that $\overline{\bm{JS}}_{_{\mathrm{IF}}}(\alpha, \beta)=0$.
It will be shown that $\alpha=\beta$.

Suppose on the contrary that $\alpha \neq \beta$. Without loss
of generality, assume $0\leq \nu_{\alpha}< \nu_{\beta}$.
From Eq.~\eqref{Z-Function-*}, it follows that
$\mathcal{Z}(\alpha, \beta)
\geq (\nu_{\alpha}+\nu_{\beta})\cdot \zeta (\frac{\nu_{\alpha}}{\nu_{\alpha}
+\nu_{\beta}})\geq \nu_{\beta}\cdot \zeta (\frac{\nu_{\alpha}}{\nu_{\alpha}
+\nu_{\beta}})$.
This, together with the fact that $\zeta(\_)$ is strictly decreasing on $[0, 0.5]$ and
$\frac{\nu_{\alpha}}{\nu_{\alpha}+\nu_{\beta}}\in [0, 0.5)$, implies that
$\mathcal{Z}(\alpha, \beta)\geq \nu_{\beta}\cdot \zeta (\frac{\nu_{\alpha}}{\nu_{\alpha}
+\nu_{\beta}})>\nu_{\beta} \cdot \zeta(0.5)=0$. Therefore,
$$
0=\overline{\bm{JS}}_{_{\mathrm{IF}}}(\alpha, \beta)=
\sqrt{\frac{\mathcal{Z}(\alpha, \beta)}{2}}>0,
$$
which is a contradiction.
\end{IEEEproof}

\begin{property}
\label{Pro-5-IFV}
Let $\alpha$, $\beta$, $\gamma \in \Theta$.

(1) If $\alpha \subset \beta \subset \gamma$, then
$\overline{\bm{JS}}_{_{\mathrm{IF}}}(\alpha, \beta) \leq
\overline{\bm{JS}}_{_{\mathrm{IF}}}(\alpha, \gamma)$
and $\overline{\bm{JS}}_{_{\mathrm{IF}}}(\beta, \gamma)\leq
\overline{\bm{JS}}_{_{\mathrm{IF}}}(\alpha, \gamma)$;

(2) If $\alpha\subsetneqq \beta \subsetneqq \gamma$, then
$\overline{\bm{JS}}_{_{\mathrm{IF}}}(\alpha, \beta) <
\overline{\bm{JS}}_{_{\mathrm{IF}}}(\alpha, \gamma)$
and $\overline{\bm{JS}}_{_{\mathrm{IF}}}(\beta, \gamma)<
\overline{\bm{JS}}_{_{\mathrm{IF}}}(\alpha, \gamma)$.
\end{property}

\begin{IEEEproof}
(1) Fix the IFV $\alpha$. It follows directly from Lemma~\ref{Lemma-1-*}
that $\mathcal{Z}(\alpha, \beta) \leq \mathcal{Z}(\alpha, \gamma)$.
Thus, $\overline{\bm{JS}}_{_{\mathrm{IF}}}(\alpha, \beta)
=\sqrt{\frac{1}{2}\mathcal{Z}(\alpha, \beta)}\leq
\sqrt{\frac{1}{2}\mathcal{Z}(\alpha, \gamma)}=
\overline{\bm{JS}}_{_{\mathrm{IF}}}(\alpha, \gamma)$. Similarly, we can prove
$\overline{\bm{JS}}_{_{\mathrm{IF}}}(\beta, \gamma)\leq
\overline{\bm{JS}}_{_{\mathrm{IF}}}(\alpha, \gamma)$
by Lemma~\ref{Lemma-2-*}.

(2) Similarly, we can prove
$\overline{\bm{JS}}_{_{\mathrm{IF}}}(\alpha, \beta) <
\overline{\bm{JS}}_{_{\mathrm{IF}}}(\alpha, \gamma)$
and $\overline{\bm{JS}}_{_{\mathrm{IF}}}(\beta, \gamma)<
\overline{\bm{JS}}_{_{\mathrm{IF}}}(\alpha, \gamma)$
by Lemmas~\ref{Lemma-1-*-strict} and \ref{Lemma-2-*-strict}.
\end{IEEEproof}

\begin{property}[\textrm{Triangle inequality}]
\label{Pro-6-IFV}
Let $\alpha$, $\beta$, $\gamma\in \Theta$. Then,
$\overline{\bm{JS}}_{_{\mathrm{IF}}}(\alpha, \beta)+
\overline{\bm{JS}}_{_{\mathrm{IF}}}(\beta, \gamma)\geq
\overline{\bm{JS}}_{_{\mathrm{IF}}}(\alpha, \gamma)$.
\end{property}

\begin{IEEEproof}
For convenience, denote $\sqrt{L(1-\mu_{\alpha}, 1-\mu_{\beta})}=\xi_1$,
$\sqrt{L(\nu_{\alpha}, \nu_{\beta})}=\eta_1$, $\sqrt{L(1-\mu_{\beta}, 1-\mu_{\gamma})}=\xi_2$,
$\sqrt{L(\nu_{\beta}, \nu_{\gamma})}$ $=\eta_2$, $\sqrt{L(1-\mu_{\alpha}, 1-\mu_{\gamma})}=\xi_3$,
 and $\sqrt{L(\nu_{\alpha}, \nu_{\gamma})}=\eta_3$.
By Lemma~\ref{Lemma-ES}, we have
\begin{equation}
\label{Eq-Wu-4-**}
\xi_3\leq \xi_1+\xi_2, \
\eta_3\leq \eta_1+\eta_2.
\end{equation}
Meanwhile, from $(\xi_1^2+\eta_1^2)\cdot (\xi_2^2+\eta_2^2)-(\xi_1\xi_2+ \eta_1\eta_2)^{2}
=\xi_1^2 \eta_2^2+$ $\xi_2^2\eta_1^2-2\xi_1\xi_2\eta_1\eta_2\geq 0$, it follows that
$(\sqrt{\xi_1^{2}+\eta_1^{2}}+ \sqrt{\xi_2^{2}
+\eta_2^{2}})^{2}=\xi_1^{2}+\eta_1^{2}+\xi_2^{2}+\eta_2^{2}+2 \sqrt{\xi_1^{2}+\eta_1^{2}}
\sqrt{\xi_2^{2}+\eta_2^{2}}\geq \xi_1^{2}+\eta_1^{2}+\xi_2^{2}+\eta_2^{2}+ 2
(\xi_1\xi_2+\eta_1\eta_2)=(\xi_1+\xi_2)^{2}+(\eta_1+\eta_2)^{2}$. This, together with
Eq.~\eqref{Eq-Wu-4-**}, implies that
\begin{align*}
\overline{\bm{JS}}_{_{\mathrm{IF}}}(\alpha, \gamma)=& \sqrt{\frac{1}{2}
(\xi_3^2+\eta_3^{2})}\leq \sqrt{\frac{1}{2}[(\xi_1+\xi_2)^{2}+(\eta_1+\eta_2)^{2}]}\\
\leq & \sqrt{\frac{1}{2}\left(\sqrt{\xi_1^{2}+\eta_1^{2}}+ \sqrt{\xi_2^{2}
+\eta_2^{2}}\right)^{2}}\\
=& \sqrt{\frac{1}{2}(\xi_1^2+\eta_1^2)}+
\sqrt{\frac{1}{2}(\xi_2^2+\eta_2^2)}\\
=& \overline{\bm{JS}}_{_{\mathrm{IF}}}(\alpha, \beta)+
\overline{\bm{JS}}_{_{\mathrm{IF}}}(\beta, \gamma).
\end{align*}
\end{IEEEproof}

Summing Properties~\ref{Pro-1-IFV}--\ref{Pro-6-IFV}, we have the following results.

\begin{theorem}
(1) The divergence measure $\overline{\bm{JS}}_{_{\mathrm{IF}}}$ is a
SIFDisM on $\Theta$.

(2) The function $\mathbf{S}_{_{\mathrm{IF}}}(\alpha, \beta)=1-\overline{\bm{JS}}_{_{\mathrm{IF}}}
(\alpha, \beta)$ is a SIFSimM on $\Theta$.
\end{theorem}

\begin{theorem}
\label{Entropy-Thm-IFV}
The mapping $E$ defined by
\begin{equation}
\label{Entropy-eq-1}
\begin{split}
E: \Theta & \longrightarrow [0, 1], \\
\alpha & \longmapsto 1-\overline{\bm{JS}}_{_{\mathrm{IF}}}
(\alpha, \alpha^{\complement}),
\end{split}
\end{equation}
is an IFEM on $\Theta$.
\end{theorem}

\begin{IEEEproof}
(E1), (E2), and (E3) follow directly from Properties \ref{Pro-3-IFV} and \ref{Pro-4-IFV},
and Eq.~\eqref{Entropy-eq-1}, respectively.

(E4) For $\alpha$, $\beta\in \Theta$, consider the following two cases:

 E4-1) If $\mu_{\alpha}\leq \mu_{\beta} \leq \nu_{\beta} \leq \nu_{\alpha}$, then
 $\alpha\subset \beta\subset \beta^{\complement} \subset \alpha^{\complement}$. This,
 together with Property~\ref{Pro-5-IFV}, implies that $E(\alpha)=
 1-\overline{\bm{JS}}_{_{\mathrm{IF}}}(\alpha, \alpha^{\complement})
 \leq 1-\overline{\bm{JS}}_{_{\mathrm{IF}}}(\alpha, \beta^{\complement})
 \leq 1-\overline{\bm{JS}}_{_{\mathrm{IF}}}(\beta, \beta^{\complement})
 =E(\beta)$;

 E4-2) If $\mu_{\alpha}\geq \mu_{\beta} \geq \nu_{\beta} \geq \nu_{\alpha}$, then
 $\alpha^{\complement} \subset \beta^{\complement} \subset \beta \subset \alpha$.
 This,
 together with Property~\ref{Pro-5-IFV}, implies that $E(\alpha)=
 1-\overline{\bm{JS}}_{_{\mathrm{IF}}}(\alpha, \alpha^{\complement})
 \leq 1-\overline{\bm{JS}}_{_{\mathrm{IF}}}(\alpha, \beta^{\complement})
 \leq 1-\overline{\bm{JS}}_{_{\mathrm{IF}}}(\beta, \beta^{\complement})
 =E(\beta)$.
\end{IEEEproof}

\begin{figure}[H]
\centering
{\scalebox{0.33}{\includegraphics[]{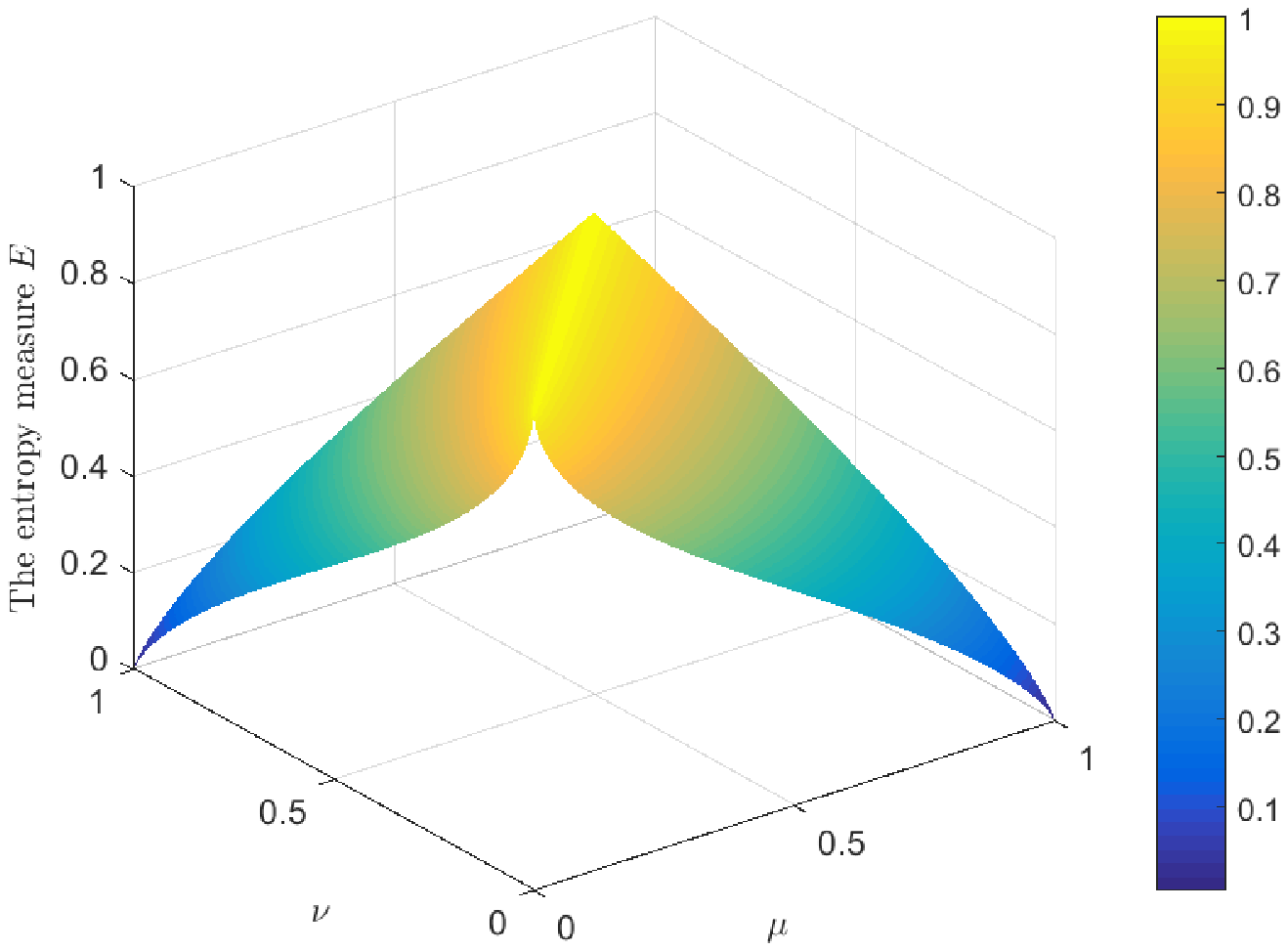}}}
\caption{The graph of IFEM $E$ in Theorem~\ref{Entropy-Thm-IFV}}
\label{Fig-Entropy-IFV}
\end{figure}

Fig.~\ref{Fig-Entropy-IFV} shows the graph of IFEM
in Theorem~\ref{Entropy-Thm-IFV}.

\subsection{A distance/similarity measure on IFSs}

Let $X=\{x_1, x_2, \ldots, x_n\}$ be a finite UOD and $I_1=\Big\{\frac{\alpha_j^{(1)}}{x_j}
\mid 1\leq j\leq n,  \alpha_j^{(1)}\in \Theta\Big\}$ and $I_2=\Big\{\frac{\alpha_j^{(2)}}{x_j}
\mid 1\leq j\leq n, \alpha_j^{(2)}\in \Theta\Big\}$ be two IFSs on $X$.

Applying the the normalized Jensen-Shannon IF divergence measure
$\overline{\bm{JS}}_{_{\mathrm{IF}}}$ on IFVs, we define a new normalized
distance measure $\bm{d}_{_{\mathrm{Wu}}}(I_1, I_2)$ as follows:
\begin{small}
\begin{equation}
\label{Eq-Wu-1a}
\begin{split}
& \bm{d}_{_{\mathrm{Wu}}}(I_1, I_2)=
\sum_{j=1}^{n}\omega_j\cdot \overline{\bm{JS}}_{_{\mathrm{IF}}}(\alpha_{j}^{(1)},
\alpha_{j}^{(2)})\\
=& \sum_{j=1}^{n}
\omega_j \left[ \frac{1}{2}\left( (1-\mu_{I_1}(x_j))  \cdot  \log_{2}
\frac{2 (1-\mu_{I_1}(x_j))}{(1-\mu_{I_1}(x_j))+(1-\mu_{I_2}(x_j))}\right.\right.\\
& \quad \quad \quad +(1-\mu_{I_2}(x_j)) \cdot \log_{2}\frac{2 (1-\mu_{I_2}(x_j))}
{(1-\mu_{I_1}(x_j))+(1-\mu_{I_2}(x_j))}\\
& \quad \quad \quad +\nu_{I_1}(x_j) \cdot \log_{2}\frac{2 \nu_{I_1}(x_j)}{\nu_{I_1}(x_j)+\nu_{I_2}(x_j)}\\
& \left.\left. \quad \quad \quad +\nu_{I_2}(x_j) \cdot \log_{2}\frac{2 \nu_{I_2}(x_j)}{\nu_{I_1}(x_j)+\nu_{I_2}(x_j)}\right)\right]^{0.5}.
\end{split}
\end{equation}
\end{small}
where $\alpha_j^{(1)}=\langle \mu_{I_1}(x_j),
\nu_{I_1}(x_j)\rangle$, $\alpha_j^{(2)}=\langle \mu_{I_2}(x_j),
\nu_{I_2}(x_j) \rangle$, and $\omega=(\omega_1, \omega_2, \ldots, \omega_n)^{\top}$ is the weight
vector of $x_{j}$ ($j=1, 2, \ldots, n$) with $\omega_j\in (0, 1]$ and
$\sum_{j=1}^{n}\omega_j=1$. Clearly,
\begin{equation}
\label{eq-Wu=}
\bm{d}_{_{\mathrm{Wu}}}(I_1, I_2)=
\sum_{j=1}^{n}\omega_j \cdot \sqrt{\frac{\mathcal{Z}(\alpha_{j}^{(1)}, \alpha_{j}^{(2)})}{2}}.
\end{equation}

By applying Properties~\ref{Pro-1-IFV}--\ref{Pro-6-IFV},
one can easily verify that the function $\bm{d}_{_{\mathrm{Wu}}}$
has the following properties.

\begin{property}
\label{Pro-1-IFV-2}
$\bm{d}_{_{\mathrm{Wu}}}(I_1, I_2)=\bm{d}_{_{\mathrm{Wu}}}(I_2, I_1)$.
\end{property}

\begin{property}
$0\leq \bm{d}_{_{\mathrm{Wu}}}(I_1, I_2)\leq 1$.
\end{property}

\begin{property}
\label{Pro-3-IFV-2}
$\bm{d}_{_{\mathrm{Wu}}}(I_1, I_2)=1$ if and only if, for any $1\leq j\leq n$,
($\alpha_{j}^{(1)}=\langle 0, 1\rangle$ and $\alpha_{j}^{(2)}=\langle 1, 0 \rangle$) or
($\alpha_{j}^{(1)}=\langle 1, 0\rangle$ and $\alpha_{j}^{(2)}=\langle 0, 1 \rangle$).
\end{property}

\begin{property}
$\bm{d}_{_{\mathrm{Wu}}}(I_1, I_2)=0$ if and only if $I_1=I_2$.
\end{property}

\begin{property}
\label{Pro-6-IFV-2}
Let $X=\{x_1, x_2, \ldots, x_3\}$ and $I_1$, $I_2$, $I_3\in \mathrm{IFS}(X)$.

(1) If $I_1 \subset I_2 \subset I_3$, then
$\bm{d}_{_{\mathrm{Wu}}}(I_1, I_2) \leq \bm{d}_{_{\mathrm{Wu}}}(I_1, I_3)$
and $\bm{d}_{_{\mathrm{Wu}}}(I_2, I_3)\leq \bm{d}_{_{\mathrm{Wu}}}(I_1, I_3)$.

(2) (1) If $I_1\subsetneqq I_2\subsetneqq I_3$, then
$\bm{d}_{_{\mathrm{Wu}}}(I_1, I_2) < \bm{d}_{_{\mathrm{Wu}}}(I_1, I_3)$
and $\bm{d}_{_{\mathrm{Wu}}}(I_2, I_3)< \bm{d}_{_{\mathrm{Wu}}}(I_1, I_3)$.
\end{property}

Summing Properties~\ref{Pro-1-IFV-2}--\ref{Pro-6-IFV-2}, similarly to the proof
of Theorem~\ref{Entropy-Thm-IFV}, we have the following results.

\begin{theorem}
(1) The distance measure $\bm{d}_{_{\mathrm{Wu}}}$ defined by
Eq.~\eqref{Eq-Wu-1a} is a SIFDisM on $\mathrm{IFS}(X)$.

(2) The function $\mathbf{S}_{_{\mathrm{Wu}}}(I_1, I_2)
=1-\bm{d}_{_{\mathrm{Wu}}}(I_1, I_2)$ is a SIFSimM on $\mathrm{IFS}(X)$.
\end{theorem}

\begin{theorem}
\label{Entropy-Thm-IFS}
The mapping $E$ defined by
\begin{equation}
\label{Entropy-eq-2}
\begin{split}
E: \Theta & \longrightarrow [0, 1], \\
\alpha & \longmapsto 1-\bm{d}_{_{\mathrm{Wu}}}
(\alpha, \alpha^{\complement}),
\end{split}
\end{equation}
is an IFEM on $\mathrm{IFS}(X)$.
\end{theorem}

\section{Comparative analysis}
Xiao~\cite{Xiao2021} showed that the distance measure $\bm{d}_{\widetilde{\chi}}$
is better than other existing distance measures proposed in
\cite{SK2000,Gr2004,WX2005,HPK2012,YC2012,SMLZC2018,SWQH2019} by some
numerical comparisons. However, for the nonlinear distance measure $d_{_{\mathrm{YC}}}$
introduced by Yang and Chiclana~\cite{YC2009}, only one figure
(see \cite[Fig.~4]{Xiao2021}) was used to show that the curve of $\bm{d}_{\widetilde{\chi}}$
is sharper than that of $d_{_{\mathrm{YC}}}$ for some special cases.
This does not convincingly explain the superiority of the distance
measure $\bm{d}_{\widetilde{\chi}}$. This section demonstrates
that our proposed distance measure $\bm{d}_{_{\mathrm{Wu}}}$
is completely superior to Xiao's distance measure $d_{\widetilde{\chi}}$
and Yang and Chiclana's distance measure $d_{_{\mathrm{YC}}}$. Because of
the duality of distance and similarity measures, for brevity we only compare and
analyze distance measures.

Xiao~\cite{Xiao2021} used some examples to illustrate the superiority
of the distance measure $\bm{d}_{\widetilde{\chi}}$. The following shows
that our proposed distance measure has the same superiority for the same
numerical examples.

\begin{example}[{\textrm{\protect\cite[Examples~3]{Xiao2021}}}]
\label{Exm-3-Xiao}
Assume that the IFSs $I_1$ and $I_2$ on UOD $X=\{x\}$ are given by
$I_1=\left\{\frac{\langle \mu, \nu\rangle}{x}\right\}$ and
$I_2=\left\{\frac{\langle \nu, \mu\rangle}{x}\right\}.$

Fig.~\ref{Fig-1} shows the changing trend of distance
$\bm{d}_{_{\mathrm{Wu}}}(I_1, I_2)$ with varying parameters
$\mu$ and $\nu$ satisfying $\langle \mu,
\nu \rangle \in \Theta$. Observing from Fig.~\ref{Fig-1}, it can be seen
that our proposed distance measure $\bm{d}_{_{\mathrm{Wu}}}$ has
a form similar to the distance measure $\bm{d}_{\widetilde{\chi}}$
shown in~\cite[Fig.~2]{Xiao2021} with $\overline{\bm{JS}}_{_{\mathrm{IF}}}
(\langle 0, 1 \rangle, \langle 1, 0 \rangle)
=\overline{\bm{JS}}_{_{\mathrm{IF}}}(\langle 1, 0 \rangle, \langle 0, 1 \rangle)
=1$.
\begin{figure}[H]
\centering
{\scalebox{0.33}{\includegraphics[]{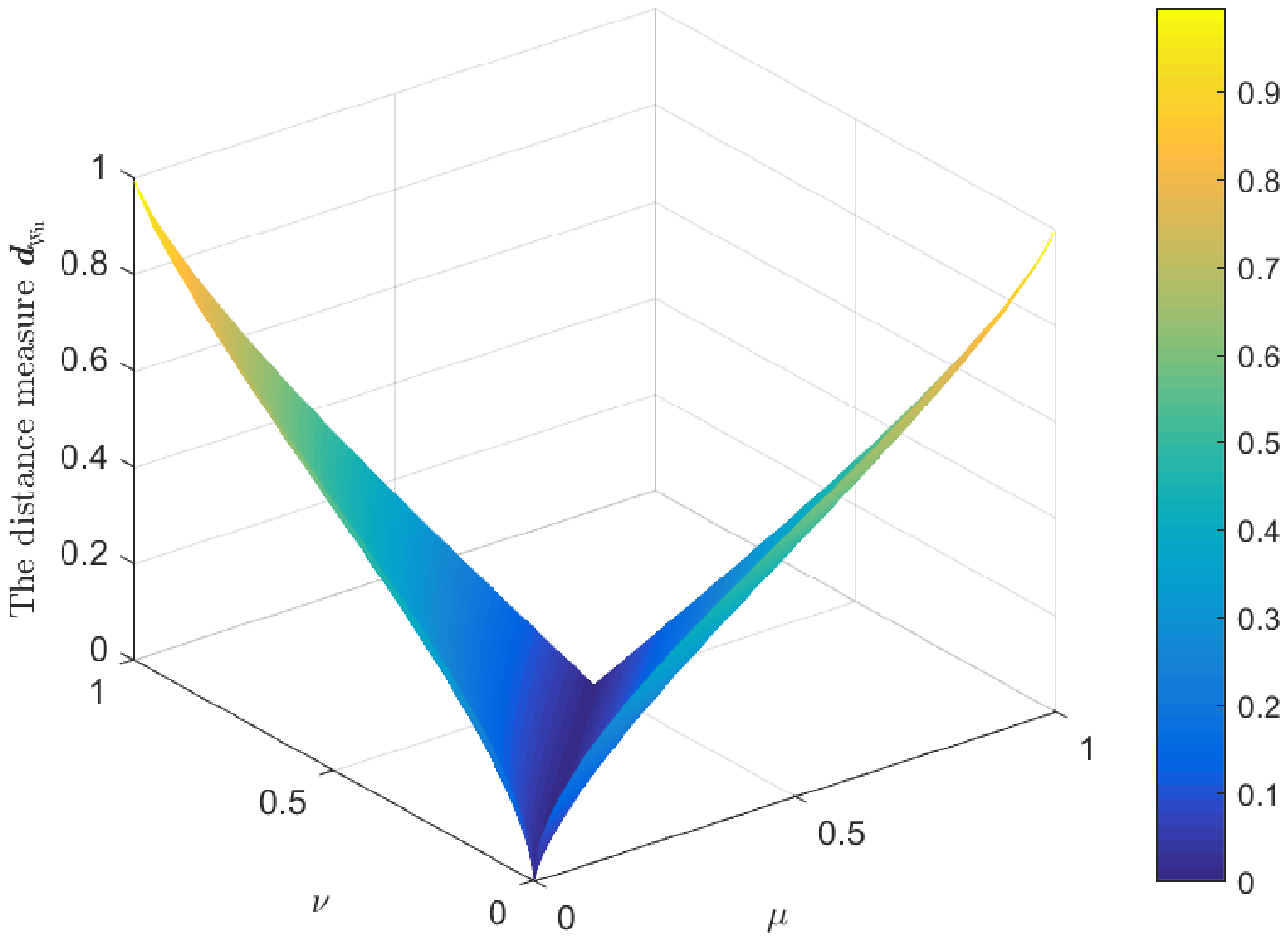}}}
\caption{The distance measure $\bm{d}_{_{\mathrm{Wu}}}(I_1, I_2)$ in Example~\ref{Exm-3-Xiao}}
\label{Fig-1}
\end{figure}
\end{example}


\begin{example}[{\textrm{\protect\cite[Examples~5 and 6]{Xiao2021}}}]
\label{Exm-5-Xiao}
Assume that the IFSs $A_i$ and $B_i$ on the UOD $X=\{x_1, x_2\}$
in Case~$i$ ($i=1, 2, 3, 4, 5$) are given as shown in Table~\ref{Tab-1}.
\begin{table}[H]	
	\centering
	\caption{Two IFSs $A_i$ and $B_{i}$ in Case~$i$
($i=1, 2, 3, 4, 5$) in Example~\ref{Exm-5-Xiao}}
	\label{Tab-1}
     \scalebox{1.0}{
	\begin{tabular}{cccccc}
		\toprule
		IFSs & Case~1 &  Case~2 \\
		\midrule
		$A_{i}$ & $\left\{\frac{\langle 0.30, 0.20\rangle}{x_1},
\frac{\langle 0.40, 0.30 \rangle}{x_2}\right\}$ &
$\left\{\frac{\langle 0.30, 0.20\rangle}{x_1},
\frac{\langle 0.40, 0.30 \rangle}{x_2}\right\}$ \\
$B_{i}$ & $\left\{\frac{\langle 0.15, 0.25\rangle}{x_1},
\frac{\langle 0.25, 0.35 \rangle}{x_2}\right\}$ &
$\left\{\frac{\langle 0.16, 0.26\rangle}{x_1},
\frac{\langle 0.26, 0.36 \rangle}{x_2}\right\}$ \\
		\bottomrule
		IFSs & Case~3 &  Case~4 \\
		\midrule
		$A_{i}$ & $\left\{\frac{\langle 0.50, 0.40\rangle}{x_1},
\frac{\langle 0.40, 0.30 \rangle}{x_2}\right\}$ &
$\left\{\frac{\langle 0.50, 0.40\rangle}{x_1},
\frac{\langle 0.40, 0.30 \rangle}{x_2}\right\}$ \\
$B_{i}$ & $\left\{\frac{\langle 0.15, 0.25\rangle}{x_1},
\frac{\langle 0.25, 0.35 \rangle}{x_2}\right\}$ &
$\left\{\frac{\langle 0.16, 0.26\rangle}{x_1},
\frac{\langle 0.26, 0.36 \rangle}{x_2}\right\}$ \\
		\bottomrule
		IFSs & Case~5  \\
		\midrule
		$A_{i}$ & $\left\{\frac{\langle 0.30, 0.20\rangle}{x_1},
\frac{\langle 0.40, 0.30 \rangle}{x_2}\right\}$ \\
$B_{i}$ & $\left\{\frac{\langle 0.45, 0.15\rangle}{x_1},
\frac{\langle 0.55, 0.25 \rangle}{x_2}\right\}$ \\
		\bottomrule
	\end{tabular}
      }
\end{table}

The comparative results produced by Xiao's distance measure $\bm{d}_{\widetilde{\chi}}$
and our distance measure $\bm{d}_{_{\mathrm{Wu}}}$ are displayed in Table~\ref{Tab-2},
which indicates that both distance measures can effectively distinguish $A_{i}$
and $B_{i}$ in Cases~1--5.

\begin{table}[H]	
	\centering
	\caption{Two distance measures in Cases 1-5}
	\label{Tab-2}
     \scalebox{1.0}{
	\begin{tabular}{cccccc}
		\toprule
		Distance & Case~1 &  Case~2 & Case~3 & Case~4 & Case~5 \\
		\midrule
        $\bm{d}_{_{\mathrm{\widetilde{\chi}}}}$ & $0.14614$ & $0.13531$ &
           $0.17210$ & $0.13352$ & $0.13224$ \\
		$\bm{d}_{_{\mathrm{Wu}}}$ & $0.08563$ & $0.08568$ &
           $0.07462$ & $0.09802$ & $0.09615$ \\
		\bottomrule
	\end{tabular}
      }
\end{table}
\end{example}

\subsection{Comparative analysis between Xiao's distance
measure $d_{\widetilde{\chi}}$ and our proposed distance
$d_{_{\mathrm{Wu}}}$}

\begin{example}[\textrm{Continuation of Example~\protect\ref{Exm-Wu-2-Sec2}}]
\label{Exm-Wu-2}
Let IFSs $I_1$, $I_1^{\prime}$, $I_2^{(\lambda)}$, and $I_3^{(\lambda)}$ on UOD $X=\{x\}$
be given as in Example~\ref{Exm-Wu-2-Sec2}.
By direct calculation, it follows from Eq.~\eqref{Eq-Xiao-1}
that
{\small
\begin{equation}
\begin{split}
&\bm{d}_{\widetilde{\chi}}(I_1^{\prime}, I_3^{(\lambda)})\\
=& \sqrt{\frac{1}{2}\left[\lambda \log_{2}\frac{2\lambda}{\lambda}
+1\cdot \log_{2} \frac{2}{2-\lambda}+(1-\lambda)\log_{2}\frac{2(1-\lambda)}{2-\lambda}\right]}\\
=& \sqrt{\frac{1}{2}\left[\lambda+ \log_{2} \frac{2}{2-\lambda}+(1-\lambda)\log_{2}\frac{2(1-\lambda)}{2-\lambda}\right]}.
\end{split}
\end{equation}}

(iii) Observing from Example~\ref{Exm-Wu-2-Sec2} (ii), in contrast to
Property~\ref{Pro-3-IFV}, there exist infinite IFSs $I_{2}^{(\lambda)}$
($\lambda\in [0, 1]$), such that the distance from $I_1^{\prime}$ is equal
to the maximum value $1$. This result is completely inferior to
Property~\ref{Pro-3-IFV} for IFVs.

By applying
our distance measure $\bm{d}_{_{\mathrm{Wu}}}$ defined in Eq.~\eqref{Eq-Wu-1a},
we obtain
{\small\begin{equation}
\label{Eq-Wu-Exm8-1}
\begin{split}
\bm{d}_{_{\mathrm{Wu}}}(I_1, I_2^{(\lambda)})
= \sqrt{\frac{1}{2}\cdot (1-\lambda)\cdot \log_{2}\frac{2(1-\lambda)}{1-\lambda}}
=\sqrt{\frac{1-\lambda}{2}},
\end{split}
\end{equation}
\begin{equation}
\label{Eq-Wu-Exm8-2}
\begin{split}
&\bm{d}_{_{\mathrm{Wu}}}(I_1^{\prime}, I_2^{(\lambda)})\\
= & \sqrt{\frac{1}{2}\left[\log_{2}\frac{2}{2-\lambda}+(1-\lambda)\cdot \log_{2}
\frac{2(1-\lambda)}{2-\lambda}+(1-\lambda)+1\right]},
\end{split}
\end{equation}}
and
{\small \begin{equation}
\label{Eq-Wu-Exm8-3}
\begin{split}
& \bm{d}_{_{\mathrm{Wu}}}(I_1, I_3^{(\lambda)})\\
= & \sqrt{\frac{1}{2}\left[(1-\lambda)\cdot \log_{2}\frac{2(1-\lambda)}{1-\lambda}
+(1-\lambda)\cdot \log_{2}\frac{2(1-\lambda)}{1-\lambda}\right]}\\
= & \sqrt{1-\lambda},
\end{split}
\end{equation}
\begin{equation}
\label{Eq-Wu-Exm8-4}
\begin{split}
\bm{d}_{_{\mathrm{Wu}}}(I_1^{\prime}, I_3^{(\lambda)})
= \sqrt{\log_{2}\frac{2}{2-\lambda}+(1-\lambda)\cdot \log_{2}
\frac{2(1-\lambda)}{2-\lambda}}.
\end{split}
\end{equation}}

Directly observing the values of $I_1$, $I_1^{\prime}$,
$I_2^{(\lambda)}$, and $I_3^{(\lambda)}$, we find that
$I_1^{\prime}\subset I_3^{(\lambda)}\subset I_2^{(\lambda)} \subset I_1$, and thus
$I_2^{(\lambda)}$ is more similar to $I_1$ than $I_{3}^{(\lambda)}$, and
$I_3^{(\lambda)}$ is more similar to $I_1^{\prime}$ than $I_{2}^{(\lambda)}$,
which are consistent with our computed results, since $\mathbf{S}_{_{\mathrm{Wu}}}
(I_1, I_3^{(\lambda)})
=1-\bm{d}_{_{\mathrm{Wu}}}(I_1, I_3^{(\lambda)})=1-\sqrt{1-\lambda}
<1-\sqrt{\frac{1-\lambda}{2}}=1-\bm{d}_{_{\mathrm{Wu}}}(I_1, I_2^{(\lambda)})
=\mathbf{S}_{_{\mathrm{Wu}}}(I_1, I_2^{(\lambda)})$ and
$\mathbf{S}_{_{\mathrm{Wu}}}
(I_1^{\prime}, I_3^{(\lambda)})
=1-\bm{d}_{_{\mathrm{Wu}}}(I_1^{\prime}, I_3^{(\lambda)})
=1-\sqrt{\log_{2}\frac{2}{2-\lambda}+(1-\lambda)\cdot \log_{2}
\frac{2(1-\lambda)}{2-\lambda}}
>1-\sqrt{\frac{1}{2}\left[\log_{2}\frac{2}{2-\lambda}+(1-\lambda)\cdot \log_{2}
\frac{2(1-\lambda)}{2-\lambda}+(1-\lambda)+1\right]}
=1-\bm{d}_{_{\mathrm{Wu}}}(I_1^{\prime}, I_2^{(\lambda)})
=\mathbf{S}_{_{\mathrm{Wu}}}(I_1^{\prime}, I_2^{(\lambda)})$.

By varying the parameter $\lambda$ from $0$ to $1$,
the results in Figs.~\ref{Fig-1-1} and \ref{Fig-1-2} visualize
the changing trend of distances between $I_{2}^{(\lambda)}$
and $I_1$ ($I_{1}^{\prime}$) and between $I_{3}^{(\lambda)}$ and $I_1$ ($I_{1}^{\prime}$)
discussed in Example~\ref{Exm-Wu-2} by using Xiao's distance measure $\bm{d}_{\widetilde{\chi}}$
and our proposed distance measure $\bm{d}_{_{\mathrm{Wu}}}$. The simulation results are
consistent with our calculation. This example indicates that our proposed
distance is far superior to Xiao's distance measure $\bm{d}_{\widetilde{\chi}}$ in~\cite{Xiao2021}.
\begin{figure}[H]
\centering
{\scalebox{0.33}{\includegraphics[]{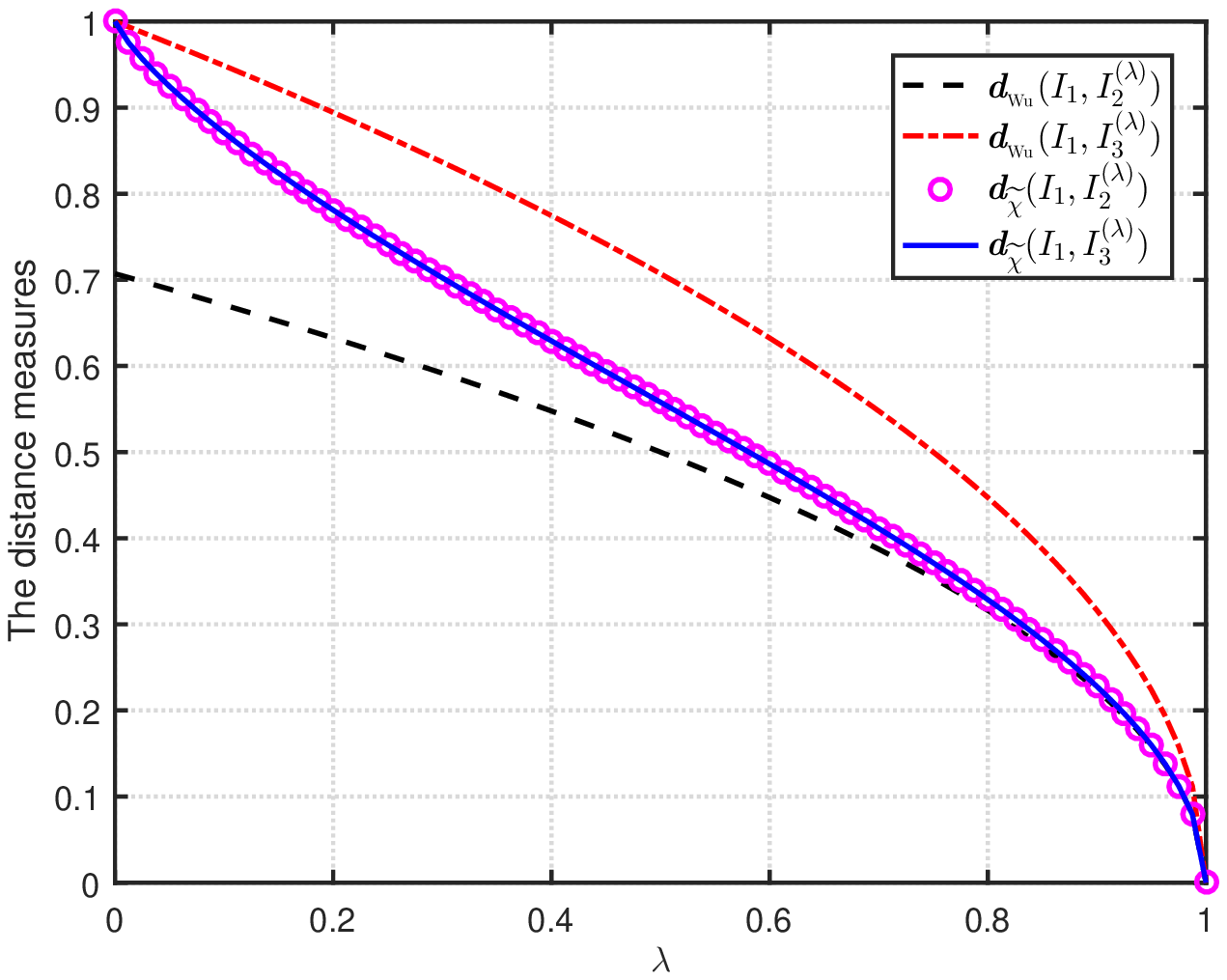}}}
\caption{The distances between $I_{2}^{(\lambda)}$
and $I_1$, and between $I_{3}^{(\lambda)}$ and $I_1$
in Example~\ref{Exm-Wu-2}}
\label{Fig-1-1}
\end{figure}
\begin{figure}[H]
\centering
{\scalebox{0.33}{\includegraphics[]{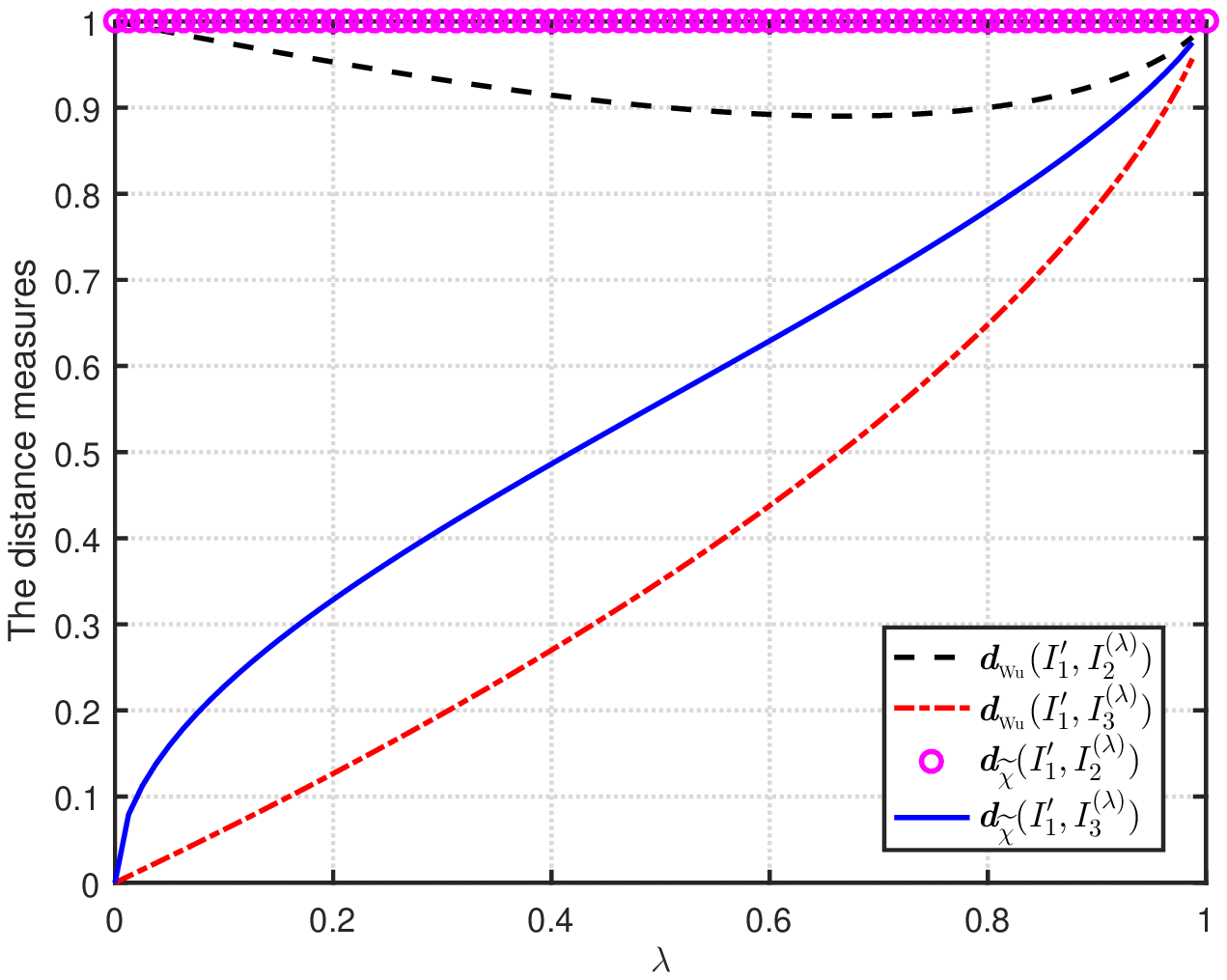}}}
\caption{The distances between $I_{2}^{(\lambda)}$
and $I_1^{\prime}$, and between $I_{3}^{(\lambda)}$ and $I_1^{\prime}$
in Example~\ref{Exm-Wu-2}}
\label{Fig-1-2}
\end{figure}
\end{example}

\begin{example}
\label{Exm-Wu-3-Sec5}
Figs.~\ref{Fig-Xiao} and \ref{Fig-Wu} visualize the changing trend of
distances between $I_{2}=\left\{\frac{\langle \mu, \nu \rangle}{x}\right\}$
and $I_1=\left\{\frac{\langle 1, 0 \rangle}{x}\right\}$, and
between $I_{2}=\left\{\frac{\langle \mu, \nu \rangle}{x}\right\}$
and $I_1^{\prime}=\left\{\frac{\langle 0, 1 \rangle}{x}\right\}$
with varying the parameters $\mu$ and $\nu$ satisfying $\langle \mu,
\nu \rangle \in \Theta$, by using Xiao's distance measure
$\bm{d}_{\widetilde{\chi}}$ and our proposed distance measure
$\bm{d}_{_{\mathrm{Wu}}}$, respectively.

From Fig.~\ref{Fig-Xiao} (e) and (f), we can observe that
(1) the distance $\bm{d}_{\widetilde{\chi}}(I_1, I_2)$ between
$I_1=\left\{\frac{\langle 1, 0 \rangle}{x}\right\}$ and
$I_2=\left\{\frac{\langle \mu, \nu \rangle}{x}\right\}$
remains unchanged when the membership degree $\mu$ is fixed;
(2) the distance $\bm{d}_{\widetilde{\chi}}(I_1^{\prime}, I_2)$ between
$I_1^{\prime}=\left\{\frac{\langle 0, 1 \rangle}{x}\right\}$ and
$I_2=\left\{\frac{\langle \mu, \nu \rangle}{x}\right\}$
remains unchanged when the non-membership degree $\nu$ is fixed.
These are consistent with the following calculation results:
\begin{align*}
\bm{d}_{\widetilde{\chi}}(I_1, I_2)=\sqrt{\frac{1}{2}\left[\log_{2}\frac{2}{1+\mu}
+\mu\cdot \log_{2}\frac{2\mu}{1+\mu}+(1-\mu)\right]},
\end{align*}
and
\begin{align*}
\bm{d}_{\widetilde{\chi}}(I_1^{\prime}, I_2)= \sqrt{\frac{1}{2}\left[\log_{2}\frac{2}{1+\nu}
+\nu\cdot \log_{2}\frac{2\nu}{1+\nu}+(1-\nu)\right]}.
\end{align*}

From Figs.~\ref{Fig-Wu} (c) and (d), we can observe that
(1) the distance $\bm{d}_{_{\mathrm{Wu}}}(I_1, I_2)$ between
$I_1=\left\{\frac{\langle 1, 0 \rangle}{x}\right\}$ and
$I_2=\left\{\frac{\langle \mu, \nu \rangle}{x}\right\}$
increases strictly with the increase of $\nu \in [0, 1-\mu]$
when the membership degree $\mu$ is fixed;
(2) the distance $\bm{d}_{_{\mathrm{Wu}}}(I_1^{\prime}, I_2)$ between
$I_1^{\prime}=\left\{\frac{\langle 0, 1 \rangle}{x}\right\}$ and
$I_2=\left\{\frac{\langle \mu, \nu \rangle}{x}\right\}$
increases strictly with the increase of $\mu \in [0, 1-\nu]$
when the membership degree $\nu$ is fixed. These are reasonable
and consistent
with our results in Lemmas~\ref{Lemma-1-*}--\ref{Lemma-2-*-strict}.

The above results also show the unreasonableness of Xiao's distance measure
$\bm{d}_{\widetilde{\chi}}$ and indicate the superiority of
our proposed distance measure.
\begin{figure}[H]
\centering
\subfigure[The distance $\bm{d}_{\widetilde{\chi}}(I_1, I_2)$ in Example~\ref{Exm-Wu-3-Sec5}]
{\scalebox{0.295}{\includegraphics[]{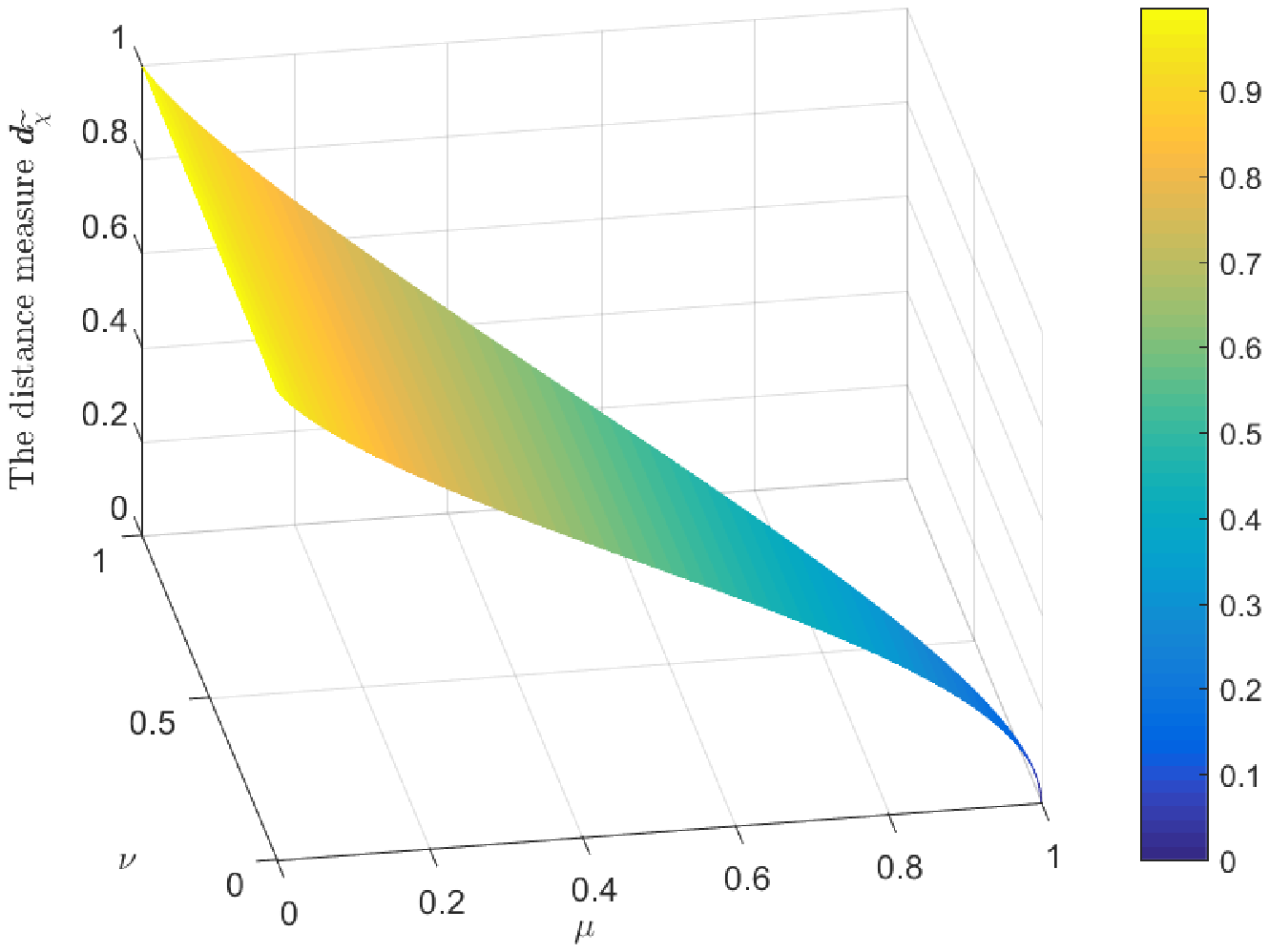}}}
\subfigure[The distance $\bm{d}_{\widetilde{\chi}}(I_1^{\prime}, I_2)$
in Example~\ref{Exm-Wu-3-Sec5}]{\scalebox{0.295}{\includegraphics[]{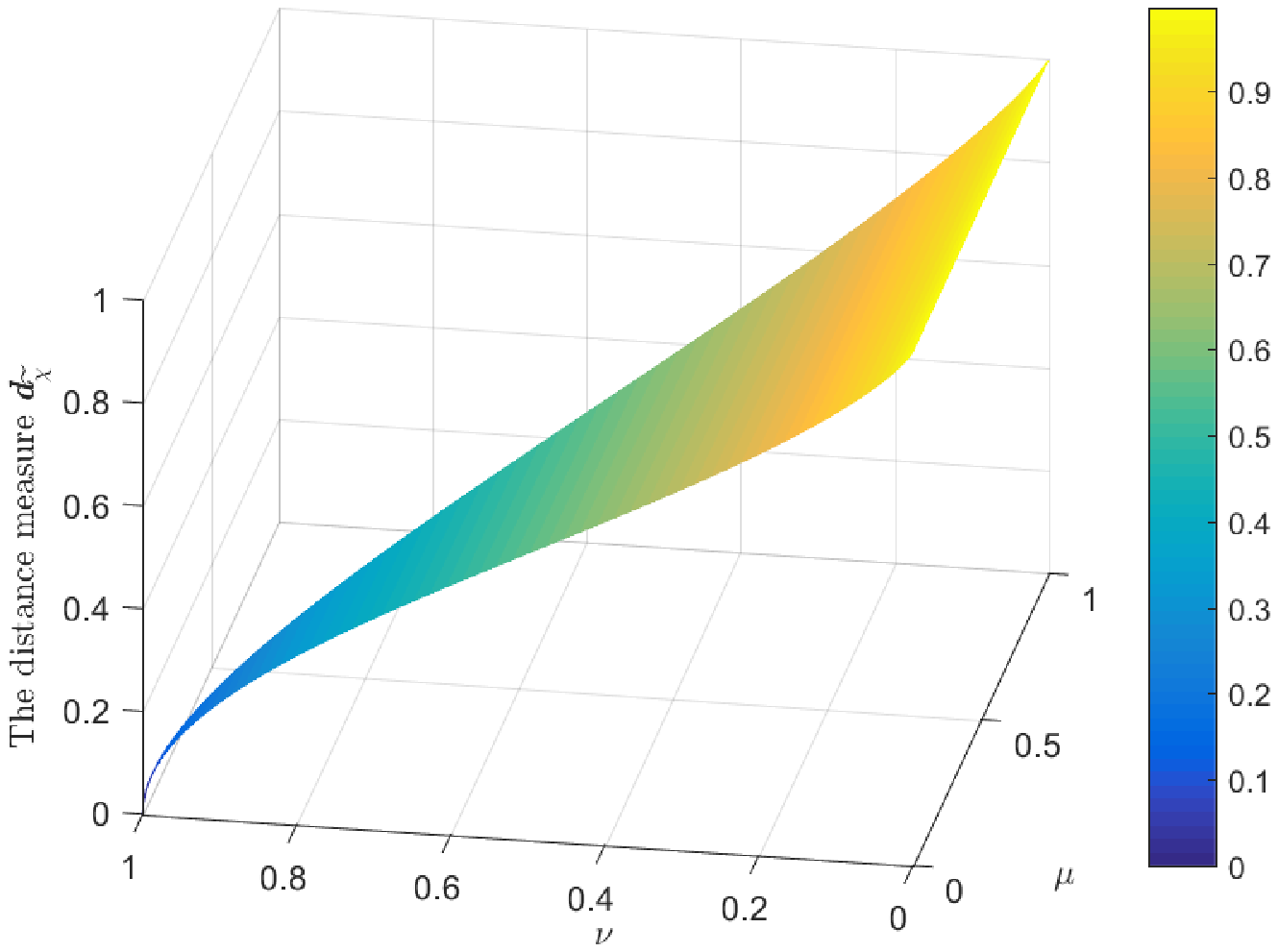}}}
\subfigure[Geometric distribution of $\bm{d}_{\widetilde{\chi}}(I_1,$ $I_2)$
in Example~\ref{Exm-Wu-3-Sec5}]{\scalebox{0.295}{\includegraphics[]{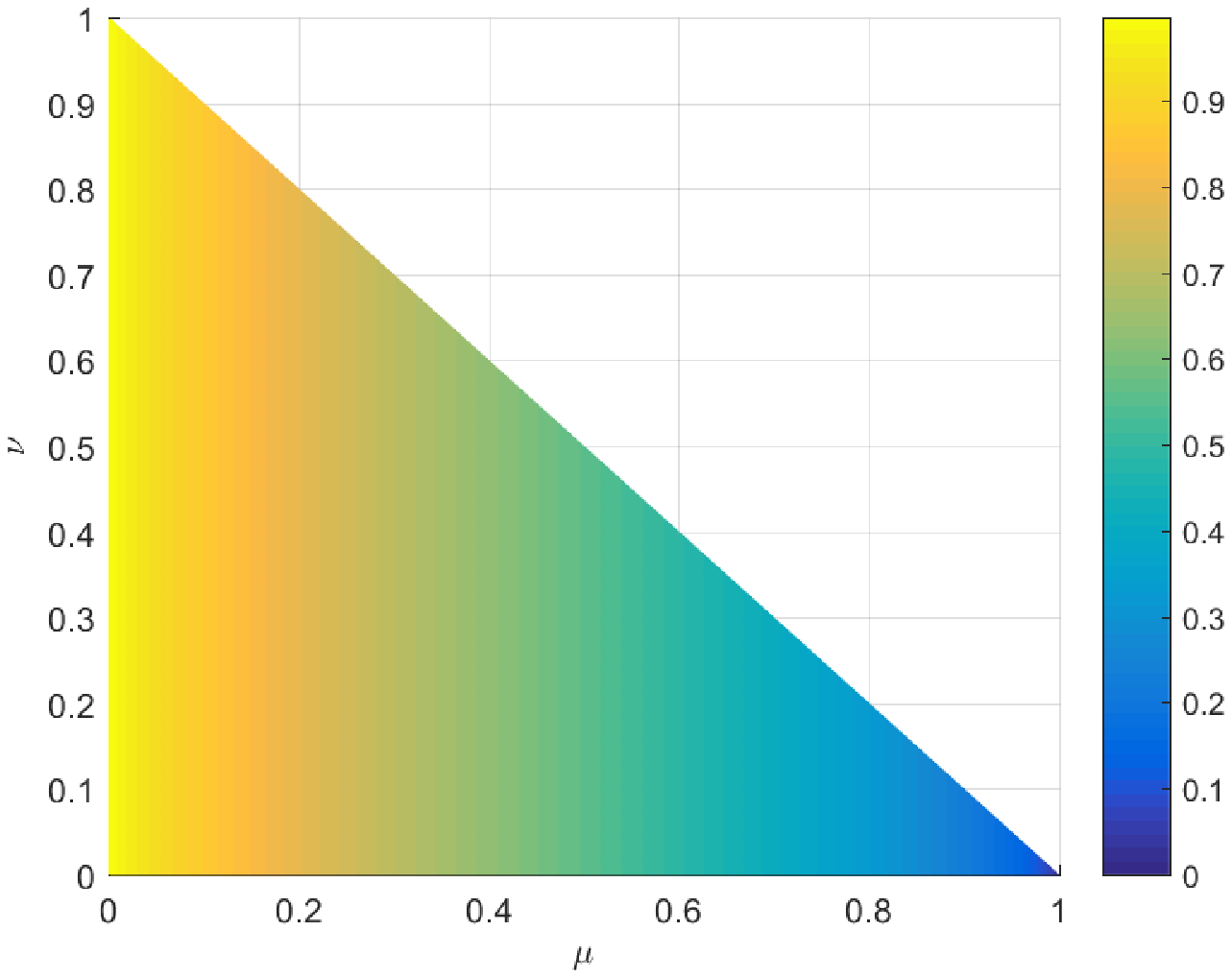}}}
\subfigure[Geometric distribution of $\bm{d}_{\widetilde{\chi}}(I_1^{\prime},$ $I_2)$
in Example~\ref{Exm-Wu-3-Sec5}]{\scalebox{0.295}{\includegraphics[]{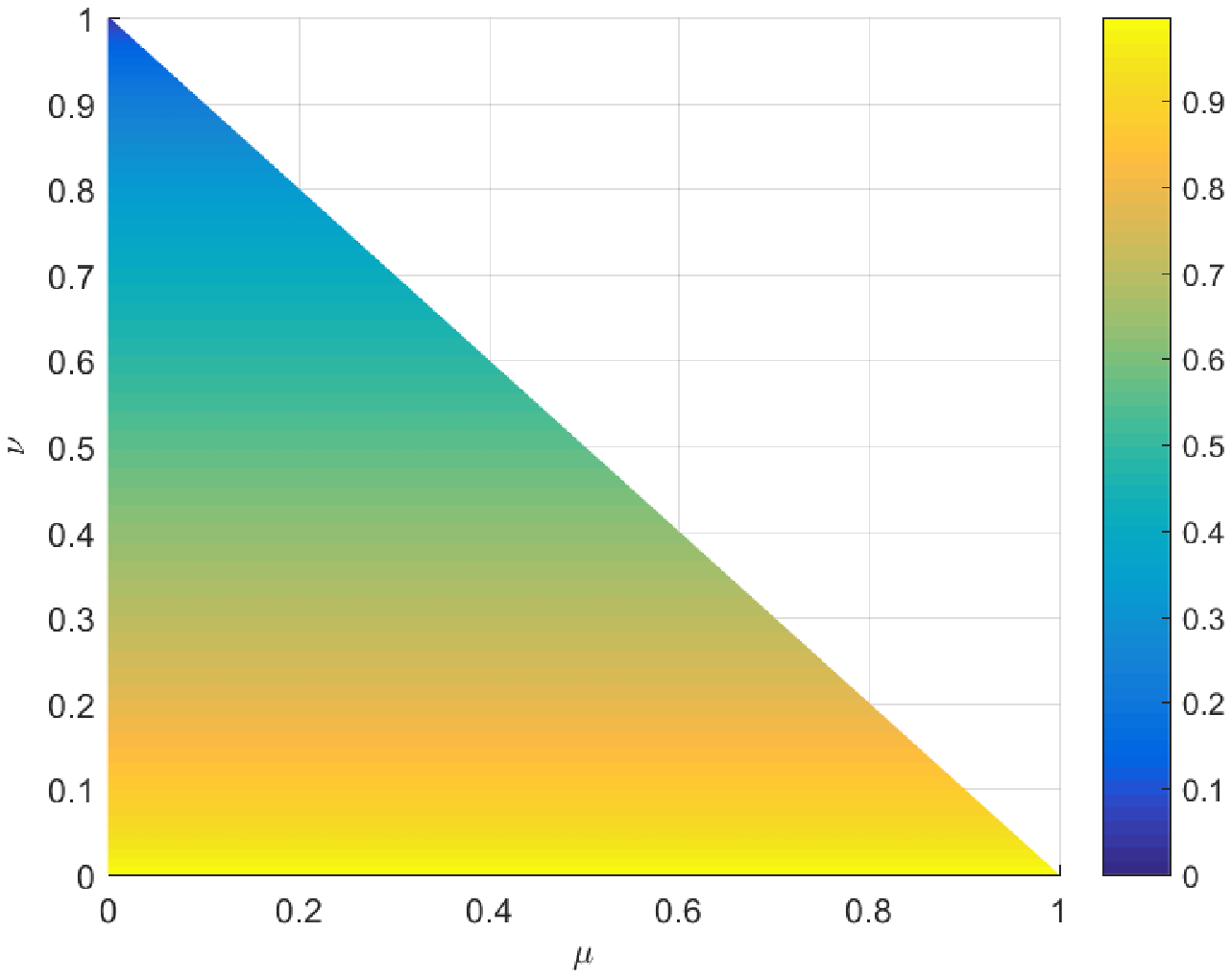}}}
\caption{The distance measures $\bm{d}_{\widetilde{\chi}}(I_1, I_2)$ and
$\bm{d}_{\widetilde{\chi}}(I_1^{\prime}, I_2)$ in Example~\ref{Exm-Wu-3-Sec5}}
\label{Fig-Xiao}
\end{figure}

\begin{figure}[H]
\centering
\subfigure[The distance $\bm{d}_{_{\mathrm{Wu}}}(I_1, I_2)$ in Example~\ref{Exm-Wu-3-Sec5}]
{\scalebox{0.295}{\includegraphics[]{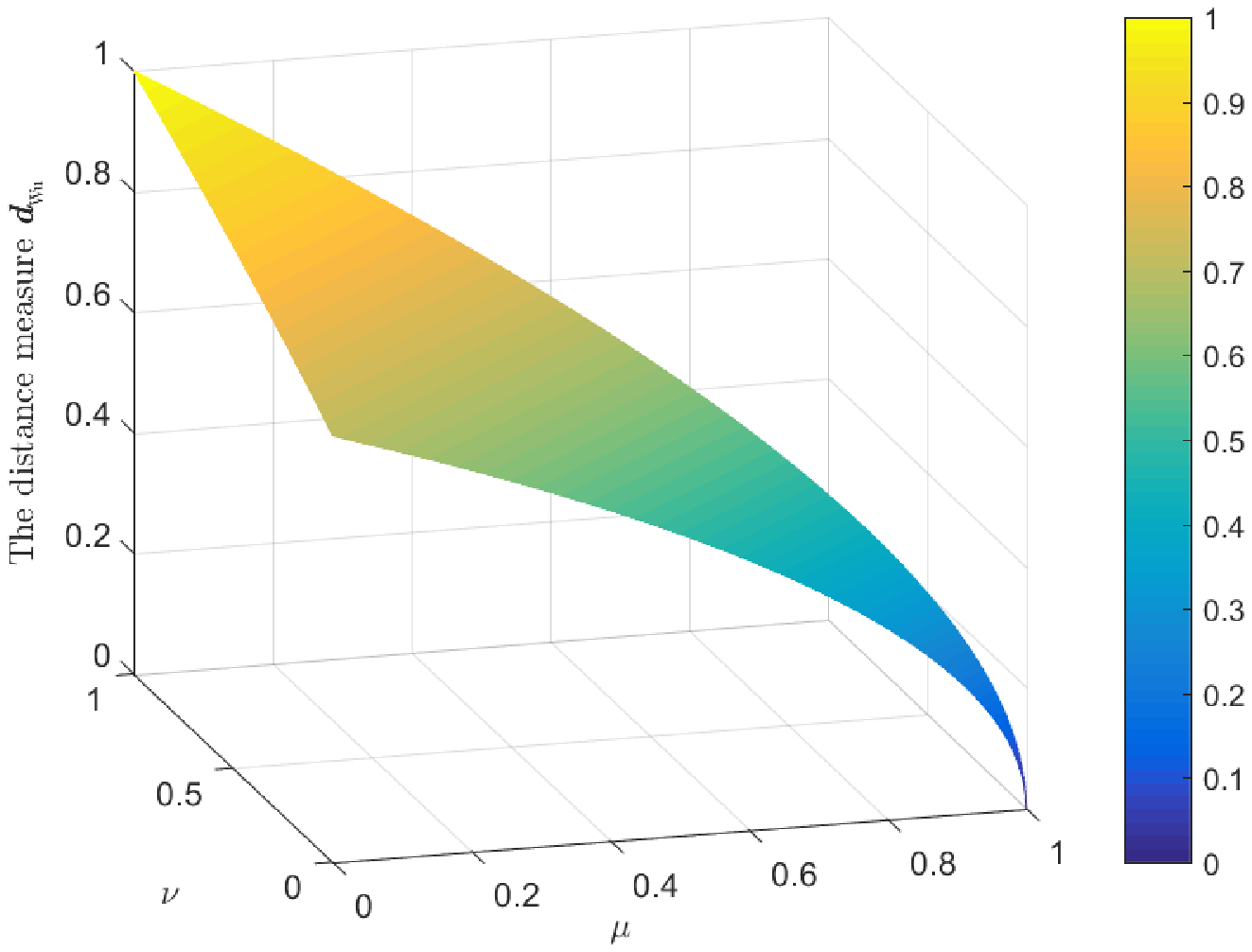}}}
\subfigure[The distance $\bm{d}_{_{\mathrm{Wu}}}(I_1^{\prime}, I_2)$ in Example~\ref{Exm-Wu-3-Sec5}]
{\scalebox{0.295}{\includegraphics[]{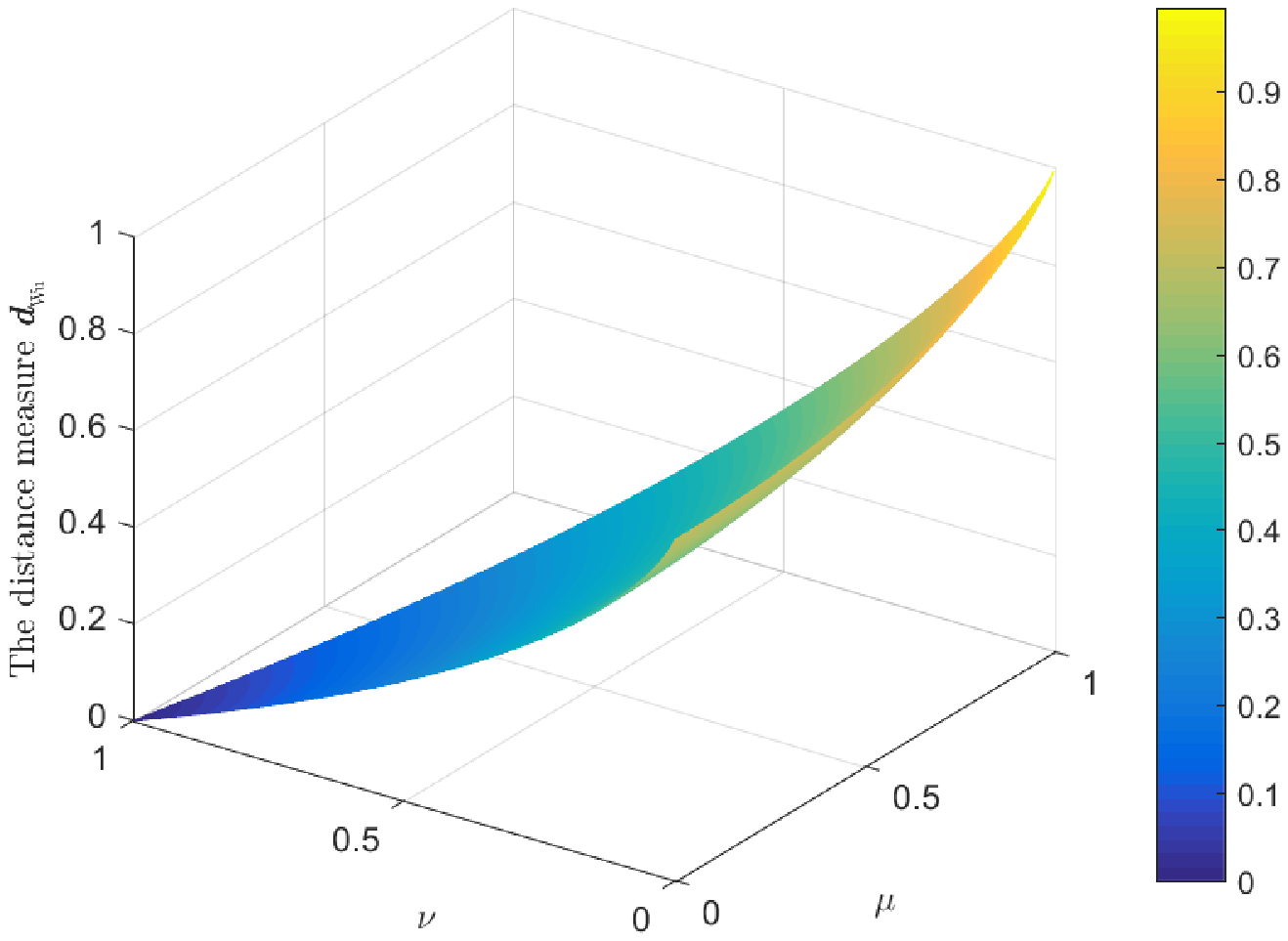}}}
\subfigure[Geometric distribution of $\bm{d}_{_{\mathrm{Wu}}}(I_1,$ $I_2)$ in Example~\ref{Exm-Wu-3-Sec5}]
{\scalebox{0.295}{\includegraphics[]{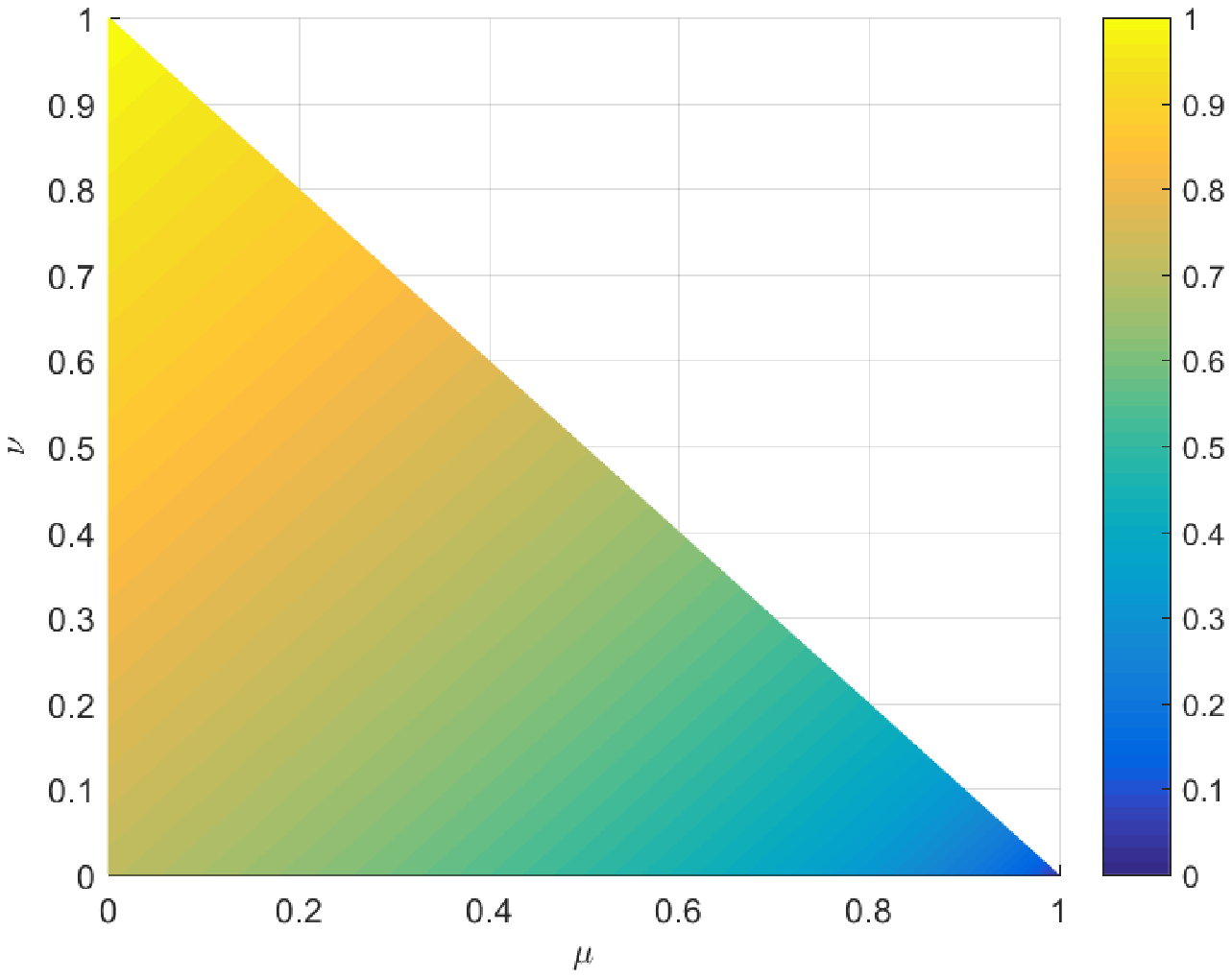}}}
\subfigure[Geometric distribution of $\bm{d}_{_{\mathrm{Wu}}}$ $(I_1^{\prime}, I_2)$ in Example~\ref{Exm-Wu-3-Sec5}]
{\scalebox{0.295}{\includegraphics[]{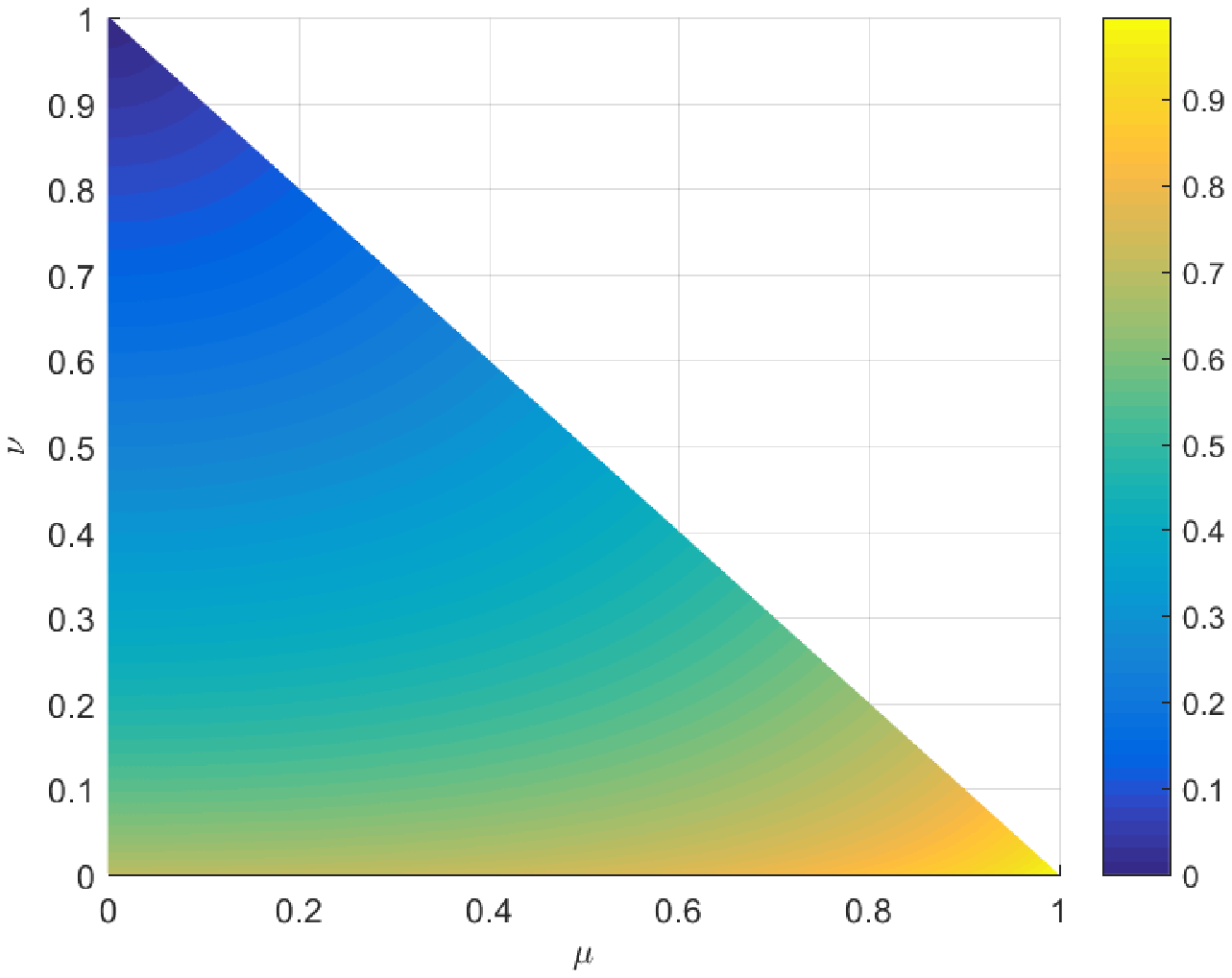}}}
\caption{The distance measures $\bm{d}_{_{\mathrm{Wu}}}(I_1, I_2)$ and
$\bm{d}_{_{\mathrm{Wu}}}(I_1^{\prime}, I_2)$ in Example~\ref{Exm-Wu-3-Sec5}}
\label{Fig-Wu}
\end{figure}
\end{example}

\subsection{Comparative analysis between Yang and Chiclana's spherical distance
 $d_{_{\mathrm{YC}}}$ and our proposed distance $d_{_{\mathrm{Wu}}}$}

\begin{example}[\textrm{Continuation of Example~\protect\ref{Exm-YC-2-Sec2}}]
\label{Exm-Wu-2-1-Sec5}
Let IFSs $I_1$, $I_2^{(\lambda)}$, and $I_3^{(\lambda)}$ on UOD
$X=\{x\}$ be given as in Example~\ref{Exm-YC-2-Sec2}. Together with Eqs.~\eqref{eq-Exm-YC-1},
\eqref{Eq-Wu-Exm8-1}, and \eqref{Eq-Wu-Exm8-3}, by varying the parameter $\lambda$ from $0$ to $1$,
Fig.~\ref{Fig-4} visually shows the changing trend of distances between
$I_{2}^{(\lambda)}$ and $I_1$, and between $I_{3}^{(\lambda)}$ and $I_1$
in Example~\ref{Exm-Wu-2-1-Sec5}, by using the distance measure $d_{_{\mathrm{YC}}}$
and our proposed distance measure $\bm{d}_{_{\mathrm{Wu}}}$. The simulation results are
consistent with our calculation. This example demonstrates that our proposed
distance is far superior to Yang and Chiclana's spherical distance in \cite{YC2009}.
\begin{figure}[H]
\centering
{\scalebox{0.33}{\includegraphics[]{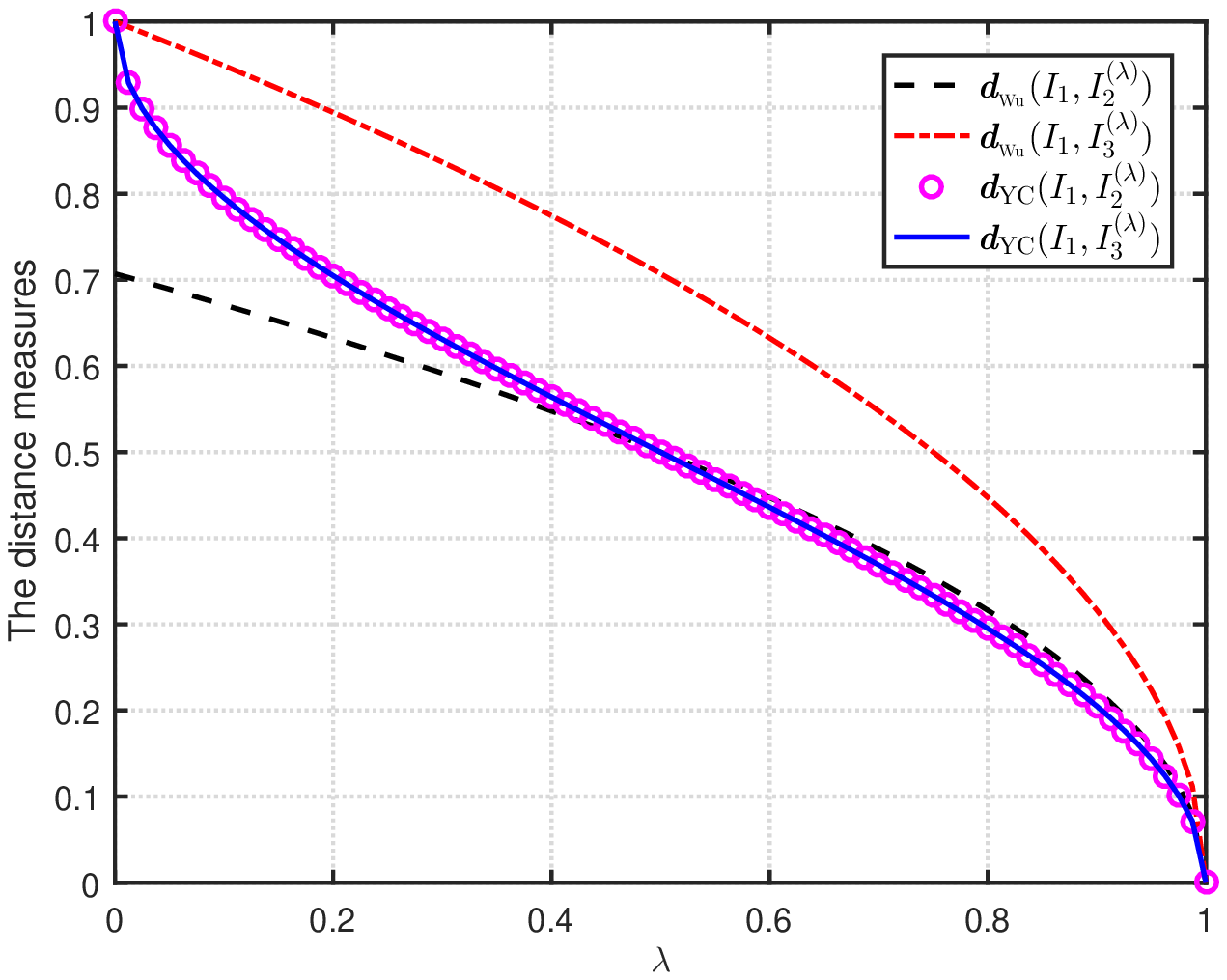}}}
\caption{The distances between $I_{2}^{(\lambda)}$
and $I_1$, and between $I_{3}^{(\lambda)}$ and $I_1$
in Example~\ref{Exm-Wu-2-1-Sec5}}
\label{Fig-4}
\end{figure}
\end{example}

\begin{example}
\label{Exm-Wu-5-Sec5}
Fig.~\ref{Fig-YC} visualizes the changing trend of
distances between $I_{2}=\left\{\frac{\langle \mu, \nu \rangle}{x}\right\}$
and $I_1=\left\{\frac{\langle 1, 0 \rangle}{x}\right\}$, and
between $I_{2}=\left\{\frac{\langle \mu, \nu \rangle}{x}\right\}$
and $I_1^{\prime}=\left\{\frac{\langle 0, 1 \rangle}{x}\right\}$
with varying the parameters $\mu$ and $\nu$ satisfying $\langle \mu,
\nu \rangle \in \Theta$, by using Yang and Chiclana's spherical distance
$\bm{d}_{_{\mathrm{YC}}}$.

From Figs.~\ref{Fig-YC} (c) and (d), we can observe that
(1) the distance $d_{_{\mathrm{YC}}}(I_1, I_2)$ between
$I_1=\left\{\frac{\langle 1, 0 \rangle}{x}\right\}$ and
$I_2=\left\{\frac{\langle \mu, \nu \rangle}{x}\right\}$
remains unchanged when the membership degree $\mu$ is fixed;
(2) the distance $d_{_{\mathrm{YC}}}(I_1^{\prime}, I_2)$ between
$I_1^{\prime}=\left\{\frac{\langle 0, 1 \rangle}{x}\right\}$ and
$I_2=\left\{\frac{\langle \mu, \nu \rangle}{x}\right\}$
remains unchanged when the non-membership degree $\nu$ is fixed.
These are consistent with the following calculation results:
\begin{align*}
d_{_{\mathrm{YC}}}(I_1, I_2)=
\frac{2}{\pi}\arccos \sqrt{\mu},
\end{align*}
and
\begin{align*}
d_{_{\mathrm{YC}}}(I_1^{\prime}, I_2)=
\frac{2}{\pi}\arccos \sqrt{\nu}.
\end{align*}

Contrarily to Fig.~\ref{Fig-Wu}, the above results also show the
unreasonableness of Yang and Chiclana's spherical distance
$d_{_{\mathrm{YC}}}$ and indicate the superiority of
our distance measure.
\begin{figure}[H]
\centering
\subfigure[The distance $d_{_{\mathrm{YC}}}(I_1, I_2)$ in Example~\ref{Exm-Wu-5-Sec5}]
{\scalebox{0.295}{\includegraphics[]{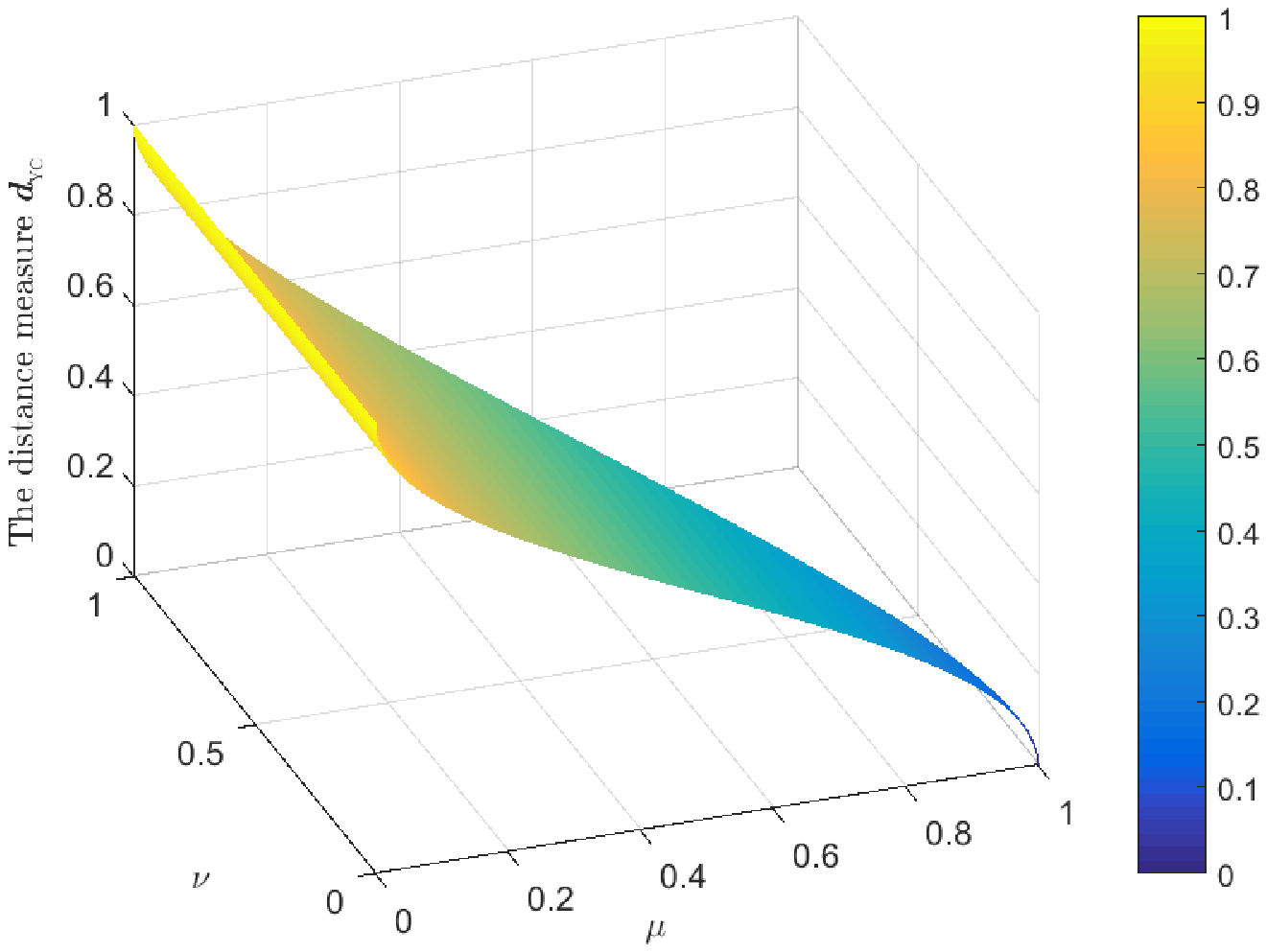}}}
\subfigure[The distance $d_{_{\mathrm{YC}}}(I_1^{\prime}, I_2)$ in Example~\ref{Exm-Wu-5-Sec5}]
{\scalebox{0.295}{\includegraphics[]{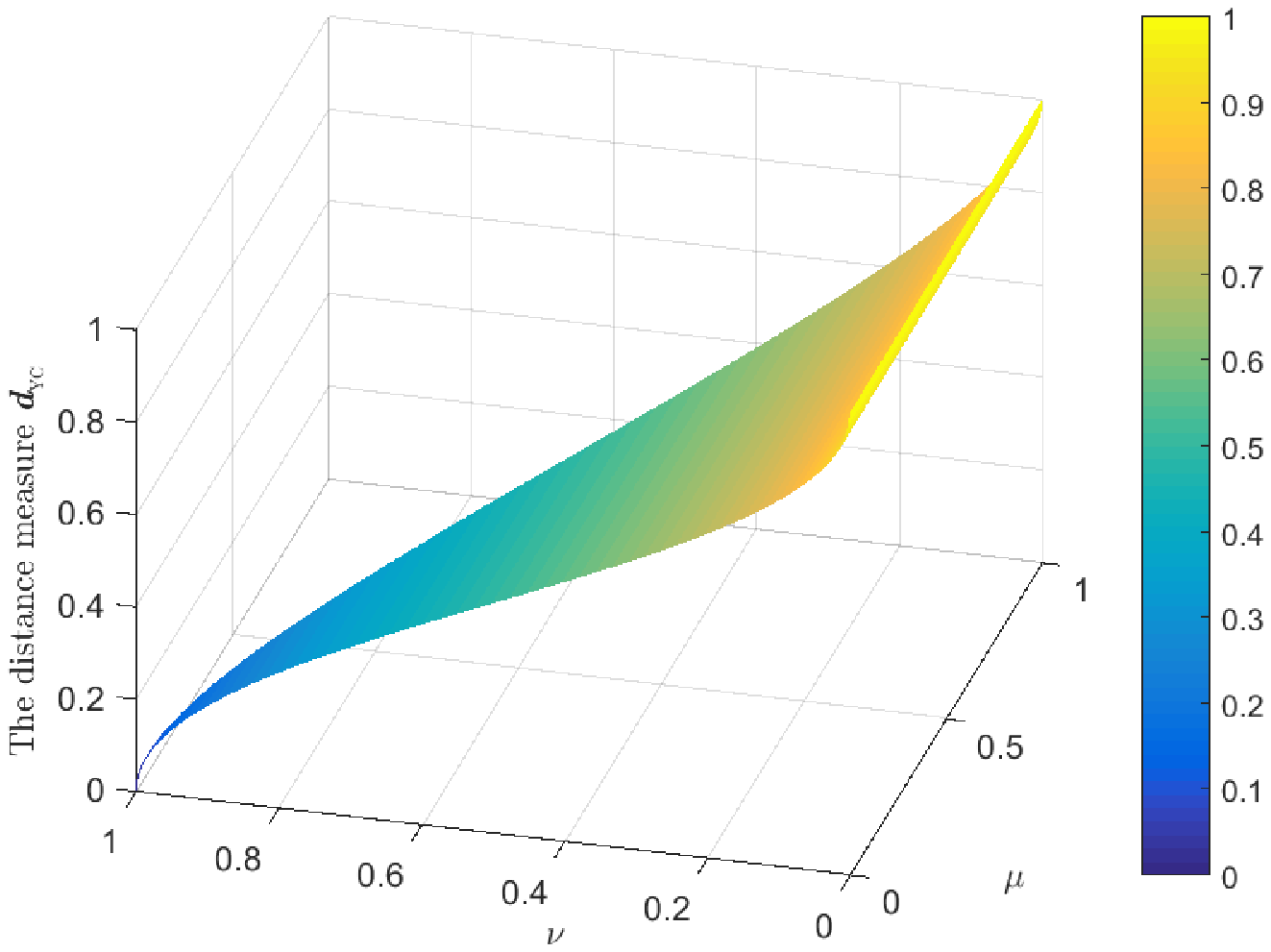}}}
\subfigure[Geometric distribution of $d_{_{\mathrm{YC}}}(I_1,$ $I_2)$
in Example~\ref{Exm-Wu-5-Sec5}]{\scalebox{0.295}{\includegraphics[]{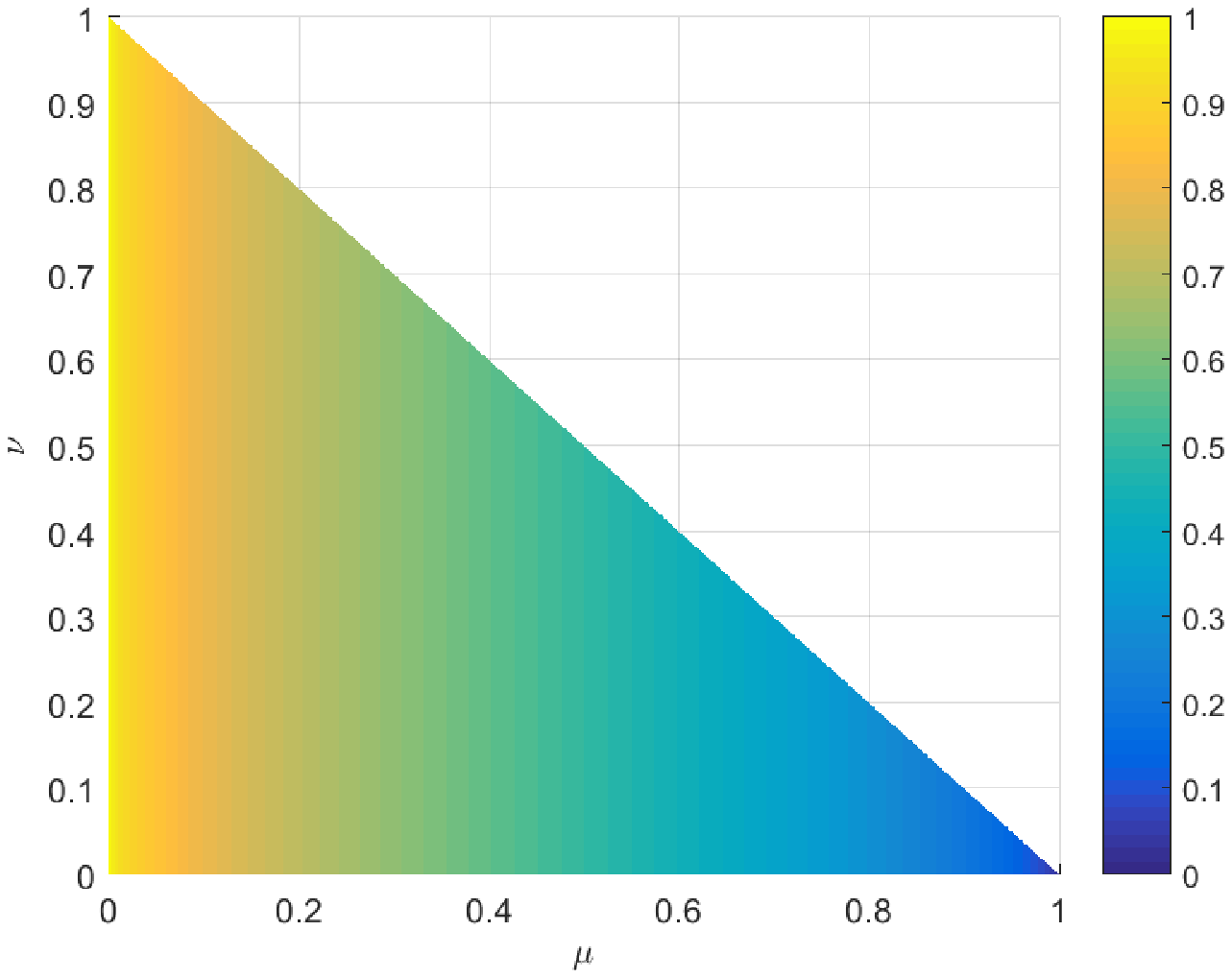}}}
\subfigure[Geometric distribution of $d_{_{\mathrm{YC}}}(I_1^{\prime},$ $I_2)$
in Example~\ref{Exm-Wu-5-Sec5}]{\scalebox{0.295}{\includegraphics[]{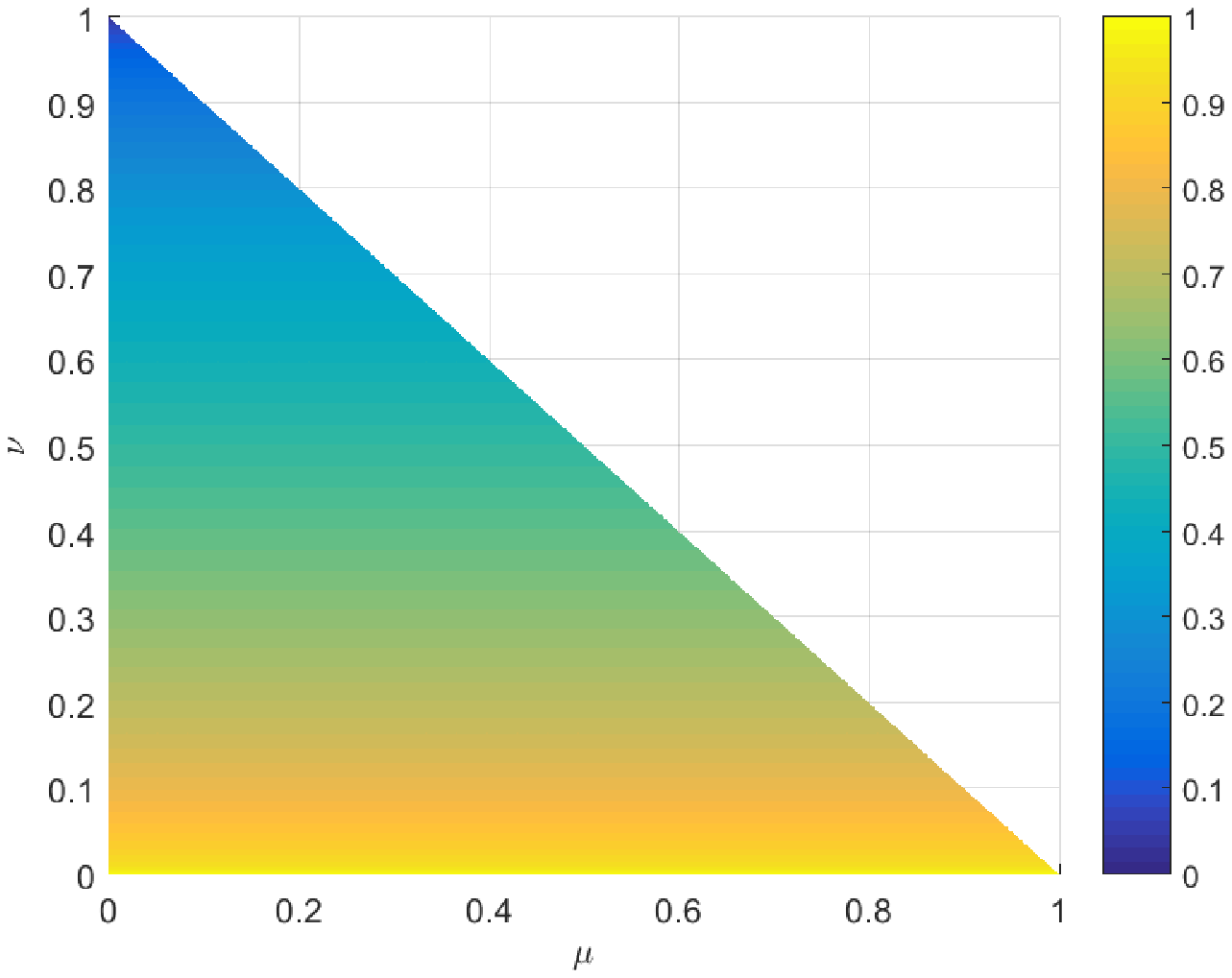}}}
\caption{The distance measures $d_{_{\mathrm{YC}}}(I_1, I_2)$ and
$d_{_{\mathrm{YC}}}(I_1^{\prime}, I_2)$ in Example~\ref{Exm-Wu-5-Sec5}}
\label{Fig-YC}
\end{figure}
\end{example}

\section{Application to pattern recognition}

In practical applications, in order to better distinguish highly similar but inconsistent IFSs,
we introduce the following parametric distance and similarity measures for IFSs.

Let $X=\{x_1, x_2, \ldots, x_n\}$ be a finite UOD and {\small $I_1=\Big\{\frac{\alpha_j^{(1)}}{x_j}
\mid 1\leq j\leq n,  \alpha_j^{(1)}\in \Theta\Big\}$} and {\small $I_2=\Big\{\frac{\alpha_j^{(2)}}{x_j}
\mid 1\leq j\leq n, \alpha_j^{(2)}\in \Theta\Big\}$} be two IFSs on $X$. For $\lambda>0$,
define
\begin{small}
\begin{equation}
\label{Eq-Wu-Para-1}
\begin{split}
& \bm{d}_{_{\mathrm{Wu}}}^{(\lambda)}(I_1, I_2)= \sum_{j=1}^{n}
\omega_j \left\{ \frac{1}{2}\right. \\
& \quad \left[ (1-(\mu_{I_1}(x_j))^{\lambda})  \cdot  \log_{2}
\frac{2 (1-(\mu_{I_1}(x_j))^{\lambda})}{(1-(\mu_{I_1}(x_j))^{\lambda})
+(1-(\mu_{I_2}(x_j))^{\lambda})}\right.\\
& \quad +(1-(\mu_{I_2}(x_j))^{\lambda}) \cdot \log_{2}\frac{2 (1-(\mu_{I_2}(x_j))^{\lambda})}
{(1-(\mu_{I_1}(x_j))^{\lambda})+(1-(\mu_{I_2}(x_j))^{\lambda})}\\
& \quad +(\nu_{I_1}(x_j))^{\lambda} \cdot \log_{2}\frac{2 (\nu_{I_1}(x_j))^{\lambda}}
{(\nu_{I_1}(x_j))^{\lambda}+(\nu_{I_2}(x_j))^{\lambda}}\\
& \left.\left. \quad +(\nu_{I_2}(x_j))^{\lambda} \cdot \log_{2}\frac{2 (\nu_{I_2}(x_j))^{\lambda}}{(\nu_{I_1}(x_j))^{\lambda}+(\nu_{I_2}(x_j))^{\lambda}
}\right]\right\}^{0.5},
\end{split}
\end{equation}
\end{small}
and
\begin{equation}
\label{Eq-Wu-Para-2}
\mathbf{S}_{_{\mathrm{Wu}}}^{(\lambda)}(I_1, I_2)=
1-\bm{d}_{_{\mathrm{Wu}}}^{(\lambda)}(I_1, I_2),
\end{equation}
where $\alpha_j^{(1)}=\langle \mu_{I_1}(x_j),
\nu_{I_1}(x_j)\rangle$, $\alpha_j^{(2)}=\langle \mu_{I_2}(x_j),
\nu_{I_2}(x_j) \rangle$, and $\omega=(\omega_1, \omega_2, \ldots,
\omega_n)^{\top}$ is the weight vector of $x_{j}$ ($j=1, 2, \ldots, n$)
with $\omega_j\in (0, 1]$ and $\sum_{j=1}^{n}\omega_j=1$.

Similarly to the discussions in Section~\ref{Sec-V}, we have the following result.

\begin{theorem} For $\lambda>0$,
\begin{enumerate}[(1)]
  \item the function $\bm{d}_{_{\mathrm{Wu}}}^{(\lambda)}$ defined by
Eq.~\eqref{Eq-Wu-Para-1} is a SIFDisM on $\mathrm{IFS}(X)$;
  \item the function $\mathbf{S}_{_{\mathrm{Wu}}}^{(\lambda)}$
defined by Eq.~\eqref{Eq-Wu-Para-2} is a SIFSimM on $\mathrm{IFS}(X)$.
\end{enumerate}
\end{theorem}

Now, we utilize a practical example to illustrate the effectiveness
of our proposed distance measure.

\begin{example}
[{\textrm{\protect\cite[Application 2]{Xiao2021},\cite[Example 4.3]{LZ2018}}}]
\label{Exm-Appl-1}
Consider a pattern classification problem with three classes and three attributes
$A=\{x_1, x_2, x_3\}$, described by three patterns $P=\{P_1, P_2, P_3\}$
and a test sample $S_1$ expressed by the IFSs listed in Table~\ref{Tab-PRP}.
\begin{table}[H]	
	\centering
	\caption{Pattern classification problem with three-classes and three-attributes in Example~\ref{Exm-Appl-1}}
	\label{Tab-PRP}
\scalebox{0.8}{
	\begin{tabular}{cccccccc}
		\toprule
		 \multirow{3}*{} & \multirow{3}*{} & \multicolumn{6}{c}{Attribute}\\
\cline{3-8}
         ~ & ~ & \multicolumn{2}{c}{$x_{1}$} & \multicolumn{2}{c}{$x_{2}$} & \multicolumn{2}{c}{$x_{3}$}\\
\cline{3-8}
         ~ & ~ & $\mu_{P}(x_1)$ & $\nu_{P}(x_1)$ & $\mu_{P}(x_2)$ & $\nu_{P}(x_2)$ & $\mu_{P}(x_3)$ & $\nu_{P}(x_3)$\\
		\midrule
		\multirow{3}*{$\textrm{Pattern}$} & $P_1$ & 0.15 & 0.25 & 0.25 & 0.35 & 0.35 & 0.45\\
		~ & $P_2$ & 0.05 & 0.15 & 0.15 & 0.25 & 0.25 & 0.35\\
        ~ & $P_3$ & 0.16 & 0.26 & 0.26 & 0.36 & 0.36 & 0.46\\
        \midrule
		$\textrm{Sample}$ & $S_1$ & 0.30 & 0.20 & 0.40 & 0.30 & 0.50 & 0.40\\
		\bottomrule
	\end{tabular}
}
\end{table}

\begin{table}[H]	
	\centering
	\caption{Pattern recognition results by different similarity measures
in Example~\ref{Exm-Appl-1}}
	\label{Tab-PR-Diff-Sim-Mesaure}
\scalebox{0.8}{
\begin{threeparttable}
	\begin{tabular}{ccccc}
		\toprule
		 \multirow{2}*{Method} & \multicolumn{3}{c}{Similarity measure} & \multirow{2}*{Classification}\\
\cline{2-4}
         ~ & $1-\textrm{dis}(P_1, S_1)$ & $1-\textrm{dis}(P_2, S_1)$ & $1-\textrm{dis}(P_3, S_1)$ & ~\\
		\midrule
		$d_{_{\mathrm{SK}}}^{\mathrm{H}}$ & 0.85 & 0.70 & 0.86 & $P_{3}$\\
		$d_{_{\mathrm{SK}}}^{\mathrm{E}}$ & 0.87 & 0.72 & 0.88 & $P_{3}$\\
        $d_{_{\mathrm{G}}}$ & 0.85 & 0.75 & 0.86 & $P_{3}$\\
		$d_{_{\mathrm{W1}}}$ & 0.90 & 0.80 & 0.91 & $P_{3}$\\
        $d_{_{\mathrm{W2}}}$ & 0.90 & 0.85 & 0.90 & \XSolid\\
        $d_{_{\mathrm{P}}}$ & 0.85 & 0.70 & 0.86 & $P_{3}$\\
        $d_{_{\mathrm{Y}}}$ & 0.85 & 0.70 & 0.86 & $P_{3}$\\
        $d_{_{\mathrm{H}}}^{\mathrm{T}}$ & 0.95 & 0.88 & 0.95 & \XSolid\\
        $d_{_{\mathrm{H}}}^{\mathrm{R}}$ & 0.96 & 0.93 & 0.96 & \XSolid\\
        $d_{_{\mathrm{H}}}^{\mathrm{L}}$ & $1-3.70\times10^{-17}$ & $1-3.70\times10^{-17}$ & $1-3.70\times10^{-17}$ & \XSolid\\
        $d_{_{\mathrm{H}}}^{\mathrm{KD}}$ & 0.90 & 0.85 & 0.90 & \XSolid\\
        $d_{_{\mathrm{H}}}^{\mathrm{M}}$ & 0.90 & 0.85 & 0.90 & \XSolid\\
        $d_{_{\mathrm{H}}}^{\mathrm{LA}}$ & 0.93 & 0.92 & 0.93 & \XSolid\\
        $d_{_{\mathrm{H}}}^{\mathrm{G}}$ & 0.95 & 0.92 & 0.95 & \XSolid\\
        $d_{_{\mathrm{SW}}}$ & 0.99 & 0.95 & 0.99 & \XSolid\\
        $d_{_{\mathrm{SM}}}$ & 0.86 & 0.81 & 0.90 & $P_{3}$\\
        $d_{_{\mathrm{L}}}$ & 0.80 & 0.60 & 0.81 & $P_{3}$\\
        $d_{_{\mathrm{YC}}}$ & 0.89 & 0.77 & 0.90 & $P_{3}$\\
        $d_{\widetilde{\chi}}$ & 0.85 & 0.69 & 0.86 & $P_{3}$\\
        $\bm{d}_{_{\mathrm{Wu}}}^{(1/3)}$ & 0.91 & 0.84 & 0.92 & $P_{3}$\\
		\bottomrule
	\end{tabular}
\begin{tablenotes}
\item[] \XSolid ~denotes that it cannot be determined;
\item[] The details for distance measures in
Table~\ref{Tab-PR-Diff-Sim-Mesaure} can be found in \cite[Section~III]{Xiao2021}.
\end{tablenotes}
\end{threeparttable}
}
\end{table}

\begin{figure*}
\centering
{\scalebox{0.55}{\includegraphics[]{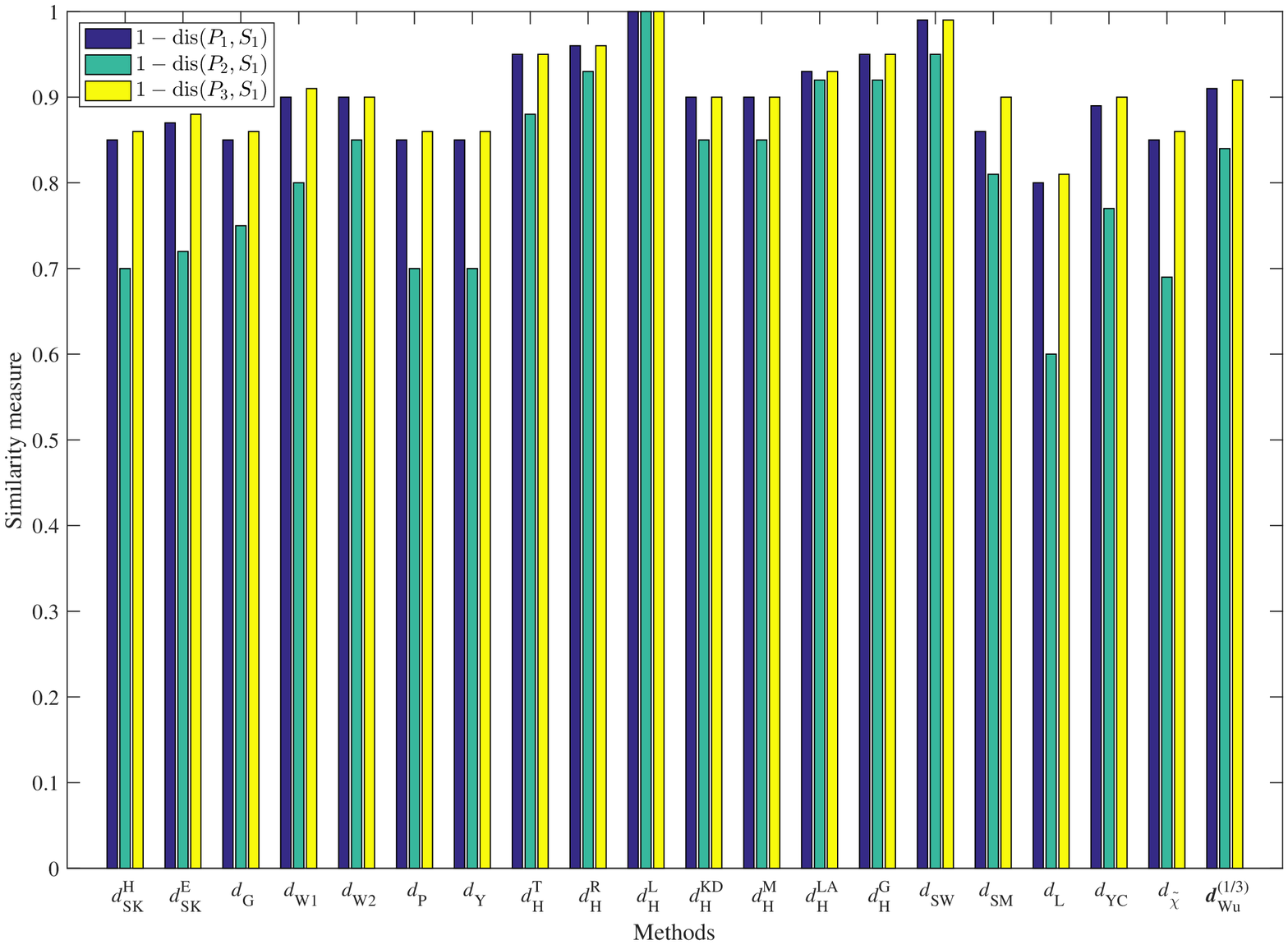}}}
\caption{Comparison results of different similarity measures in Example~\ref{Exm-Appl-1}}
\label{Fig-Application}
\end{figure*}

If we take the weight vector $\omega$ of three attributes as $\omega=(\frac{1}{3},
\frac{1}{3}, \frac{1}{3})^{\top}$, then by the principle of maximum degree of
similarity measures, the pattern classification results obtained by using
different distance measures are listed in Table~\ref{Tab-PR-Diff-Sim-Mesaure}
and Fig.~\ref{Fig-Application}. Observing from Table~\ref{Tab-PR-Diff-Sim-Mesaure}
and shown in Fig.~\ref{Fig-Application}, one can see that the test sample $S_{1}$ is classified to
the pattern $P_3$ by our proposed distance measure, which is consistent with the
results obtained by the distance measures $d_{_{\mathrm{SK}}}^{\mathrm{E}}$,
$d_{_{\mathrm{G}}}$, $d_{_{\mathrm{W1}}}$, $d_{_{\mathrm{W2}}}$, $d_{_{\mathrm{P}}}$,
$d_{_{\mathrm{Y}}}$, $d_{_{\mathrm{SW}}}$, $d_{_{\mathrm{SM}}}$, $d_{_{\mathrm{L}}}$,
$d_{_{\mathrm{YC}}}$, and $d_{\widetilde{\chi}}$; However, the methods by using
the distance measures $d_{_{\mathrm{W2}}}$, $d_{_{\mathrm{H}}}^{\mathrm{T}}$,
$d_{_{\mathrm{H}}}^{\mathrm{R}}$, $d_{_{\mathrm{H}}}^{\mathrm{L}}$,
$d_{_{\mathrm{H}}}^{\mathrm{KD}}$, $d_{_{\mathrm{H}}}^{\mathrm{M}}$,
$d_{_{\mathrm{H}}}^{\mathrm{LA}}$, $d_{_{\mathrm{H}}}^{\mathrm{G}}$,
and $d_{_{\mathrm{SW}}}$, cannot determine to which pattern the test
sample $S_1$ belongs.
\end{example}


\section{Conclusion}
This paper is devoted to the construction of SIFDisM and SIFSimM, which can effectively
measure the differences between IFSs. First, we show some examples to demonstrate that Xiao's
distance measure in \cite{Xiao2021} and Yang and Chiclana's spherical distance in
\cite{YC2009} have some shortcomings, which may cause counter-intuitive
results. To overcome these shortcomings, we present the concepts of strict
intuitionistic fuzzy distance measure (SIFDisM) and strict intuitionistic fuzzy similarity
measure (SIFSimM), and propose a novel IFDisM based on Jensen-Shannon divergence. Moreover, we
prove that the dual similarity measure of our proposed distance measure is an SIFSimM and
its induced entropy measure is an IF entropy measure. Meanwhile, we perform some comparative
analysis to illustrate that our proposed distance measure is completely superior to the
existing IFDisMs; in particular, it is much better than Xiao's distance measure in
\cite{Xiao2021}, Hung and Yang's $J_{\alpha}$-divergence in~\cite{HY2008},
Joshi and Kumar's dissimilarity divergence in~\cite{JK2018}, and Yang and
Chiclana's spherical distance in \cite{YC2009}; consequently, it is better than distance
measures in \cite{SK2000,Gr2004,WX2005,HPK2012,YC2012,SMLZC2018,SWQH2019}. Finally, to
illustrate the availability of our proposed IFSM, we apply it to a practical pattern
recognition problem. In the future, we will apply our methods to
establish new distance/similarity measures for Pythagorean fuzzy sets, q-rung
orthopair fuzzy sets, spherical fuzzy sets, picture fuzzy sets, and
some other interval-valued fuzzy sets.


\begin{thebibliography}{10}
\providecommand{\url}[1]{#1}
\csname url@samestyle\endcsname
\providecommand{\newblock}{\relax}
\providecommand{\bibinfo}[2]{#2}
\providecommand{\BIBentrySTDinterwordspacing}{\spaceskip=0pt\relax}
\providecommand{\BIBentryALTinterwordstretchfactor}{4}
\providecommand{\BIBentryALTinterwordspacing}{\spaceskip=\fontdimen2\font plus
\BIBentryALTinterwordstretchfactor\fontdimen3\font minus
  \fontdimen4\font\relax}
\providecommand{\BIBforeignlanguage}[2]{{%
\expandafter\ifx\csname l@#1\endcsname\relax
\typeout{** WARNING: IEEEtran.bst: No hyphenation pattern has been}%
\typeout{** loaded for the language `#1'. Using the pattern for}%
\typeout{** the default language instead.}%
\else
\language=\csname l@#1\endcsname
\fi
#2}}
\providecommand{\BIBdecl}{\relax}
\BIBdecl

\bibitem{Sz2014}
E.~Szmidt, \emph{Distances and {S}imilarities in {I}ntuitionistic {F}uzzy
  {S}ets}, ser. Studies in Fuzziness and Soft Computing.\hskip 1em plus 0.5em
  minus 0.4em\relax Springer, Berlin, Heidelberg, 2014, vol. 307.

\bibitem{Xu2007a}
Z.~Xu, ``Some similarity measures of intuitionistic fuzzy sets and their
  applications to multiple attribute decision making,'' \emph{Fuzzy Optim.
  Decis. Making}, vol.~6, no.~2, pp. 109--121, 2007.

\bibitem{Li2014}
D.-F. Li, \emph{Decision and {G}ame {T}heory in {M}anagement with
  {I}ntuitionistic {F}uzzy {S}ets}, ser. Studies in Fuzziness and Soft
  Computing.\hskip 1em plus 0.5em minus 0.4em\relax {S}pringer-{V}erlag
  {B}erlin {H}eidelberg, 2014, vol. 308.

\bibitem{WCZ2022}
\BIBentryALTinterwordspacing
X.~Wu, Z.~Zhu, G.~Chen, and P.~Liu, ``A monotonous intuitionistic fuzzy
  {TOPSIS} method with linear orders under two novel admissible distance
  measures,'' \emph{submitted to IEEE Trans. Fuzzy Syst.}, 2022. [Online].
  Available: \url{https://arxiv.org/abs/2206.02567}
\BIBentrySTDinterwordspacing

\bibitem{SK2000}
E.~Szmidt and J.~Kacprzyk, ``Distances between intuitionistic fuzzy sets,''
  \emph{Fuzzy Sets Syst.}, vol. 114, no.~3, pp. 505--518, 2000.

\bibitem{Gr2004}
P.~Grzegorzewski, ``Distances between intuitionistic fuzzy sets and/or
  interval-valued fuzzy sets based on the {H}ausdorff metric,'' \emph{Fuzzy
  Sets Syst.}, vol. 148, no.~2, pp. 319--328, 2004.

\bibitem{HY2004}
W.-L. Hung and M.-S. Yang, ``Similarity measures of intuitionistic fuzzy sets
  based on {H}ausdorff distance,'' \emph{Pattern Recognit. Lett.}, vol.~25,
  no.~14, pp. 1603--1611, 2004.

\bibitem{WX2005}
W.~Wang and X.~Xin, ``Distance measure between intuitionistic fuzzy sets,''
  \emph{Pattern Recognit. Lett.}, vol.~26, no.~13, pp. 2063--2069, 2005.

\bibitem{XC2008}
Z.~Xu and J.~Chen, ``An overview of distance and similarity measures of
  intuitionistic fuzzy sets,'' \emph{Int. J. Uncertain. Fuzziness. Knowl.-Based
  Syst.}, vol.~16, no.~04, pp. 529--555, 2008.

\bibitem{HY1981}
C.-L. Hwang and K.~Yoon, \emph{Multiple Attribute Decision Making: Methods and
  Applications A State-of-the-Art Survey}, ser. Lecture Notes in Economics and
  Mathematical Systems.\hskip 1em plus 0.5em minus 0.4em\relax
  {S}pringer-{V}erlag {B}erlin {H}eidelberg, 1981, vol. 186.

\bibitem{CCL2016}
S.-M. Chen, S.-H. Cheng, and T.-C. Lan, ``A novel similarity measure between
  intuitionistic fuzzy sets based on the centroid points of transformed fuzzy
  numbers with applications to pattern recognition,'' \emph{Inf. Sci.}, vol.
  343, pp. 15--40, 2016.

\bibitem{YC2009}
Y.~Yang and F.~Chiclana, ``Intuitionistic fuzzy sets: spherical representation
  and distances,'' \emph{Int. J. Intell. Syst.}, vol.~24, no.~4, pp. 399--420,
  2009.

\bibitem{HY2008}
W.-L. Hung and M.-S. Yang, ``On the {$J$}-divergence of intuitionistic fuzzy
  sets with its application to pattern recognition,'' \emph{Inf. Sci.}, vol.
  178, no.~6, pp. 1641--1650, 2008.

\bibitem{JK2018}
R.~Joshi and S.~Kumar, ``A dissimilarity {J}ensen--{S}hannon divergence measure
  for intuitionistic fuzzy sets,'' \emph{Int. J. Intell. Syst.}, vol.~33,
  no.~11, pp. 2216--2235, 2018.

\bibitem{Xiao2021}
F.~Xiao, ``A distance measure for intuitionistic fuzzy sets and its application
  to pattern classification problems,'' \emph{IEEE Trans. Syst., Man, Cybern.,
  Syst.}, vol.~51, no.~6, pp. 3980--3992, 2021.

\bibitem{HPK2012}
A.~G. Hatzimichailidis, G.~A. Papakostas, and V.~G. Kaburlasos, ``A novel
  distance measure of intuitionistic fuzzy sets and its application to pattern
  recognition problems,'' \emph{Int. J. Intell. Syst.}, vol.~27, no.~4, pp.
  396--409, 2012.

\bibitem{YC2012}
Y.~Yang and F.~Chiclana, ``Consistency of {2D} and {3D} distances of
  intuitionistic fuzzy sets,'' \emph{Expert Syst. Appl.}, vol.~39, no.~10, pp.
  8665--8670, 2012.

\bibitem{SMLZC2018}
F.~Shen, X.~Ma, Z.~Li, Z.~Xu, and D.~Cai, ``An extended intuitionistic fuzzy
  {TOPSIS} method based on a new distance measure with an application to credit
  risk evaluation,'' \emph{Inf. Sci.}, vol. 428, pp. 105--119, 2018.

\bibitem{SWQH2019}
Y.~Song, X.~Wang, W.~Quan, and W.~Huang, ``A new approach to construct
  similarity measure for intuitionistic fuzzy sets,'' \emph{Soft Comput.},
  vol.~23, no.~6, pp. 1985--1998, 2019.

\bibitem{XC2012}
Z.~Xu and X.~Cai, \emph{{I}ntuitionistic {F}uzzy {I}nformation {A}ggregation:
  {T}heory and {A}pplications}, ser. Mathematics Monograph Series.\hskip 1em
  plus 0.5em minus 0.4em\relax {S}cience {P}ress, 2012, vol.~20.

\bibitem{Ata1999}
K.~T. Atanassov, \emph{{I}ntuitionistic {F}uzzy {S}ets: {T}heory and
  {A}pplications}, ser. Studies in Fuzziness and Soft Computing.\hskip 1em plus
  0.5em minus 0.4em\relax {S}pringer-{V}erlag {B}erlin {H}eidelberg, 1999,
  vol.~35.

\bibitem{Xu2007}
Z.~Xu, ``Intuitionistic fuzzy aggregation operators,'' \emph{IEEE Trans. Fuzzy
  Syst.}, vol.~15, no.~6, pp. 1179--1187, 2007.

\bibitem{SK2001}
E.~Szmidt and J.~Kacprzyk, ``Entropy for intuitionistic fuzzy sets,''
  \emph{Fuzzy Sets Syst.}, vol. 118, no.~3, pp. 467--477, 2001.

\bibitem{ES2003}
D.~M. Endres and J.~E. Schindelin, ``A new metric for probability
  distributions,'' \emph{IEEE Trans. Inf. Theory}, vol.~49, no.~7, pp.
  1858--1860, 2003.

\bibitem{LZ2018}
M.~Luo and R.~Zhao, ``A distance measure between intuitionistic fuzzy sets and
  its application in medical diagnosis,'' \emph{Artif. Intell. Med.}, vol.~89,
  pp. 34--39, 2018.

\end{thebibliography}


\begin{thebibliography}{10}
\providecommand{\url}[1]{#1}
\csname url@samestyle\endcsname
\providecommand{\newblock}{\relax}
\providecommand{\bibinfo}[2]{#2}
\providecommand{\BIBentrySTDinterwordspacing}{\spaceskip=0pt\relax}
\providecommand{\BIBentryALTinterwordstretchfactor}{4}
\providecommand{\BIBentryALTinterwordspacing}{\spaceskip=\fontdimen2\font plus
\BIBentryALTinterwordstretchfactor\fontdimen3\font minus
  \fontdimen4\font\relax}
\providecommand{\BIBforeignlanguage}[2]{{%
\expandafter\ifx\csname l@#1\endcsname\relax
\typeout{** WARNING: IEEEtran.bst: No hyphenation pattern has been}%
\typeout{** loaded for the language `#1'. Using the pattern for}%
\typeout{** the default language instead.}%
\else
\language=\csname l@#1\endcsname
\fi
#2}}
\providecommand{\BIBdecl}{\relax}
\BIBdecl

\bibitem{Za1965}
L.~A. Zadeh, ``Fuzzy sets,'' \emph{Inf. Control}, vol.~8, no.~3, pp. 338--353,
  1965.

\bibitem{Ata1986}
K.~T. Atanassov, ``Intuitionistic fuzzy set,'' \emph{Fuzzy Sets Syst.},
  vol.~20, pp. 87--96, 1986.

\bibitem{Ata1999}
------, \emph{{I}ntuitionistic {F}uzzy {S}ets: {T}heory and {A}pplications},
  ser. Studies in Fuzziness and Soft Computing.\hskip 1em plus 0.5em minus
  0.4em\relax {S}pringer-{V}erlag {B}erlin {H}eidelberg, 1999, vol.~35.

\bibitem{AG1989}
K.~T. Atanassov and G.~Gargov, ``Interval valued intuitionistic fuzzy sets,''
  \emph{Fuzzy Sets Syst.}, vol.~31, pp. 343--349, 1989.

\bibitem{Ata2020}
K.~T. Atanassov, \emph{Interval-{V}alued {I}ntuitionistic {F}uzzy {S}ets}, ser.
  Studies in Fuzziness and Soft Computing.\hskip 1em plus 0.5em minus
  0.4em\relax Springer, Berlin, Heidelberg, 2020, vol. 388.

\bibitem{DGM2017}
S.~Das, D.~Guha, and R.~Mesiar, ``Extended {B}onferroni mean under
  intuitionistic fuzzy environment based on a strict t-conorm,'' \emph{IEEE
  Trans. Syst., Man, Cybern., Syst.}, vol.~47, no.~8, pp. 2083--2099, 2017.

\bibitem{Xu2007a}
Z.~Xu, ``Some similarity measures of intuitionistic fuzzy sets and their
  applications to multiple attribute decision making,'' \emph{Fuzzy Optim.
  Decis. Making}, vol.~6, no.~2, pp. 109--121, 2007.

\bibitem{CCL2016a}
S.-M. Chen, S.-H. Cheng, and T.-C. Lan, ``Multicriteria decision making based
  on the {TOPSIS} method and similarity measures between intuitionistic fuzzy
  values,'' \emph{Inf. Sci.}, vol. 367, pp. 279--295, 2016.

\bibitem{ZZ2021}
W.-B. Zhang and G.-Y. Zhu, ``A multiobjective optimization of {PCB} prototyping
  assembly with {OFA} based on the similarity of intuitionistic fuzzy sets,''
  \emph{IEEE Trans. Fuzzy Syst.}, vol.~29, no.~7, pp. 2054--2061, 2021.

\bibitem{JLZJHZY2022}
Q.~Jiang, S.~Lee, X.~Zeng, X.~Jin, J.~Hou, W.~Zhou, and S.~Yao, ``A multi-focus
  image fusion scheme based on similarity measure of transformed isosceles
  triangles between intuitionistic fuzzy sets,'' \emph{IEEE Trans. Instrum.
  Meas.}, vol.~71, p. 5013115 (15 pages), 2022, doi:{10.1109/TIM.2022.3169571}.

\bibitem{LC2002}
D.-F. Li and C.~Cheng, ``New similarity measures of intuitionistic fuzzy sets
  and application to pattern recognitions,'' \emph{Pattern Recognit. Lett.},
  vol.~23, pp. 221--225, 2002.

\bibitem{LS2003}
Z.~Liang and P.~Shi, ``Similarity measures on intuitionistic fuzzy sets,''
  \emph{Pattern Recognit. Lett.}, vol.~24, no.~15, pp. 2687--2693, 2003.

\bibitem{HYHL2012}
C.-M. Hwang, M.-S. Yang, W.-L. Hung, and M.-G. Lee, ``A similarity measure of
  intuitionistic fuzzy sets based on the sugeno integral with its application
  to pattern recognition,'' \emph{Inf. Sci.}, vol. 189, pp. 93--109, 2012.

\bibitem{PHK2013}
G.~A. Papakostas, A.~G. Hatzimichailidis, and V.~G. Kaburlasos, ``Distance and
  similarity measures between intuitionistic fuzzy sets: A comparative analysis
  from a pattern recognition point of view,'' \emph{Pattern Recognit. Lett.},
  vol.~34, no.~14, pp. 1609--1622, 2013.

\bibitem{CCL2016}
S.-M. Chen, S.-H. Cheng, and T.-C. Lan, ``A novel similarity measure between
  intuitionistic fuzzy sets based on the centroid points of transformed fuzzy
  numbers with applications to pattern recognition,'' \emph{Inf. Sci.}, vol.
  343, pp. 15--40, 2016.

\bibitem{Ngu2016}
H.~Nguyen, ``A novel similarity/dissimilarity measure for intuitionistic fuzzy
  sets and its application in pattern recognition,'' \emph{Expert Syst. Appl.},
  vol.~45, pp. 97--107, 2016.

\bibitem{JJLY2019}
Q.~Jiang, X.~Jin, S.-J. Lee, and S.~Yao, ``A new similarity/distance measure
  between intuitionistic fuzzy sets based on the transformed isosceles
  triangles and its applications to pattern recognition,'' \emph{Expert Syst.
  Appl.}, vol. 116, pp. 439--453, 2019.

\bibitem{Xiao2021}
F.~Xiao, ``A distance measure for intuitionistic fuzzy sets and its application
  to pattern classification problems,'' \emph{IEEE Trans. Syst., Man, Cybern.,
  Syst.}, vol.~51, no.~6, pp. 3980--3992, 2021.

\bibitem{LZ2018}
M.~Luo and R.~Zhao, ``A distance measure between intuitionistic fuzzy sets and
  its application in medical diagnosis,'' \emph{Artif. Intell. Med.}, vol.~89,
  pp. 34--39, 2018.

\bibitem{BZLW2019}
X.~Bai, Y.~Zhang, H.~Liu, and Y.~Wang, ``Intuitionistic center-free {FCM}
  clustering for {MR} brain image segmentation,'' \emph{IEEE J. Biomed. Health
  Informat.}, vol.~23, no.~5, pp. 2039--2051, 2019.

\bibitem{LSM2018}
Q.~D. Lohani, R.~Solanki, and P.~K. Muhuri, ``Novel adaptive clustering
  algorithms based on a probabilistic similarity measure over {A}tanassov
  intuitionistic fuzzy set,'' \emph{IEEE Trans. Fuzzy Syst.}, vol.~26, no.~6,
  pp. 3715--3729, 2018.

\bibitem{LLC2014}
N.~Li, F.~Liu, and Z.~Chen, ``A texture measure defined over intuitionistic
  fuzzy set theory for the detection of built-up areas in high-resolution sar
  images,'' \emph{IEEE J. Sel. Top. Appl. Earth Obs. Remote Sens.}, vol.~7,
  no.~10, pp. 4255--4265, 2014.

\bibitem{Mit2003}
H.~B. Mitchell, ``On the {D}engfeng--{C}huntian similarity measure and its
  application to pattern recognition,'' \emph{Pattern Recognit. Lett.},
  vol.~24, no.~16, pp. 3101--3104, 2003.

\bibitem{Sz2014}
E.~Szmidt, \emph{Distances and {S}imilarities in {I}ntuitionistic {F}uzzy
  {S}ets}, ser. Studies in Fuzziness and Soft Computing.\hskip 1em plus 0.5em
  minus 0.4em\relax Springer, Berlin, Heidelberg, 2014, vol. 307.

\bibitem{Li2014}
D.-F. Li, \emph{Decision and {G}ame {T}heory in {M}anagement with
  {I}ntuitionistic {F}uzzy {S}ets}, ser. Studies in Fuzziness and Soft
  Computing.\hskip 1em plus 0.5em minus 0.4em\relax {S}pringer-{V}erlag
  {B}erlin {H}eidelberg, 2014, vol. 308.

\bibitem{WCZ2022}
X.~Wu, Z.~Zhu, C.~Chen, G.~Chen, and P.~Liu, ``A monotonous intuitionistic
  fuzzy {TOPSIS} method under general linear orders via admissible distance
  measures,'' \emph{IEEE Trans. Fuzzy Syst.}, p. 14 pages, 2022,
  doi:10.1109/TFUZZ.2022.3205435.

\bibitem{SK2000}
E.~Szmidt and J.~Kacprzyk, ``Distances between intuitionistic fuzzy sets,''
  \emph{Fuzzy Sets Syst.}, vol. 114, no.~3, pp. 505--518, 2000.

\bibitem{Gr2004}
P.~Grzegorzewski, ``Distances between intuitionistic fuzzy sets and/or
  interval-valued fuzzy sets based on the {H}ausdorff metric,'' \emph{Fuzzy
  Sets Syst.}, vol. 148, no.~2, pp. 319--328, 2004.

\bibitem{HY2004}
W.-L. Hung and M.-S. Yang, ``Similarity measures of intuitionistic fuzzy sets
  based on {H}ausdorff distance,'' \emph{Pattern Recognit. Lett.}, vol.~25,
  no.~14, pp. 1603--1611, 2004.

\bibitem{WX2005}
W.~Wang and X.~Xin, ``Distance measure between intuitionistic fuzzy sets,''
  \emph{Pattern Recognit. Lett.}, vol.~26, no.~13, pp. 2063--2069, 2005.

\bibitem{XC2008}
Z.~Xu and J.~Chen, ``An overview of distance and similarity measures of
  intuitionistic fuzzy sets,'' \emph{Int. J. Uncertain. Fuzziness. Knowl.-Based
  Syst.}, vol.~16, no.~04, pp. 529--555, 2008.

\bibitem{HY1981}
C.-L. Hwang and K.~Yoon, \emph{Multiple Attribute Decision Making: Methods and
  Applications A State-of-the-Art Survey}, ser. Lecture Notes in Economics and
  Mathematical Systems.\hskip 1em plus 0.5em minus 0.4em\relax
  {S}pringer-{V}erlag {B}erlin {H}eidelberg, 1981, vol. 186.

\bibitem{YC2009}
Y.~Yang and F.~Chiclana, ``Intuitionistic fuzzy sets: spherical representation
  and distances,'' \emph{Int. J. Intell. Syst.}, vol.~24, no.~4, pp. 399--420,
  2009.

\bibitem{HY2008}
W.-L. Hung and M.-S. Yang, ``On the {$J$}-divergence of intuitionistic fuzzy
  sets with its application to pattern recognition,'' \emph{Inf. Sci.}, vol.
  178, no.~6, pp. 1641--1650, 2008.

\bibitem{JK2018}
R.~Joshi and S.~Kumar, ``A dissimilarity {J}ensen--{S}hannon divergence measure
  for intuitionistic fuzzy sets,'' \emph{Int. J. Intell. Syst.}, vol.~33,
  no.~11, pp. 2216--2235, 2018.

\bibitem{HPK2012}
A.~G. Hatzimichailidis, G.~A. Papakostas, and V.~G. Kaburlasos, ``A novel
  distance measure of intuitionistic fuzzy sets and its application to pattern
  recognition problems,'' \emph{Int. J. Intell. Syst.}, vol.~27, no.~4, pp.
  396--409, 2012.

\bibitem{YC2012}
Y.~Yang and F.~Chiclana, ``Consistency of {2D} and {3D} distances of
  intuitionistic fuzzy sets,'' \emph{Expert Syst. Appl.}, vol.~39, no.~10, pp.
  8665--8670, 2012.

\bibitem{SMLZC2018}
F.~Shen, X.~Ma, Z.~Li, Z.~Xu, and D.~Cai, ``An extended intuitionistic fuzzy
  {TOPSIS} method based on a new distance measure with an application to credit
  risk evaluation,'' \emph{Inf. Sci.}, vol. 428, pp. 105--119, 2018.

\bibitem{SWQH2019}
Y.~Song, X.~Wang, W.~Quan, and W.~Huang, ``A new approach to construct
  similarity measure for intuitionistic fuzzy sets,'' \emph{Soft Comput.},
  vol.~23, no.~6, pp. 1985--1998, 2019.

\bibitem{XC2012}
Z.~Xu and X.~Cai, \emph{{I}ntuitionistic {F}uzzy {I}nformation {A}ggregation:
  {T}heory and {A}pplications}, ser. Mathematics Monograph Series.\hskip 1em
  plus 0.5em minus 0.4em\relax {S}cience {P}ress, 2012, vol.~20.

\bibitem{Xu2007}
Z.~Xu, ``Intuitionistic fuzzy aggregation operators,'' \emph{IEEE Trans. Fuzzy
  Syst.}, vol.~15, no.~6, pp. 1179--1187, 2007.

\bibitem{SK2001}
E.~Szmidt and J.~Kacprzyk, ``Entropy for intuitionistic fuzzy sets,''
  \emph{Fuzzy Sets Syst.}, vol. 118, no.~3, pp. 467--477, 2001.

\bibitem{ES2003}
D.~M. Endres and J.~E. Schindelin, ``A new metric for probability
  distributions,'' \emph{IEEE Trans. Inf. Theory}, vol.~49, no.~7, pp.
  1858--1860, 2003.

\end{thebibliography}


\end{document}